\documentclass[12pt]{article}
\usepackage{amssymb}
\usepackage{amsmath}
\usepackage{amsfonts}
\usepackage{mitpress}

\newtheorem{theorem}{Theorem}

\newtheorem{definition}[theorem]{Definition}

\newtheorem{lemma}[theorem]{Lemma}

\newtheorem{proposition}[theorem]{Proposition}

\newdimen\dummy
\dummy=\oddsidemargin
\input{tcilatex}
\begin{document}

\title{$p-$adic Abel-Jacobi map and $p$-adic Gross-Zagier formula for
Hilbert modular forms}
\author{Ivan Blanco-Chacon and Ignacio Sols \\
University College Dublin and Universidad Complutense}
\maketitle

\begin{abstract}
We compute the $p-$adic Abel-Jacobi map of the product of a Hilbert modular
surface and a modular curve at a null-homologous (modified) embedding of the
modular curve in this product, evaluated on differentials associated to a
Hilbert cuspidal form $f$ of weight $(2,2)$ and a cuspidal form of weight $2$%
. We generalize this computation to suitable null-homological cycles in the
fibre products of the universal families on the surface and the curve,
evaluated at differentials associated to $f$ and $g$ of higher weights.

We express the values of the \ $p-$adic Abel-Jacobi map at these weights in
terms of a $p$-adic $L-$ function associated to a Hida family of Hilbert
modular forms and a Hida family of cuspidal forms. Our\ function is \ a
Hilbert modular analogue of the $p-$ adic $L-$function introduced in \cite%
{DR}
\end{abstract}

\section{\textbf{Introduction}}

The present work starts a series devoted to studying null-homologous cycles
on the product of a modular curve and a Hilbert modular surface, and higher
dimensional generalizations of this, with a view to producing distinguished
cohomology classes in the Bloch-Kato Selmer group of the representation
associated to an elliptic curve $E/\mathbb{Q}$ twisted by an Artin
representation associated to a Hilbert modular form, and to relate these
classes with the dimension of the Bloch-Kato Selmer group. Following the
philosophy of \cite{DR}, of which our program is a prolongation to the
Hilbert modular setting, our aim is to obtain a $p-$adic Gross-Zagier type
formula relating the image under the $p-$adic Abel-Jacobi map of certain
null-homologous cycles with the special values of a $p-$adic $L-$function
arising out of a $p-$adic family of Hilbert cuspidal forms passing by a
given one $f$ of weight $(2,2)$ and a $p-$adic family of cuspidal forms
passing by a rational one $g$ of weight $2$, for a prime $p$, generic in the
sense we will specify.

To state the main result, we need some preliminary notions and notations.
Let $\mathfrak{K}$ be a real quadratic field and fix an inclusion $inc:%
\mathfrak{K\hookrightarrow \mathbb{R}\ }$, and assume, just for convenience,
that $\mathfrak{K}$ has a unit of norm $-1$, so that $\mathfrak{u}^{+}/(%
\mathfrak{u}^{+})^{2}$ is trivial (otherwise, the content of this article
holds after replacing our cohomologies by their invariant part under the
natural action of this two-elements group on the moduli\ schemes we
consider). Let $\mathfrak{o}$ be the ring of integers of $\mathfrak{K}$ and
choose a narrow class of it and $\mathfrak{a\subseteq o}$ a primitive
integral ideal representing it. Assume that $N_{2}\geq 4$, the discriminant $%
d_{\mathfrak{K}}=Nm(\mathfrak{d}_{\mathfrak{K}}\mathfrak{)}$, and $N_{1}=Nm(%
\mathfrak{a)}$ are mutually coprime, and write $N=d_{\mathfrak{K}}N_{1}N_{2}$
. Consider a ring $R\supseteq \mathfrak{o}$ in which $N$ is invertible.
Denote by $Y_{R}(N_{2})$ the Hilbert modular surface of level $N_{2}$,
chosen toroidal desigularization of the Satake minimal compactification $%
\overline{\mathcal{M}_{R}(\mathfrak{a},N_{2})}$ of the geometrically normal
and irreducible surface $\mathcal{M}_{R}(\mathfrak{a},N_{2})$ over $R$ which
represents the functor assigning to schemes $S$ over $R$ the set of abelian
surfaces over $R$ with $\mathfrak{o-}$real multiplication, ordered $%
\mathfrak{a-}$polarization and $N_{2}-$level structure (see the definitions
in section 2). The Rapoport rank one $\mathfrak{o\otimes }_{\mathbb{Z}}%
\mathcal{O}_{\mathcal{M}_{R}(N_{2})}$ -sheaf $\mathcal{R}_{\mathcal{M}%
_{R}(N_{2})}:=pr_{\mathcal{M}_{R}(N_{2})\ast }\Omega _{A_{\mathcal{\mathcal{M%
}}_{R}(N_{2})}^{U}/\mathcal{\mathcal{M}}_{R}(N_{2})}$ of relative
differentials of the universal object $A_{\mathcal{\mathcal{M}}%
_{R}(N_{2})}^{U}$ is locally free, because $d_{\mathfrak{K}}$ is invertible
in $R$. The inclusion and the conjugate inclusion of $\mathfrak{o}$ in $R$
associate to $\mathcal{R}_{\mathcal{M}_{R}(N_{2})}$ two line bundles whose $%
k,k^{\prime }$ powers are the modular line bundles $L_{\mathcal{M}%
_{R}(N_{2})}^{(k,k^{\prime })},$ whose sections are the $\mathfrak{a}$%
-Hilbert modular forms of weight $(k,k^{\prime })$ over $R$ . These extend
to line bundles $L_{Y_{R}(N_{2})}^{(k,k^{\prime })}$ on the whole $%
Y_{R}(N_{2})$ because the Rapoport sheaf $\mathcal{R}_{\mathcal{M}%
_{R}(N_{2})}$ extends to a rank one locally free $\mathfrak{o\otimes }_{%
\mathbb{Z}}\mathcal{O}_{Y_{R}(N_{2})}$ -sheaf $\mathcal{R}_{Y_{R}(N_{2})}$.
Hilbert modular forms of parallel weight $(k,k^{\prime })$ can be defined
equivalently as sections of these line bundles on $\mathcal{M}_{R}(\mathfrak{%
a},N_{2})$ or on the whole $Y_{R}(N_{2})$ (Koecher principle), and those
vanishing on the transform divisor $D^{c}$ of the cusps $\overline{\mathcal{M%
}_{R}(N_{2})}-\mathcal{M}_{R}(N_{2})$ are called cuspidal. The complex
manifold $\mathcal{M}_{\mathbb{C}}(\mathfrak{a},N_{2})(\mathbb{C})$ is the
quotient of the Poincar\'{e} square half-plane $\mathfrak{H}^{2}$ by the
Hilbert modular group $\Gamma _{1}(\mathfrak{a,}N_{2})$ consisting of
matrices of $\left( 
\begin{array}{cc}
\mathfrak{o} & \mathfrak{a^{\ast }} \\ 
\mathfrak{ad}_{\mathfrak{K}}N_{2} & \mathfrak{o}%
\end{array}%
\right) $ whose determinant belongs to $\mathfrak{u}^{+}.$ Since we assume $%
\mathfrak{u}^{+}/(\mathfrak{u}^{+})^{2}$ trivial, it is also the quotient by
the action of the subgroup $\Gamma _{1}^{1}(\mathfrak{a,}N_{2})\subseteq
\Gamma _{1}(\mathfrak{a,}N_{2})$ of unimodular matrices.

\bigskip The bundle of logarithmic differentials $\Omega
_{Y_{R}(N_{2})/R}(\log D^{c})\approx L_{Y_{R}(N_{2})}^{(2,0)}\oplus
L_{Y_{R}(N_{2})}^{(0,2)}$ has determinant $\omega
_{Y_{R}(N_{2})}(D^{c})\approx L_{Y_{R}(N_{2})}^{(2,2)}$ so that a cuspidal $%
\mathfrak{a}$-Hilbert modular form $f$ of weight $(2,2)$ corresponds
bijectively to a 2-differential $\omega _{f}$ on $Y_{R}(N_{2})$ (also
denoted $\omega (f)$ if the expression of $f$ is long), this generalizing to
an $L_{Y_{R}(N_{2})}^{(k,k^{\prime })}$-valued $\log D^{c}-$differential $%
\omega _{f}$ if $f$ has weight $(k,k^{\prime })\geq (2,2)$.

\bigskip\ There is a generically injective embedding $j_{X_{R}(N)}:X_{R}(N)%
\longrightarrow Y_{R}(N_{2})$ of the modular curve $X_{R}(N)$ into $%
Y_{R}(N_{2})$ for the congruence group $\Gamma (N)=\Gamma _{0}(d_{\mathfrak{K%
}}N_{1})\cap \Gamma _{1}(N_{2})$ after compactification, desingularization
and base change to $Spec$ $R$, and the restriction $j_{X_{R}(N)}^{\ast
}L_{Y_{R}(N_{2})}^{(k,k^{\prime })}$ is $L_{0,R}^{k+k^{\prime }},$ for the
modular line bundle $L_{0}$ on $X(N)$ (denoted \underline{$\omega $} in the
literature) such that $L_{0}^{2}$ $\approx \omega _{X(N)}(cusps)$.\ 

\bigskip Let $K\supseteq $ $\mathfrak{K}$ be a number field having ring of
integers $R$ $\supseteq $ $\mathfrak{o}$. Choose $p$\ splitting in $K$ so
that it splits in $\mathfrak{K}$ as $p=\pi \pi ^{\prime }$. Then\textit{\ }$%
p $ does not divide the discriminant $d_{\mathfrak{K}}$\textit{\ }, and we
choose it prime with\textit{\ }$N=d_{\mathfrak{K}}N_{1}N_{2}$ . The
embedding $\mathfrak{K\hookrightarrow \mathfrak{K}}_{\pi }\approx \mathbb{Q}%
_{p}$ extends to an embedding $K\hookrightarrow \overline{\mathbb{Q}}_{p}$
which restricts to $R_{\mathfrak{m}}\approx \mathbb{Z}_{p}$ for a maximal
ideal $\mathfrak{m}$ of $R$ lying on $\pi $, so that $\mathbb{Z}\subseteq 
\mathfrak{o\subseteq }R$ induces isomorphisms $\mathbb{Z}_{p}\approx 
\mathfrak{o}_{\pi }\approx R_{\mathfrak{m}}$ and $\mathbb{Q}_{p}\approx 
\mathfrak{K}_{\pi }\approx K_{\mathfrak{m}}$;\ and we analogously obtain $%
\mathbb{Z}_{p}\approx \mathfrak{o}_{\pi ^{\prime }}\approx R_{\mathfrak{m}%
^{\prime }}$ and $\mathbb{Q}_{p}\approx \mathfrak{K}_{\pi ^{\prime }}\approx
K_{\mathfrak{m}^{\prime }}$ . The locus of points in $\mathcal{M}_{\mathbb{Q}%
_{p}}(N_{2})$ corresponding to Abelian surfaces with ordinary $\mathcal{%
\mathbb{F}}_{p}-$ specialization is $\mathcal{A\cap M}_{\mathbb{Z}%
_{p}}(N_{2})$ for the mod. $p-$nonvanishing locus $\mathcal{A}$ of the Hasse
section of the Hasse modular line bundle on $Y_{\mathbb{Q}_{p}}(N_{2})$.
This is an open set of the associated $\mathbb{Q}_{p}-$rigid space and has
as base of strict neighborhoods the mod. $p^{n}$-nonvanishing loci $\mathbb{%
\mathcal{W}}_{p^{-n}}$ of the Hasse section. The complement of the Hasse
divisor $D^{h}$ where this section vanishes is $Y_{\mathbb{Q}_{p}}^{\prime
}(N_{2})=\dbigcup\limits_{\varepsilon >0}\mathbb{\mathcal{W}}_{\varepsilon }$%
, a Zariski open set which generalizes to the Hilbert-modular context the
open set $X_{\mathbb{Q}_{p}}^{\prime }(N)=j_{X_{R}(N)}^{-1}(Y_{\mathbb{Q}%
_{p}}^{\prime }(N_{2}))$. The $p-$adic Hilbert modular forms defined as
rigid sections $f$ of the restricted modular line bundles $L_{\mathcal{A}%
}^{(k,k^{\prime })}$ are also $p-$adic Hilbert modular forms in the sense of
Katz in \textbf{\cite{Split Katz}} so they have a $q-$expansion at the
nonramified cusps, which at the one called "standard cusp" is denoted $%
f(q):=\dsum\limits_{\nu \in \mathfrak{(a}^{-1})^{+}\cup \{0\}}^{{}}a_{\nu
}q^{\nu }$ , with $a_{\nu }\in \mathbb{Z}_{p}$ and $a_{0}$ null when $f$ is
cuspidal, being $f\mapsto f(q)$ an injective map;\ those extending to some $%
\mathbb{\mathcal{W}}_{\varepsilon }$ are called overconvergent. The
Frobenius morphism $\phi $ on the ordinary locus $Y_{\mathbb{F}_{p}}^{\prime
}(N_{2})\subseteq Y_{\mathbb{F}_{p}}(N_{2})$ has a natural cover $\phi :$ $%
\phi ^{\ast }L_{Y_{\mathbb{F}_{p}}^{\prime }(N_{2})}^{(k,k^{\prime
})}\longrightarrow L_{Y_{\mathbb{F}_{p}}^{\prime }(N_{2})}^{(k,k^{\prime })}$%
, and a lifting $\phi :\mathcal{A\longrightarrow \mathcal{A}}$ with natural
cover$\ \phi :\phi ^{\ast }L_{\mathcal{A}}^{(k,k^{\prime })}\longrightarrow
L_{\mathcal{A}}^{(k,k^{\prime })}$ so it acts on $p-$adic $\mathfrak{a-}$%
Hilbert modular forms $f$, being $\dsum\limits_{\nu \in \mathfrak{(a}%
^{-1})^{+}\cup \{0\}}^{{}}a_{\nu }q^{p\nu }$ the $q-$expansion of the
transform; by a recent result of Conrad \cite{Conrad}, this lifting extends
to a morphism between open sets $\mathbb{\mathcal{W}}_{\varepsilon }$, so
that $\phi $ preserves overconvergence.

Let $o$ be a rational point of $X_{\mathbb{Q}}(N)$ and let $X_{\mathbb{Q}%
}(N)_{o}$ be the modified diagonal divisor 
\begin{equation*}
X_{\mathbb{Q}}(N)-o\times X_{\mathbb{Q}}(N)-X_{\mathbb{Q}}(N)\times o
\end{equation*}%
whose complexification $X_{\mathbb{C}}(N)_{o}$ is null-homologous in $X_{%
\mathbb{C}}(N)\times Y_{\mathbb{C}}(N_{2})$ as a consequence of the fact
that $Y_{\mathbb{C}}$ has null $1-$homology (we sometimes abuse notation by
identifying complex schemes with the analytic space of their complex
points). Let $f$\ be an $\mathfrak{a-}$Hilbert cuspidal form of parallel
weight $(2,2)$ and level $N_{2}$ over the integer ring\textbf{\ }$R\supseteq 
\mathfrak{o}$ of a number field $K\supseteq \mathfrak{K}$ . Assume
furthermore that $f$ is an eigenform for the Hecke operators as defined for
instance in \textbf{\cite{VdG} }VI 1.

Assume we are given a rational cuspidal form $g$ of weight $2$ for the
modular group $\Gamma (N)$ , normalized eigenform for the Hecke operators $%
T_{n}$ and assume that the prime $p$ is chosen ordinary for $g$ and that $f$
is non-ordinary at $\pi $ and $\pi ^{\prime },$ i.e. it has Hecke eigenvalue 
$a_{\pi }$ being non-unit of $R_{\mathfrak{m}}\approx \mathfrak{o}_{\pi
}\approx \mathbb{Z}_{p}$ and $a_{\pi ^{\prime }}$ being non-unit of $R_{%
\mathfrak{m}^{\prime }}\approx \mathfrak{o}_{\pi ^{\prime }}\approx \mathbb{Z%
}_{p}$ , i.e. having positive "slope" $\sigma =$ $ord_{\mathfrak{m}}(a_{\pi
})$ and $\sigma ^{\prime }=ord_{\mathfrak{m}^{\prime }}(a_{\pi ^{\prime }})$
.

\bigskip \textbf{\ }Take a model of the above curve, surface and
null-homologous cycle over $\mathbb{Q}_{p}$, as smooth and projective over $%
\mathbb{Z}_{p}$ (and fix an inclusion $\mathbb{Q}_{p}\hookrightarrow \mathbb{%
C}$ once and for all\textbf{).} We can use then the Besser theory developed
in\textbf{\ }\cite{Besser inventi}, to compute in section 3 $\ $the $p-$adic
Abel-Jacobi map of $X_{\mathbb{Q}_{p}}(N)\times Y_{\mathbb{Q}_{p}}(N_{2})$
at the null-homologous cycle $X_{\mathbb{Q}_{p}}(N)_{o}$ evaluated on the de
Rham class of $\omega _{f}\otimes \eta _{g}^{u-r}$ for the differential $%
\omega _{f}$ on $Y_{\mathbb{Q}_{p}}(N_{2})$ attached to $f$ and the unique
"unit root" differential $\eta _{g}^{u-r}$ on $X_{\mathbb{Q}_{p}}(N)$ having
cup product $<\omega _{g},\eta _{g}^{u-r}>=1.$ Denote, just as in \cite{DR} $%
e_{ord}$ the projection to the "ordinary subspace" of the de Rham cohomology
space $H_{dR}^{1}$ $(X_{\mathbb{Q}_{p}}^{\prime }(N)/\mathbb{Q}_{p})$ where $%
\phi $ acts with eigenvalue $p$ , and $e_{g}$ the projection to the
2-dimensional isotypic component of $\omega _{g}$ , this lying in the
subspace $H_{dR}^{1}(X_{\mathbb{Q}_{p}}(N)/\mathbb{Q}_{p})$. The main result
in section 3 is 
\begin{equation}
AJ_{p}(X_{\mathbb{Q}_{p}}(N)_{o})(\omega _{f}\otimes \eta
_{g}^{u-r})=<e_{g}e_{ord}Q_{f}(\phi )^{-1}j_{X_{\mathbb{Q}_{p}}^{\prime
}(N)}^{\ast }\varrho ^{\sharp },\eta _{g}^{u-r}>
\label{Abel-Jacobi introduccion}
\end{equation}%
where $Q_{f}(x)$\ is the characteristic polynomial of $\phi $ acting on the
isotypic component $H^{2}(Y_{\mathbb{Q}_{p}}(N_{2})/\mathbb{Q}_{p})(f)$ of $%
\omega _{f}$ in $H^{2}(Y_{\mathbb{Q}_{p}}(N_{2})/\mathbb{Q}_{p})$, and $%
\varrho ^{\sharp }$ is a primitive of the overconvergent $2-$differential $%
Q_{f}(\phi )\omega _{f\text{ }}$ whose existence we prove in section 3.
Under the above non-ordinary conditions, there is in fact (as recalled in
sec. 4) by \cite{David et al.},\textbf{\ }a unique primitive $\varrho
^{\sharp }$ of type $(1,0)$, i.e. overconvergent section of the first direct
summand of the decomposition $\Omega _{Y_{R}(N_{2})/R}(\log D^{c})\approx
L_{Y_{R}(N_{2})}^{(2,0)}\oplus L_{Y_{R}(N_{2})}^{(0,2)}$ .

We give the proof of (\ref{Abel-Jacobi introduccion}), more in general, for
a $1-$differential $\eta $ not necessarily `of the form $\eta _{g}^{u-r}$,
but then $e_{g}e_{ord}$ does not appear in the first factor of the the
product (\ref{Abel-Jacobi introduccion}), so we force this factor to stay in
the kernel $H_{dR}^{1}(X_{\mathbb{Q}_{p}}(N)/\mathbb{Q}_{p})$ of the residue
map of $H_{dR}^{1}$ $(X_{\mathbb{Q}_{p}}^{\prime }(N)/\mathbb{Q}_{p})$ by
replacing $Q_{f}(\phi )$ by $P_{f}(\phi )=(1-p^{-1}\phi )Q_{f}(\phi )$ , as
the first factor kills the target of the residue map and the product (\ref%
{Abel-Jacobi introduccion}) remains the same while replacing $Q_{f}(\phi )$
by any multiple polynomial. \bigskip

In section 4, we generalize (\ref{Abel-Jacobi introduccion}) for the
evaluation of a $p-$ adic Abel-Jacobi at Hilbert modular forms $f$ and $g$
of higher weights $(k,k)=(2+n,2+n)$ and $k_{0}=2+n_{0}$ for $0\mathfrak{\leq 
}n_{0}=2n-2t$, with $t>0$ or $t=n=n_{0}=0$ at certain null-homologous cycles 
$\Delta _{n,n_{0}}$ of higher dimension. Choose a rational prime $p$ just as
above, i.e. not dividing $N$ and splitting in $K$ and in $\mathfrak{o}$ as $%
p=\pi \pi ^{\prime }$ being ordinary for $g$ and nonordinary, i.e. of
positive slope $ord_{\mathfrak{m}}(a_{\pi })$ and $ord_{\mathfrak{m}^{\prime
}}(a_{\pi ^{\prime }})$, for $f$. \ Consider the $L^{(n,n)}$ -valued $2-$%
differential $\omega _{f}$ of $Y(N_{2})$ and the $L_{0}^{n_{0}}$-valued $1-$%
differential $\omega _{g}$ of $X(N)$ (the definition is over $\mathbb{Q}_{p}$%
, when omitted). The product analogous to (\ref{Abel-Jacobi introduccion})
would need a covariant derivative for differentials valued in modular line
bundles, but we have at least the Gauss-Manin covariant $\log $ $D^{c}$%
-derivative $\nabla ^{GM}$ on a bundle $\mathcal{L}^{n}$ in which $\mathcal{R%
}^{n}=\dbigoplus\limits_{k,k^{\prime }\geq 0,\text{ }k+k^{\prime
}=n}^{{}}L^{(k,k^{\prime })}$ is included with a splitting projection $sp$
on just $\mathcal{A}$ . This $\mathcal{L}^{n}$ is the $n-$symmetric power of
the four -bundle $\mathcal{L}$ of relative first de Rham cohomology which
generalizes to the Hilbert modular context the direct factor $\mathcal{L}_{0}
$ of its restriction $j_{X(N)}^{\ast }\mathcal{L=\mathcal{L}}_{0}\oplus 
\mathcal{L}_{0}$ , a decomposition which follows from the fact that the
universal semiabelian surface splits on $X(N)$ as fibre product of two
copies of the universal generalized elliptic curve, because of our
assumptions. Under our non-ordinary hypothesis, by \cite{David et al.}, for
the overconvergent "$\pi ^{\prime }-$depletion" $f^{[\pi ^{\prime }]}$ of $f$
of $q-$expansion $f^{[\pi ^{\prime }]}(q)=\dsum_{\nu \in \mathfrak{(a}%
^{-1})^{+},\text{ }\pi ^{\prime }\nmid \nu }a_{\nu }q^{\nu }$, the
associated $L^{(n,n)}$-valued 2-differential $\omega _{f^{[\pi ^{\prime }]}}$
has a unique overconvergent $\mathcal{L}^{2n}-$valued differential $\varrho
_{\mathcal{F}^{[\pi ^{\prime }]}}$ of type $(1,0)$ such that $\nabla
^{GM}\varrho _{\mathcal{F}^{[\pi ^{\prime }]}}$ $=\omega _{f^{[\pi ^{\prime
}]}}$. The split projection in $\mathcal{A}$ of $\varrho _{\mathcal{F}^{[\pi
^{\prime }]}}$ provides a unique nearly overconvergent form $%
F_{[2+n,n]}^{[\pi ^{\prime }]}$ of weight $(2+n,n)$ whose associated $%
L^{(n,n)}$- differential $\varrho _{F_{[2+n,n]}^{[\pi ^{\prime }]}}$has $%
\omega _{f^{[\pi ^{\prime }]}}$ as its "split derivative" $\nabla
_{sp}^{GM}:=sp\circ \nabla ^{GM}$ . In terms of the $\theta ,\theta ^{\prime
}$ operators on $p-$adic forms acting as $\dsum \nu a_{\nu }q^{\nu }$ and $%
\dsum \nu ^{\prime }a_{\nu }q^{\nu }$ on $q-$expansions, this amounts to $%
F_{[2+n,n]}^{[\pi ^{\prime }]}=\theta ^{\prime -1}f^{[\pi ^{\prime }]}$ ,
i.e. it has $q-$expansion $\dsum_{\nu \in \mathfrak{(a}^{-1})^{+},\text{ }%
\pi ^{\prime }\nmid \nu }\nu ^{\prime -1}a_{\nu }q^{\nu }$ with coefficients
still lying in $R_{\mathfrak{m}}$ as the $\nu ^{\prime }$ inverted in them
are units of $\mathfrak{o}_{\pi }\approx R_{\mathfrak{m}}$ ; and the
analogous holds for $\pi -$depletion.

We can express the explicit value of the product (\ref{Abel-Jacobi
introduccion}) by a purely formal computation which we perform in section 5.
For this, we decompose $f$ as sum $f_{00^{\prime }}+f_{01^{\prime
}}+f_{10^{\prime }}+f_{11^{\prime }}$ of its four stabilizations or
eigenvectors of the Verschiebung operator $U_{p}$ on $p-$adic $\mathfrak{a-}$%
Hilbert modular forms expressed on $q-$expansions by $(U_{p}f)(q)=\dsum_{\nu
\in \mathfrak{(a}^{-1})^{+},}a_{p\nu }q^{\nu }$ and preserving
overconvergence. Each of the four root factors of the four dimensional
polynomial $Q_{f}(\phi )$ acts on one of these four summands as a $p-$%
depletion operator, and the other three factors become in the cup-product (%
\ref{Abel-Jacobi introduccion}) (up to summands with vanishing contribution
to the final product) the multiplication by a scalar factor, which is the
cause of the Euler factors $\mathcal{E}(f,g)$ and $\mathcal{E}_{1}(g)$
appearing in the following theorem computing suitable $p-$adic Abel-Jacobi
maps. Denoting $\mathcal{E}^{n_{0}}$ the fibred power of the universal
generalized elliptic curve $\mathcal{E}$ over $X(N)$ , and $\widehat{A}^{U,n}
$ the fibred power of the smooth compactification $\widehat{A}^{U}$ of the
universal semiabelian scheme over $Y(N_{2})$ , we construct in section 4
null-homologous cycles $\Delta _{n,n_{0}}$ in $\widehat{A}^{U,n}\times 
\mathcal{E}^{n_{0}}$ for integers $n,n_{0}$ as above. These generalize, with
the help of the Scholl operators killing the nonintermediate homology of a
product of generalized elliptic curves, the null-homologous cycle $\Delta
_{0,0}$ which is the modified diagonal $X(N)_{o}$ viewed as a cycle in $%
Y(N_{2})\times X(N)$ .

\bigskip

\begin{theorem}
\label{TEOREMA PRINCIPAL INTRODUCCION} With $\mathfrak{K},\mathfrak{o},R,%
\mathfrak{o,\mathfrak{a}}$ , and integers $N=d_{\mathfrak{K}}$ $N_{1}N_{2}$
as previously introduced, let $f$ be an $\mathfrak{a-}$Hilbert cuspidal form
of level $N_{2}$ and weight $(k,k)=(n+2,n+2)\geq (2,2)$ over $R\supseteq 
\mathfrak{o}$, and let $f(q)=\dsum_{\nu \in \mathfrak{(a}^{-1})^{+}}a_{\nu
}q^{\nu }$ be its $q-$expansion at the standard cusp. Assume we are given an 
$\mathfrak{a-}$Hilbert cuspidal Hecke eigenform $f$ of level $N_{2}\geq 4$
and weight $(k,k)=(n+2,n+2)\geq (2,2)$ over the ring of integers $R\supseteq 
\mathfrak{o}$ of a number field $K$ extending $\mathfrak{K}$, and let $%
f(q)=\dsum_{\nu \in \mathfrak{(a}^{-1})^{+}}a_{\nu }q^{\nu }$ be its $q-$%
expansion at the standard cusp; suppose also that we are given a rational
cuspidal Hecke eigenform $g$ for the modular group $\Gamma (N)=\Gamma
_{0}(d_{\mathfrak{K}}N_{1})\cap \Gamma _{1}(N_{2})$ of weight $%
k_{0}=n_{0}+2\geq 2$ with $n_{0}=2n-2t$ \ being $t>0$ or $t=$ $n_{0}=n=0$ ,
and let $g(q)=\dsum_{n\geq 1}$ $b_{n}q^{n}$ be its $q-$expansion at the
infinity cusp. Denote $\omega _{f}$ the attached $L^{(n,n)}$-differential of
order 2 on a chosen toroidal desingularization $Y_{K}(N_{2})$ of the Satake
compactification of the fine moduli $\mathcal{M}_{K}(\mathfrak{a,}N_{2})$ .
Denote $\omega _{g}\ $the attached$\ \ L_{0}^{n_{0}}$-differential on the
modular curve $X_{\mathbb{Q}}(N)=X_{\mathbb{Q}}(\Gamma (N))$ and $\eta
_{g}^{u-r}$ \ the unique class of $\mathcal{L}_{0}^{n}-$valued de Rham
cohomology such that $<\omega _{g},\eta _{g}^{u-r}>=1$.

Let $p$\ be a prime not dividing $N$ and splitting in $K$\ , and in $%
\mathfrak{K}$ as $p=\pi \pi ^{\prime }$ , such that $p$ is ordinary for $g$
and nonordinary for $f$ , i.e. of positive slope $\sigma ,\sigma ^{\prime }$%
, at $\pi $ and $\pi ^{\prime }$ . Writting $j_{X_{\mathbb{Q}_{p}}^{\prime
}(N)}:X_{\mathbb{Q}_{p}}^{\prime }(N)\hookrightarrow X_{\mathbb{Q}%
_{p}}(N)\longrightarrow Y_{R_{\mathfrak{m}}}(N_{2})$ the complement of the
locus of the points of the moduli with non-ordinary specialization to $X_{%
\mathbb{F}_{p}}(N)$, the value of the $p-$adic Abel-Jacobi map at the
null-homologous cycle $\Delta _{n,n_{0}}$ is 
\begin{eqnarray}
AJ_{p}(\Delta _{n,n_{0}})(\omega _{f}\otimes \eta _{g}^{u-r}) &=&
\label{AJ 2} \\
(-1)^{t}t!\frac{\mathcal{E}_{1}(g)}{\mathcal{E}(f,g)} &<&e_{g}e_{ord}\text{ }%
j_{X^{\prime }(N)}^{\ast }\varrho (\frac{\theta ^{\prime -1-t}f^{[\pi
^{\prime }]}}{2}-\frac{\theta ^{-1-t}f^{[\pi ]}}{2}),\eta _{g}^{u-r}>\text{ ,%
}  \notag
\end{eqnarray}%
\qquad with the non-zero factors 
\begin{equation*}
\mathcal{E}_{1}(g)=1-\beta _{1}^{2}p^{-k_{0}}\text{ and \ }\mathcal{E}%
(f,g)=\dprod\limits_{i,i^{\prime }}(1-\alpha _{i}\alpha _{i^{\prime }}\beta
_{1}p^{2-2k+t})
\end{equation*}%
where $\alpha _{0},\alpha _{1}$ and $\alpha _{0^{\prime }},\alpha
_{1^{\prime }\text{ }}$ and $\beta _{1}$ are the roots of the polynomials $%
1-a_{\pi }x+p^{k-1}x^{2}$ and $1-a_{\pi ^{\prime }}x+p^{k-1}x^{2}$ ;\ and $%
\beta _{1}$ is the root of the polynomial $1-b_{p}x+p^{k_{0}-1}x^{2}$ having
having $p$ - valuation $k_{0}-1$ .
\end{theorem}

\bigskip

\bigskip

In section 6, in theorem \ref{Main theorem copy(1)}, we prove a Gross-Zagier
type formula. For given $f$ and $g$ of weights $(2,2)$ and $2$ , and choice
of splitting prime $p=\pi \pi ^{\prime }$ of positive slope $\sigma ,\sigma
^{\prime }$, at $\pi $ and $\pi ^{\prime }$, there is a Hida family $\mathbf{%
g}$ of ordinary cuspidal forms passing by $g$ , and, assuming an essentially
nonrestrictive condition \ref{monodim. condition} on $\sigma ,\sigma
^{\prime }$ , it passes by $f$ \ a $p-$adic family (version of Yamagami \ref%
{Definicion familia yamagami}) $\mathbf{f}$ of Hilbert cuspidal forms
nonordinary at $\pi ,\pi ^{\prime }$ of the same slope $\sigma ,\sigma
^{\prime }$. We construct a $p-$adic $L$- function $L_{p}(\mathbf{f,g)}$ \
of the tensor product $\dbigwedge_{\mathbf{f}}^{{}}\otimes \dbigwedge_{%
\mathbf{g}}^{{}}$ of the two algebras defining these familes, such that for
classical weights $s$ and $r$ \ lying on $(k,k)\geq (2,2)$ \ with even $%
0\leq k_{0}<2k-2$ or $k_{0}=k=2$, the value $L_{p}(\mathbf{f},\mathbf{g)(}%
r,s)$ agrees with the $p-$adic Abel-Jacobi map computed in section 5, up to
nonzero Euler factors; while for $k_{0}\geq 2k$, the value of $L_{p}(\mathbf{%
f},\mathbf{g)(}r,s)$ interpolates, up to nonzero Euler factors, a ratio of
periods of modular forms, which, by the second case of \cite{Michino},
theorem 1.1, is a special value of the complex $L-$function of the
representation associated to $g_{r}$ tensored by the representation
associated to the overconvergent $\theta ^{\prime -1-t}f_{s}-\theta
^{-1-t}f_{s}$ \ (for $t\leq -1$) restricted to $X(N)$.

\begin{theorem}
\label{Main theorem copy(1)}\textbf{\ }Consider an $\mathfrak{a}$-Hilbert
cuspidal form $f$ and a cuspidal form $g$ of weight $(2,2)$ and $2$ and
choose $p=\pi \pi ^{\prime }$ as in theorem \ref{TEOREMA PRINCIPAL
INTRODUCCION}, ordinary for $g$ and nonordinary of positive slope $\sigma
,\sigma ^{\prime }$ at $\pi ,\pi ^{\prime }$ . Assuming condition \ref%
{monodim. condition}, we associate to $f$\ and $g$\ families $\mathbf{f}$%
\textbf{\ }and $\mathbf{g}$\ as above so that $f_{s}=f$\ and $g_{r}=g$\ \
for some classic $(r,s)\in $ $\Omega _{\mathbf{f,g}}^{cl}$ . There is a
Garrett-Hida $p-$adic $L-$function $L_{p}(\mathbf{f},\mathbf{g}):\Omega _{%
\mathbf{f,\mathbf{g}}}$\ $\dashrightarrow \mathbb{C}_{p}$\ , defined in \ref%
{PW product}, with the following property: For any pair of classical
characters $(r,s)\in $ $\Omega _{\mathbf{f,g}}^{cl}$ lying over integers $%
(k_{0},k)$ such that $k_{0}$ is even with $2k-2>k_{0}>2$ or $2k-2\geq k_{0}=2
$ , the ordinary cuspidal form $g_{r}$ and the $\mathfrak{a-}$Hilbert
cuspidal form $f_{s}$ , nonordinary of slope $\sigma ,\sigma ^{\prime }$ at $%
\pi ,\pi ^{\prime }$ , have value under the $p-$adic Abel Jacobi map 
\begin{equation}
AJ_{p}(\Delta _{k-2,k_{0}-2})(\omega _{f_{s}}\otimes \eta
_{g_{r}}^{u-r})=(-1)^{t}t!\frac{\mathcal{E}_{0}(g_{r})\mathcal{E}_{1}(g_{r})%
}{\mathcal{E}(f_{s},g_{r})}L_{p}(\mathbf{f},\mathbf{g})(r,s)\text{ ,}
\label{formula teorema cap 6}
\end{equation}%
where $t=$ $k-1-(k_{0}/2)$, and $\mathcal{E}_{0}(g_{r})$ is the non-null
factor $1-\beta _{1}^{2}p^{1-k_{0}}$, as well as $\mathcal{E}_{1}(g_{r}),%
\mathcal{E}(f_{s},g_{r})$ are the factors already defined and proved
non-null in \ref{TEOREMA PRINCIPAL INTRODUCCION} .
\end{theorem}

We have been aware of the work of Y. Liu [23], essentially complementary to
ours. We will discuss the relation of this article with Liu in a future
work. It is a pleasure to thank Henri Darmon for having proposed this
problem to us, which opens, jointly with his recent work with Massimo
Bertolini and Victor Rotger, a promising line of forthcoming reserarch. We
thank Adel Betina by critical reading and suggestions which have improved
the presentation of the article, and thank David Loeffler and Sarah Zerbes
by pointing us a mistake in a previous version of it, and how to fix it.
Thanks also to Donu Arapura, Mladen Dimitrov, Eyal Goren, Tom\'{a}s G\'{o}%
mez, Haruzo Hida, Ignasi Mundet and Fran Presas by their answers to
technical questions. The authors\footnote{%
Ivan Blanco-Chacon: ivan.blanco-chacon@ucd.ie, Ignacio Sols: isols@mat.ucm.es%
} are partially supported by the MICINN project MTM2010-17389 (Spain). The
first author is also partially supported by the Science Foundation Ireland
project 13/IA/1914.

\bigskip To end this introduction, we provide, as preparatory material, a
basic observation on rational homology. Consider an embedding of an
algebraic curve into a surface $j_{C}:C\longrightarrow S$, both smooth and
projective over $\mathbb{C}$, and assume $S$\ has Betti numbers $%
b_{1}=b_{3}=0$. For a point $o\in C$, consider in $C\times C$\ the cycle

\begin{equation*}
C_{o}=C-o\times C-C\times o
\end{equation*}%
$o-$modification of the diagonal cycle $C$ of $C\times C$\ .\ The cycle,
equally denoted $C_{o}$\ , induced in $S\times C$\ by 
\begin{equation*}
j_{C}\times C:C\times C\longrightarrow S\times C
\end{equation*}%
is then null-homologous, since%
\begin{equation}
C_{o}=\partial W_{3}+\dsum E_{i}\times F_{i}=\partial (W_{3}+\dsum
A_{i}\times F_{i})\text{ ,}  \label{key topological formula}
\end{equation}%
where the 1-cycles $E_{i}$\ represent the homology classes of a basis of $%
H_{1}(C)$, the 1-cycles $F_{i\text{ }}$\ represent the homology classes of
their intersection-dual basis, the $A_{i}$\ are 2-chains of $S\ $such that $%
\partial A_{i}=E_{i}$, and $W_{3}$\ denotes a 3-chain in $C\times C$ . In
fact, for $C_{o}$ to be null-homologous, it is not necessary to assume that $%
j_{C}$ is an embedding, but only generically injective, and it is not
necessary to assume that $S$ has null 1-homology, but just that the
generically injective morphism $j_{C}:C\longrightarrow S$ \ kills the
1-homology in the sense that the image of any 1-cycle of $C$\ is
null-homologous in $S$.

\ 

\section{\textbf{Hilbert modular forms and differentials}}

\bigskip

\subsection{\protect\bigskip Hilbert modular surfaces}

Let $\mathfrak{K}$\ be a real quadratic field with\ different ideal $%
\mathfrak{d}_{\mathfrak{K}}$, and fix an inclusion $inc:\mathfrak{K}%
\hookrightarrow \mathbb{R}$\ . Its ring of integers $\mathfrak{o}$ is
invariant under the conjugation $\sigma :\nu \mapsto \nu ^{\prime }$ of $%
\mathfrak{K}$ over $\mathbb{Q}$ . Denote by $\mathfrak{o}^{+}$ the positive
cone of $\mathfrak{o}$ consisting of the totally positive elements $\nu \in 
\mathfrak{o}$ , i.e. those $\nu $ such that $\nu ,\nu ^{\prime }\in \mathbb{R%
}^{+}$, and $\mathfrak{u}^{+}\subseteq \mathfrak{o}^{+}$ the totally
positive units. Assume, just for convenience, that there is a unity of norm $%
-1$, so that $\mathfrak{u}^{+}/(\mathfrak{u}^{+}\mathfrak{)}^{2}$ is trivial
(which can be assumed since the inclusion of the modular curve in the
Hilbert modular surface attached to some $\mathfrak{K}$ is just a tool, and
under this assumption we will avoid restricting to the $\mathfrak{u}^{+}/(%
\mathfrak{u}^{+}\mathfrak{)}^{2}-$invariant part of the cohomologies we
consider). Choose once and for all a narrow class of ideals and a primitive
ideal $\mathfrak{a\subseteq \mathfrak{o}}$ representing it, whose norm $%
N_{1}:=Nm(\mathfrak{a)}$ is prime with the discriminant $d_{\mathfrak{K}}=Nm(%
\mathfrak{d}_{\mathfrak{K}}),$ and denote $\mathfrak{a}^{+}=\mathfrak{a}\cap 
\mathfrak{o}^{+}$. Choose an integer\textit{\ }$N_{2}\geq 4$ prime to $d_{%
\mathfrak{K}}N_{1}$, and denote $N:=d_{\mathfrak{K}}N_{1}N_{2}$.

All schemes in this article are separated and of finite type over its base.
For a ring $R$, let $\mathcal{M}_{R}(\mathfrak{a},N_{2})$, or just $\mathcal{%
M}_{R}(N_{2})$ understanding $\mathfrak{a}$, the irreducible scheme of
dimension 2 (i.e. surface)\ over $R$ which is fine moduli for the Abelian
surfaces $A$ over $R$ with a "real multiplication" or ring monomorphism $%
\alpha :\mathfrak{o\hookrightarrow }End(A)$ ; an ordered $\mathfrak{a-}$%
polarization or $\mathfrak{o}$-linear isomorphism $\delta :Sym_{\mathfrak{o}%
}(A,A^{\vee })\approx \mathfrak{a}$ from symmetric $\mathfrak{o-}$%
homomorphisms from $A$ to the dual\textbf{\ }Abelian surface $A^{\vee }$
applying the cone of positive polarizations $Sym_{\mathfrak{o}}(A,A^{\vee
})^{+}$ to $\mathfrak{a}^{+}$;\ and a $N_{2}-$level structure or $\mathfrak{%
o-}$linear monomorphism of group schemes%
\begin{equation}
\varepsilon :\mu _{N_{2}/R}\otimes _{\mathbb{Z}}\mathfrak{d}_{\mathfrak{K}%
}^{-1}\hookrightarrow A\text{ \ .}  \label{level epsilon}
\end{equation}%
Here $\mu _{N_{2}/R}$ $=SpecR[t]/(t^{N_{2}}-1)$ is the group-scheme over $R$
of $N_{2}-$roots of unity \footnote{%
To be more precise, for any scheme $S$ over $R$, it is $(\mu _{N_{2}}\otimes
_{\mathbb{Z}}\mathfrak{d}_{\mathfrak{K}}^{-1})(S)=\mu _{N_{2}}(S)\otimes _{%
\mathbb{Z}}\mathfrak{d}_{\mathfrak{K}}^{-1}$ . Because of \ref{the basic
split}, under the assumption $\mathfrak{o}\subseteq R$\ , an $N_{2-}$level
structure is a subgroup scheme over $R$\ of $A$\ isomorphic to $\mu
_{N_{2}}\times \mu _{N_{2}}$\ preserved by the action of\textbf{\ \ }$%
\mathfrak{o}$\textbf{\ .}}, i.e. the kernel of the $N_{2}-$power
endomorphism of the multiplicative group scheme $\mathbb{G}_{m/R}$ (cf. 3.2,
3.3 and 5.4 of \cite{And Go} or lemma 1.1 of \cite{Goren unramified}, for
instance)\footnote{%
Since we will refer frequently to the article of Katz \cite{Split Katz}, we
warn the reader that our smooth moduli $\mathcal{M}_{R}(N_{2})$ is denoted $%
\mathcal{M}(\mathfrak{a,\Gamma }_{00}(N_{2}))_{R}$ in that reference.}. For $%
R\longrightarrow R^{\prime }$, it is 
\begin{equation*}
\mathcal{M}_{R^{\prime }}(\mathfrak{a},N_{2})=\mathcal{M}_{R}(\mathfrak{a}%
,N_{2})\times _{Spec\text{ }R}Spec\text{ }R^{\prime }\text{ .}
\end{equation*}

\bigskip For level $N_{2}\leq 3$, there is a coarse moduli for these
objects, still denoted $\mathcal{\mathcal{M}}_{R}(N_{2})$, although no
longer fine (see theorem 4.9 of \textbf{\cite{Hida libro},} for instance).
In the particular case of "no level" $\mathcal{M}_{R}(1)$, the fact that $%
\mathcal{\mathcal{M}}_{R}(N_{2})$, for $N_{2}\geq 4$, has a universal
family, provides us thus with a finite morphism%
\begin{equation}
pr_{\mathcal{M}_{R}}:\mathcal{\mathcal{M}}_{R}(N_{2})\twoheadrightarrow 
\mathcal{M}_{R}(1)  \label{etale to coarse moduli}
\end{equation}%
because of the weak universal property of coarse moduli spaces.

For $N_{2}\geq4$, denote 
\begin{equation*}
pr_{\mathcal{M}_{R}(N_{2})}:A_{\mathcal{M}_{R}(N_{2})}^{U}\longrightarrow 
\mathcal{M}_{R}(N_{2})
\end{equation*}
the representing object, i.e. the universal Abelian surface with real $%
\mathfrak{o-}$multiplication $\alpha_{R}^{U}$, ordered $\mathfrak{a-}$%
polarization $\delta_{R}^{U}$ and $N_{2}-$level structure $\varepsilon^{U}$ .

Because of the assumption that $d_{\mathfrak{K}}$is invertible in $R$, the
"Rapoport" rank one $\mathfrak{o\otimes }_{\mathbb{Z}}\mathcal{O}_{\mathcal{M%
}_{R}(N_{2})}$ -sheaf 
\begin{equation*}
\mathcal{R}_{\mathcal{M}_{R}(N_{2})}:=pr_{\mathcal{M}_{R}(N_{2})\ast }\Omega
_{A_{\mathcal{\mathcal{M}}_{R}(N_{2})}^{U}/\mathcal{\mathcal{M}}_{R}(N_{2})}%
\text{ }
\end{equation*}%
is locally free, i.e. the "Rapoport" open subscheme where it is locally free
is the whole moduli, so that $\mathcal{R}_{\mathcal{M}_{R}(N_{2})}$ is
locally free of rank two as a $\mathcal{O}_{\mathcal{\mathcal{M}}%
_{R}(N_{2})} $ $\ $-sheaf. The moduli $\mathcal{M}_{R}(N_{2})$ is smooth
over $R$, as it is so, in general, for the Rapoport locus. Assume, from now
on, that\textit{\ }$\mathfrak{o}\subseteq R$. Then the monomorphism provided
by the two natural morphisms $\mathfrak{o}\otimes _{\mathbb{Z}%
}R\longrightarrow R$ , inclusion and conjugation $\sigma $, is an isomorphism%
\textbf{:} 
\begin{equation}
\mathfrak{o}\otimes _{\mathbb{Z}}R\approx R\oplus R\text{ }
\label{the basic split}
\end{equation}%
and in fact $M\otimes _{\mathbb{Z}}\mathfrak{j\approx }M\oplus M$ for any $R$%
-module $M$ and fractional ideal $\mathfrak{j}$\ \textbf{(}cf.\textbf{\ }%
2.0.3 and 2.1.1 of\textbf{\ \cite{Split Katz}}. The two natural
homomorphisms 
\begin{equation*}
\mathfrak{o\otimes }_{\mathbb{Z}}\mathcal{O}_{\mathcal{\mathcal{M}}%
_{R}(N_{2})}\longrightarrow \mathcal{O}_{\mathcal{\mathcal{M}}_{R}(N_{2})}
\end{equation*}%
induce two homomorphisms of sheaves of multiplicative groups%
\begin{equation*}
(\mathfrak{o\otimes }_{\mathbb{Z}}\mathcal{O}_{\mathcal{\mathcal{M}}%
_{R}(N_{2})})^{\ast }\longrightarrow \mathcal{O}_{\mathcal{\mathcal{M}}%
_{R}(N_{2})}^{\ast }
\end{equation*}%
which associate to any invertible $\mathfrak{o\otimes \mathcal{O}}_{\mathcal{%
\mathcal{M}}_{R}(N_{2})}$ -sheaf, two invertible $\mathfrak{\mathcal{O}}_{%
\mathcal{\mathcal{M}}_{R}(N_{2})}$ -sheaves of which it is thus the direct
sum. Applying this to the Rapoport invertible $\mathfrak{o\otimes \mathcal{O}%
}_{\mathcal{\mathcal{M}}_{R}(N_{2})}$ -sheaf, we obtain a decomposition into
line bundles (cf. 2.0.9 of \cite{Split Katz} ) 
\begin{equation*}
\mathcal{R}_{\mathcal{M}_{R}(N_{2})}\mathcal{\approx }L_{\mathcal{\mathcal{M}%
}_{R}(N_{2})}^{(1,0)}\oplus L_{\mathcal{\mathcal{M}}_{R}(N_{2})}^{(0,1)}%
\text{ .}
\end{equation*}%
This allows us to define the modular line bundles 
\begin{equation*}
L_{\mathcal{\mathcal{M}}_{R}(N_{2})^{{}}}^{(k,k^{\prime
})}:=\dbigotimes^{k}L_{\mathcal{\mathcal{M}}_{R}(N_{2})^{{}}}^{(1,0)}\otimes
\dbigotimes^{k^{\prime }}L_{\mathcal{\mathcal{M}}_{R}(N_{2})}^{(0,1)}
\end{equation*}%
of weight $(k,k^{\prime })\in \mathbb{Z}\times \mathbb{Z}$ . The $n-$th
tensor power $\mathcal{R}_{\mathcal{\mathcal{M}}_{R}(N_{2})^{{}}}^{\otimes
n} $ is an invertible $\mathcal{O}_{\mathcal{\mathcal{M}}_{R}(N_{2})^{{}}}%
\mathcal{\otimes \mathfrak{o-}}$sheaf which thus decomposes (cf. 2.0.10 of 
\cite{Split Katz}) as a direct sum of line bundles%
\begin{equation*}
\mathcal{R}_{\mathcal{\mathcal{M}}_{R}(N_{2})^{{}}}^{\otimes n}\approx L_{%
\mathcal{\mathcal{M}}_{R}(N_{2})^{{}}}^{(n,0)}\oplus L_{\mathcal{\mathcal{M}}%
_{R}(N_{2})}^{(0,n)}\text{ .}
\end{equation*}%
This must not be confused with the $\mathcal{O}_{\mathcal{\mathcal{M}}%
_{R}(N_{2})^{{}}}-$symmetric $n-$th power of the rank two $\mathcal{O}_{%
\mathcal{\mathcal{M}}_{R}(N_{2})^{{}}}$-bundle $\mathcal{R}_{\mathcal{%
\mathcal{M}}_{R}(N_{2})^{{}}}$, which is a $\mathcal{O}_{\mathcal{\mathcal{M}%
}_{R}(N_{2})^{{}}}$-bundle of rank $n+1$ denoted 
\begin{equation}
\mathcal{R}_{\mathcal{\mathcal{M}}_{R}(N_{2})^{{}}}^{n}\mathcal{\mathcal{%
\approx }}\dbigoplus\limits_{\substack{ k,k^{\prime }\geq 0  \\ k+k^{\prime
}=n}}^{{}}L_{\mathcal{\mathcal{M}}_{R}(N_{2})}^{(k,k^{\prime })\text{ }}
\label{symmetric power Rappoport}
\end{equation}%
According to the material recalled, for instance, at the beginning of proof
11.9 of \cite{And Go}, there is a surface $\overline{\mathcal{\mathcal{M}}%
_{R}(N_{2})}$ normal and projective over $R$ , called "Satake
compactification", which minimally compactifies the moduli $\mathcal{%
\mathcal{M}}_{R}(N_{2})$. The singularities of $\overline{\mathcal{\mathcal{M%
}}_{R}(N_{2})}$ \ are the points $\overline{\mathcal{\mathcal{M}}_{R}(N_{2})}%
-\mathcal{\mathcal{M}}_{R}(N_{2})$ outside the moduli\textbf{.} We choose,
once and for all, a toroidal desingularization, over $R$,%
\begin{equation*}
Y_{R}(\mathfrak{a},N_{2})\twoheadrightarrow \overline{\mathcal{M}_{R}}(%
\mathfrak{a,}N_{2})\text{,}
\end{equation*}%
We denote by $D^{c}$ the transform divisor of the cusps (this should be
denoted $D_{Y_{R}(N_{2})}^{c}$ to recall to which surface it belongs, but we
spare, in general, subindexes of surface divisors, as always clear). The
morphism (\ref{etale to coarse moduli}) extends obviously to 
\begin{equation*}
pr_{Y_{R}}:Y_{R}(N_{2})\longrightarrow Y_{R}(1)=Y_{R}
\end{equation*}%
(we will denote $Y_{R}(\mathfrak{a},N_{2})$ by just $Y_{R}(N_{2})$, when $%
\mathfrak{a}$ is clear by the context). According to an idea of Mumford
(recalled in (4.4) of \textbf{\cite{Goren unramified} }and in the discussion
before theory. 4.1 of \textbf{\cite{Hida libro}),} the universal family $A_{%
\mathcal{\mathcal{M}}_{R}(N_{2})^{{}}}^{U}$ extends to a universal family 
\begin{equation}
pr_{Y_{R}(N_{2})}:A_{Y_{R}(N_{2})}^{U}\longrightarrow Y_{R}(N_{2})
\label{universal family on Y}
\end{equation}%
\textbf{\ }of semiabelian schemes 
\begin{equation*}
1\longrightarrow \mathbb{G}_{m,R}^{n}\longrightarrow A\longrightarrow
B\longrightarrow 1
\end{equation*}%
of dimension 2 with $\mathfrak{o-}$multiplication, and $\mathcal{\mathcal{R}}%
_{\mathcal{\mathcal{M}}_{R}(N_{2})\text{ }}$ extends to the rank one locally
free $\mathfrak{o\otimes }_{\mathbb{Z}}\mathcal{O}_{Y_{R}(N_{2})}$-sheaf 
\begin{equation}
\mathcal{R}_{Y_{R}(N_{2})}:=pr_{Y_{R}(N_{2})}^{\prime }{}_{\ast }\Omega
_{A_{Y_{R}(N_{2})}^{U}/Y_{R}(N_{2})}\subseteq pr_{Y_{R}(N_{2})}{}_{\ast
}\Omega _{A_{Y_{R}(N_{2})}^{U}/Y_{R}(N_{2})}\text{ }  \label{crisis}
\end{equation}%
of relative differentials which are $A_{Y_{R}(N_{2})}^{U}-$invariant, i.e.
invariant by the action extending the translation in the Abelian fibres (in
which the invariance needs not to be imposed, cf. \cite{Mladen}, 1.6 and 5.2
iii) of \cite{Mladen y Tilouine}). This coincides with $e^{\ast }(\Omega
_{A_{Y_{R}(N_{2})}^{U}/Y_{R}(N_{2})})$ for the unity section $e$ of the
projection $pr_{Y_{R}(N_{2})}$.

The rank-one locally free $\mathcal{O}_{Y_{R}(N_{2})}$ -sheaves $%
L_{Y_{R}(N_{2})}^{(k,k^{\prime })}$ attached to $\mathcal{R}_{Y_{R}(N_{2})}$
as above, extend the modular line bundles $L_{\mathcal{M}%
_{R}^{{}}(N_{2})}^{(k,k^{\prime })}$ (cf. (4.4) and comments preceding
theorem 4.2 in \cite{Goren unramified}). In the case of parallel weight,
these are the liftings of the line bundles $L_{\overline{\mathcal{\mathcal{M}%
}_{R}(N_{2})}}^{(k,k)}$ on the minimal compactification (cf. remark 5.6 of 
\cite{And Go}).

Extending what has been said for the non-compactified moduli, the rank two
locally free $\mathcal{O}_{Y_{R}}$-sheaf $\mathcal{R}_{Y_{R}(N_{2})}$ splits
as 
\begin{equation}
\mathcal{R}_{Y_{R}(N_{2})}\approx L_{Y_{R}(N_{2})}^{(1,0)}\oplus
L_{Y_{R}(N_{2})}^{(0,1)}\text{ }  \label{main decomposition}
\end{equation}
Again, we will consider the rank one $\mathfrak{o\otimes}_{\mathbb{Z}}%
\mathcal{O}_{Y_{R}(N_{2})}$-sheaf $\ n-$tensor power $\mathcal{R}%
_{Y_{R}(N_{2})}^{\otimes n}$ of $\mathcal{R}_{Y_{R}(N_{2})}$ which, as rank
two locally free $\mathcal{O}_{Y_{R}}$-sheaf, is 
\begin{equation*}
\mathcal{R}_{Y_{R}(N_{2})}^{\otimes n}\approx L_{Y_{R}(N_{2})}^{(n,0)}\oplus
L_{Y_{R}(N_{2})}^{(0,n)}\text{ ,}
\end{equation*}
and the $n-$symmetric power 
\begin{equation}
\text{ }\mathcal{R}_{Y_{R}(N_{2})}^{n}\approx\dbigoplus \limits_{\substack{ %
k,k^{\prime}\geq0  \\ k+k^{\prime}=n}}L_{Y_{R}(N_{2})}^{(k,k^{\prime})\text{ 
}}\text{ for }n\geq 0  \label{symmetric powers de R sobre Y}
\end{equation}
of \ $\mathcal{R}_{Y_{R}(N_{2})}$ as rank two locally free $\mathcal{O}%
_{Y_{R}}$-sheaf, this having rank $n+1$.

The $R-$module $M_{k,k^{\prime }}(R,\mathfrak{a,}N_{2})$ of $\mathfrak{a}$%
-polarized Hilbert modular forms of level $N_{2}$ over $R$, and weight $%
(k,k^{\prime })$, or just $\mathfrak{a}$-Hilbert modular forms if the rest
is understood, can be defined as the $R$-module of sections%
\begin{equation*}
M_{k,k^{\prime }}(\mathfrak{a,}N_{2},R):=\Gamma L_{\mathcal{\mathcal{M}}%
_{R}^{{}}(N_{2})}^{(k,k^{\prime })}\text{ }=\Gamma
L_{Y_{R}(N_{2})}^{(k,k^{\prime })}
\end{equation*}%
(cf. \cite{And Go} 5.1, for instance). The equality is obvious by normality
of $\overline{\mathcal{\mathcal{M}}_{R}(N_{2})}$ in the parallel case $%
k=k^{\prime }$, and is called Koecher principle in general. As this principle%
\textbf{\ }is proved in \cite{Rapoport} 4.9 by checking it on open affine
subsets of $Y_{Y_{R}(N_{2})}^{(k,k^{\prime })}$, it holds also for sections
of $L_{Y_{R}(N_{2})}^{(k,k^{\prime })}$ on any open set $U\subseteq $ $%
Y_{R}(N_{2})$ , i.e. they\texttt{\ }coincide with he sections on $U\cap 
\mathcal{\mathcal{M}}_{R}(N_{2})$ . For general weight, we can provisionally
define the $R-$submodule of the Hilbert cuspidal forms, as 
\begin{equation}
S_{k,k^{\prime }}(\mathfrak{a,}N_{2},R)=\Gamma L_{Y(N_{2})}^{(k,k^{\prime })}%
\text{ }(-D^{c})  \label{unparallel Koecher}
\end{equation}%
(which will correspond to the familiar definition in terms of $q-$%
expansions). Lemma 6.3 of \cite{VdG}\textbf{\ }asserts that $M_{k,k^{\prime
}}(\mathfrak{a,}N_{2},\mathbb{C})$ is null whenever $k$ or $k^{\prime }$ are
nonpositive.

The smooth variety $\mathcal{\mathcal{M}}_{R}(N_{2})(\mathbb{C)}$ of complex
points of $\mathcal{\mathcal{M}}_{R}(N_{2})$\textbf{, }which we sloppy
denote $\mathcal{\mathcal{M}}_{\mathbb{C}}(N_{2})$, can be seen as the
quotient $\Gamma_{1}^{1}(\mathfrak{\mathfrak{a}},N_{2})\backslash\mathfrak{H}%
_{\mathfrak{K}}$ (cf. 2.1 of \textbf{\cite{Mladen y Tilouine}, }for
instance) of

\begin{equation*}
\mathfrak{H}_{\mathfrak{K}}:=[(z_{1},z_{2})\in \mathfrak{K\otimes }_{\mathbb{%
Z}}\mathbb{C\mid }im(z_{i})\in (\mathfrak{K\otimes }_{\mathbb{Z}}\mathbb{R)}%
_{+}]\approx \mathfrak{H}^{2}
\end{equation*}%
by the left action of the subgroup $\Gamma _{1}^{1}(\mathfrak{\mathfrak{a}}%
,N_{2})$ of the group $\Gamma _{0}^{1}(\mathfrak{\mathfrak{a}},N_{2})$ of
unimodular matrices of the form 
\begin{equation*}
\gamma =\left( 
\begin{array}{cc}
a & b \\ 
c & d%
\end{array}%
\right) \in \left( 
\begin{array}{cc}
\mathfrak{o} & \mathfrak{a^{\ast }} \\ 
\mathfrak{ad}_{\mathfrak{K}}N_{2} & \mathfrak{o}%
\end{array}%
\right) _{\det =1}
\end{equation*}%
given by condition $d\equiv 1$ mod. $N_{2}\mathfrak{\ }$, a quotient which
is uniformized by an algebraic variety defined over $\mathbb{Z[}\frac{1}{N}%
\mathbb{]}$ , thus defined over $R$. \ The reason is that $N$ is divisible
by all primes whose localization do not make the subgroup the full special
linear group , cf. X. 4 \cite{VdG}. We have kept the notations $\Gamma
_{1}^{1}(\mathfrak{\mathfrak{a}},N_{2})\subseteq $ $\Gamma _{0}^{1}(%
\mathfrak{\mathfrak{a}},N_{2})$ just for the sake of generality, although
these groups are equal, in our case $\mathfrak{u}^{+}/(\mathfrak{u}%
^{+})^{2}=1$, to the groups $\Gamma _{1}(\mathfrak{\mathfrak{a}}%
,N_{2})\subseteq $ $\Gamma _{0}(\mathfrak{\mathfrak{a}},N_{2})$ which in
general are bigger, having $\mathfrak{u}^{+}/(\mathfrak{u}^{+})^{2}$ as a
two-elements quotient which acts on $\Gamma _{1}^{1}(\mathfrak{\mathfrak{a}}%
,N_{2})\backslash \mathfrak{H}_{\mathfrak{K}}$ having $\Gamma _{1}(\mathfrak{%
\mathfrak{a}},N_{2})\backslash \mathfrak{H}_{\mathfrak{K}}$ \ as a quotient
(This last surface is coarse moduli of the analogous data but replacing
"polarizations" by "classes of polarizations")

In terms of this modular group, $M_{k,k^{\prime}}(\mathfrak{a,}N_{2},\mathbb{%
C})$ consists of the homolorphic functions $f:\mathfrak{H}_{\mathfrak{K}%
}\longrightarrow\mathbb{C}$ such that%
\begin{equation}
\text{\ }f(\frac{az_{1}+b}{cz_{1}+d},\frac{a^{\prime}z_{2}+b^{\prime}}{%
c^{\prime}z_{2}+d^{\prime}})=(cz_{1}+d)^{k}(c^{\prime}z_{2}+d^{\prime
})^{k^{\prime}}f(z_{1},z_{2})
\end{equation}
for all $\gamma\in\Gamma_{1}^{1}(\mathfrak{a,}N_{2})$. For $R\subseteq 
\mathbb{C}$, the $R-$submodule $M_{k,k^{\prime}}(\mathfrak{a,}N_{2},R)$
consists of those $f$ \ whose Fourier coefficients, as periodic function,
lie in $R$.

For general $R$ as in our hypothesis, \cite{And Go}, a set of cusps of $%
\mathcal{\mathcal{M}}_{R}(N_{2})$ called "non-ramified cusps"\textbf{,} is
considered in bijection with 4-tuples $(\mathfrak{\mathfrak{b}}_{1},%
\mathfrak{\mathfrak{b}}_{2},\varepsilon ,j)$ where $\mathfrak{\mathfrak{b}}%
_{1},\mathfrak{\mathfrak{b}}_{2}$ are fractional ideals such that $\mathfrak{%
\mathfrak{b}}_{1}\mathfrak{\mathfrak{b}}_{2}^{-1}=\mathfrak{a}$ ; $\iota
:N_{2}^{-1}\mathfrak{o/\mathfrak{o\approx }}N_{2}^{-1}\mathfrak{b}_{1}^{-1}%
\mathfrak{/b}_{1}^{-1}$ is a $\mathfrak{\mathfrak{o}}$-linear isomorphism;
and $j:\mathfrak{\mathfrak{b}}_{1}\otimes R\approx $ $\mathfrak{o}\otimes R$
is an $\mathfrak{\mathfrak{o\otimes }}R$ -isomorphism (cf. \cite{And Go}
(6.4) for instance). Among them, we will always refer to the standard cusp
given by $\mathfrak{\mathfrak{b}}_{1}=\mathfrak{\mathfrak{o}}$ and $%
\mathfrak{\mathfrak{b}}_{2}=\mathfrak{\mathfrak{a}}^{-1},$ and we will
consider as isomorphism $\iota $ the identity in $N_{2}^{-1}\mathfrak{o/%
\mathfrak{o}}$ and as isomorphism $j$ the identity in $\mathfrak{o}\otimes
R. $ We associate to an $\mathfrak{a-}$Hilbert modular form $f$ of weight $%
(k,k^{\prime })$ defined over $R$, as expansion near this cusp, a formal $q-$%
expansion (cf.\cite{And Go} 6.8 and \cite{Split Katz} (1.2.12), for
instance) 
\begin{equation}
f(q):=\dsum\limits_{\nu \in (\mathfrak{a}^{-1})^{+}\cup \{0\}}^{{}}a_{\nu
}q^{\nu }\text{ , with }a_{\nu }\in R\text{ .}  \label{expansion}
\end{equation}%
\textbf{\ }where the set of subindexes $(\mathfrak{a}^{-1})^{+}\cup \{0\}$
contains $(\mathfrak{\mathfrak{o}}^{-1})^{+}\cup \{0\}$, thus $\mathfrak{o}%
^{+}\cup \{0\}$, because $\mathfrak{a\subseteq o}$ (Here $a_{\nu
}=a_{\varepsilon ^{2}\nu }$ , for all positive units $\varepsilon $, so ,
under our assumption that $\mathfrak{u}^{+}=(\mathfrak{u}^{+})^{2}$ , the
coefficient $a_{\nu \text{ }}$ only depends on the principal ideal $(\nu )$
). This $q-$expansion can be seen as attached to a nowhere vanishing
relative differential $\omega _{can}=\varepsilon _{\ast }(\frac{dt}{t})$
given by the level $\varepsilon $ in (\ref{level epsilon}) ( called \ $%
\omega _{\mathfrak{a}}(j)$ in comment to (1.2.10) of \cite{Split Katz}) in
the punctured complete neighborhood of the standard cusp of $\overline{%
\mathcal{\mathcal{M}}_{R}(N_{2})}$, i.e. a section trivializing $\mathcal{R}%
_{Y_{R}(N_{2})}$ as rank one $\mathfrak{o\otimes }_{\mathbb{Z}}\mathcal{O}%
_{Y_{R}(N_{2})}$-sheaf, which induces a canonical trivialization of each
locally free sheaf $L_{Y_{R}(N_{2})}^{(k,k^{\prime })}$ on it (so that a
section becomes a function, and the $q-$expansion $f(q)$ is the evaluation
of the extended function $f$ at the Tate point. In fact, $f(q)$ belongs to
the complete local ring of the standard cusp in $\overline{\mathcal{\mathcal{%
M}}_{R}(N_{2})}$ and is independent of the toroidal embedding chosen to
desingularize the standard cusp, cf. 6.5 to 6.8 of \cite{And Go}).\ Now the
submodule $S_{k,k^{\prime }}(\mathfrak{a,}N_{2},R)\subseteq M_{k,k^{\prime
}}(\mathfrak{a,}N_{2},R)$ of cuspidal forms can be recovered by the
condition $a_{0}=0$ in the expansion of $f(q)$ at the standard cusp, and the
analogous vanishing at\ all the\ other cusps. In fact all Hilbert modular
forms of unparallel weight are cuspidal, as immediate consequence of 11.1
iii) in \cite{And Go}(cf. comment before theorem 4.2 of \cite{Goren
unramified}). Being $R\supseteq \mathfrak{o}$ of zero characteristic, the
homomorphism from $M_{k,k^{\prime }}(\mathfrak{a,}N_{2},R)$ to the $R-$%
module of formal series indexed by $(\mathfrak{a}^{-1})^{+}\cup \{0\}$ is
injective.

\subsection{Derivatives}

The Kodaira-Spencer isomorphism (obtained, for instance, from the
Gauss-Manin connection, which we will recall in section 4) takes, under our
assumption that $d_{\mathfrak{K}}$\ is invertible in\textbf{\ }$R$ , the
form of an isomorphism of invertible $\mathfrak{o\otimes }\mathcal{\mathcal{%
\mathcal{O}}}_{\mathcal{M}(N_{2})}-$sheaves 
\begin{equation*}
KS:\Omega _{\mathcal{M}_{R}(N_{2})/R}\approx \mathcal{R}_{\mathcal{M}%
_{R}(N_{2})}\otimes _{\mathfrak{o\otimes \mathcal{O}}_{\mathcal{M}%
_{R}(N_{2})}}(\mathcal{R}_{\mathcal{M}_{R}(N_{2})}\otimes _{\mathfrak{o}}%
\mathfrak{d}_{\mathfrak{K}}\mathfrak{a}^{-1})\text{ .}
\end{equation*}%
In case $\mathfrak{d}_{\mathfrak{K}}\mathfrak{a}^{-1}$ is prime to any
integer prime $p$ which is a non-unit in $R$, this is isomorphic as $%
\mathcal{\mathcal{\mathcal{O}}}_{\mathcal{M}_{R}(N_{2})}-$sheaf to $\mathcal{%
R}_{\mathcal{M}_{R}(N_{2})}^{\otimes 2}$, so we have the equally denoted
isomorphism of locally free $\mathcal{\mathcal{\mathcal{O}}}_{\mathcal{M}%
(N_{2})}-$sheaves%
\begin{equation*}
KS:\Omega _{\mathcal{M}_{R}(N_{2})/R}\approx L_{\mathcal{M}%
_{R}(N_{2})}^{(2,0)\text{ }}\oplus \text{ }L_{\mathcal{M}_{R}(N_{2})}^{(0,2)%
\text{ }}
\end{equation*}%
(cf. 2.1.1 to 2.1.3 in \cite{Split Katz} and \cite{Mladen y Tilouine} 5.2
iii) )\footnote{%
The philosophy here is that, although $\mathcal{\mathcal{M}}_{R}(N_{2})$
does not split as described by local coordinates $z_{\pi },z_{\pi ^{\prime
}} $ -as if it were the cartesian product of two modular curves-, at the
infinitesimal level it does so, i.e. the tangent bundle $T_{\mathcal{%
\mathcal{M}}_{R}(N_{2})/R}$ splits as direct sum $L_{\mathcal{\mathcal{M}}%
_{R}(N_{2})}^{(-2,0)\text{ }}\oplus $ $L_{\mathcal{\mathcal{M}}%
_{R}(N_{2})}^{(0,-2)\text{ }}$of what can be thought of the tangent bundles
of such modular curves.}. We will make from now on this assumption which
will be empty in case $R$ is a field, or equivalent to $p$ being coprime
with $\mathfrak{a}$ in case $R$ is isomorphic to $\mathbb{Z}_{p}$ , the two
cases considered in our applications. The Kodaira-Spencer isomorphism
induces 
\begin{equation*}
\mathcal{\omega }_{\mathcal{M}(N_{2})/R}\approx L_{\mathcal{M}(N_{2})}^{(2,2)%
\text{ }}\text{ .}
\end{equation*}

This extends to the whole of $Y_{R}(N_{2})$ as the sheaf 
\begin{equation}
\Omega _{Y_{R}(N_{2})/R}(\log D^{c})\approx \mathcal{R}_{Y_{R}(N_{2})}^{%
\otimes 2}\approx L_{Y_{R}(N_{2})}^{(2,0)\text{ }}\oplus
L_{Y_{R}(N_{2})}^{(0,2)\text{ }}\text{ ,}  \label{log dif splits}
\end{equation}%
of $\log $ $D^{c}-$differentials on $Y_{R}(N_{2})$ over $R$, i.e. those
having, as well as their exterior derivatives, possible simple poles at $%
D^{c}$ (cf. for instance \cite{VdG} p. 62 and \cite{Mladen y Tilouine} 5.2
iii). As a consequence of (\ref{log dif splits}) and of the fact that $%
M_{2,0}(\mathfrak{a,}N_{2},\mathbb{C})$ and $M_{2,0}(\mathfrak{a,}N_{2},%
\mathbb{C})$ are null (\cite{VdG} ch. I\textbf{),} the irregularity $%
q=h^{1,0}=h^{0,1}$ of the complex surface $Y_{\mathbb{C}}(N_{2})(\mathbb{C)}$
is null, so that this surface has null Betti numbers $b_{1}=b_{3}=2q=0$ the
key fact of our program.

The "residue" epimorphism which assigns, to each $\log $ $D^{c}-$%
differential $1-$form, its eventual simple pole, provides a short exact
sequence 
\begin{equation*}
0\longrightarrow \Omega _{Y_{R}(N_{2})/R}\longrightarrow \Omega
_{Y_{R}(N_{2})/R}(\log D^{c})\longrightarrow \mathcal{O}_{D^{c}}%
\longrightarrow 0\text{ \ .}
\end{equation*}%
From this we obtain, by comparing first Chern classes,%
\begin{equation*}
\omega _{Y_{R}(N_{2})/R}\approx L_{Y_{R}(N_{2})}^{(2,2)\text{ }}(-D^{c})%
\text{ .}
\end{equation*}%
This associates bijectively a 2-differential of the scheme $Y_{R}(N_{2})$
over $R$ to each $\mathfrak{a}$-Hilbert cuspidal form of weight $(2,2)$ and
level $N_{2}$ over $R$ .

Recall that the exterior derivative 
\begin{equation}
d=(\partial ,\partial ^{\prime })=\mathcal{O}_{\mathcal{M}%
(N_{2})}\longrightarrow \Omega _{\mathcal{M}(N_{2})}\approx L_{\mathcal{M}%
(N_{2})}^{(2,0)\text{ }}\oplus L_{\mathcal{M}(N_{2})}^{(0,2)\text{ }}\text{ }
\label{derivada}
\end{equation}%
is the unique homomorphism of sheaves of $R$-modules satisfying, on local
sections, the Leibniz law (this way of writing derivatives is a slight abuse
of notation in that we identify vector bundles with the associated locally
free sheaves of $R-$modules). It acts on the $q-$expansion of local sections
of these bundles near (in the punctured complete neighborhood of) the
standard cusp,\ as 
\begin{equation*}
f(q)=\dsum\limits_{\nu \in (\mathfrak{a}^{-1})^{+}\cup \{0\}}^{{}}a_{\nu
}q^{\nu }\mapsto (\theta f(q),\theta ^{\prime }f(q))\text{ ,}
\end{equation*}%
where%
\begin{equation}
\theta f(q)=\dsum\limits_{\nu \in (\mathfrak{a}^{-1})^{+}\cup \{0\}}^{{}}\nu
a_{\nu }q^{\nu }\text{ and }\theta ^{\prime }f(q)=\dsum\limits_{\nu \in (%
\mathfrak{a}^{-1})^{+}\cup \{0\}}^{{}}\nu ^{\prime }a_{\nu }q^{\nu }\text{ .}
\label{expansion first der. of diff}
\end{equation}%
This derivative (\ref{derivada}) has null composition with the exterior
derivative%
\begin{equation}
d=\partial ^{\prime }-\partial :\Omega _{\mathcal{M}(N_{2})}\approx L_{%
\mathcal{M}(N_{2})}^{(2,0)\text{ }}\oplus \text{ }L_{\mathcal{M}%
(N_{2})}^{(0,2)\text{ }}\longrightarrow \omega _{\mathcal{M}(N_{2})}\approx
L_{\mathcal{M}(N_{2})}^{(2,2)\text{ }}\text{ ,}  \label{derivada logaritmica}
\end{equation}%
sum of 
\begin{equation*}
\partial ^{\prime }:L_{\mathcal{M}(N_{2})}^{(2,0)\text{ }}\longrightarrow L_{%
\mathcal{M}(N_{2})}^{(2,2)\text{ }}\text{ \ and }-\partial :\text{ }L_{%
\mathcal{M}(N_{2})}^{(0,2)\text{ }}\longrightarrow L_{\mathcal{M}%
(N_{2})}^{(2,2)\text{ }}
\end{equation*}%
It acts on $q-$expansions as%
\begin{equation}
\begin{array}{c}
\begin{array}{c}
(\dsum\limits_{\nu \in (\mathfrak{a}^{-1})^{+}\cup \{0\}}^{{}}b_{\nu }q^{\nu
},\dsum\limits_{\nu \in (\mathfrak{a}^{-1})^{+}\cup \{0\}}^{{}}c_{\nu
}q^{\nu })\mapsto \dsum\limits_{\nu \in (\mathfrak{a}^{-1})^{+}\cup
\{0\}}(\nu ^{\prime }b_{\nu }-\nu c_{\nu })q^{\nu }= \\ 
\theta ^{\prime }(\dsum\limits_{\nu \in (\mathfrak{a}^{-1})^{+}\cup
\{0\}}b_{\nu }q^{\nu })-\theta (\dsum\limits_{\nu \in (\mathfrak{a}%
^{-1})^{+}\cup \{0\}}c_{\nu }q^{\nu })%
\end{array}%
\end{array}
\label{expansion second derivative of dif2}
\end{equation}%
thus providing the complex $\Omega _{\mathcal{M}(N_{2})}^{\bullet }$of
sheaves of $R-$modules on $\mathcal{M}(N_{2})$%
\begin{equation*}
\mathcal{O}_{\mathcal{M}(N_{2})}\overset{d=(\partial ,\partial ^{\prime })}{%
\longrightarrow }\Omega _{\mathcal{M}(N_{2})}=L_{\mathcal{M}(N_{2})}^{(2,0)%
\text{ }}\oplus \text{ }L_{\mathcal{M}(N_{2})}^{(0,2)\text{ }}\overset{%
d=\partial ^{\prime }-\partial }{\longrightarrow }\omega _{\mathcal{M}%
(N_{2})}=L_{\mathcal{M}(N_{2})}^{(2,2)\text{ }}\text{ .}
\end{equation*}%
This extends to the complex $\mathcal{\Omega }_{Y_{R}(N_{2})}^{\bullet
}(\log D^{c})$\textbf{\ }of $\log D^{c}-$differentials on $Y_{R}(N_{2})$

\begin{eqnarray*}
\mathcal{O}_{Y_{R}(N_{2})}\overset{d=(\partial ,\partial ^{\prime })}{%
\longrightarrow }\mathcal{\Omega }_{Y_{R}(N_{2})}(\log D^{c}) &=& \\
L_{Y_{R}(N_{2})}^{(2,0)\text{ }}\oplus \text{ }L_{Y_{R}(N_{2})}^{(0,2)\text{ 
}}\overset{d=\partial ^{\prime }-\partial }{\longrightarrow }\mathcal{\omega 
}_{Y_{R}(N_{2})}(D^{c}) &=&L_{Y_{R}(N_{2})}^{(2,2)\text{ }}\text{ .}
\end{eqnarray*}

\bigskip

\subsection{\textbf{\ Morphisms from modular curves}}

Recall we denote rather sloppy by $\mathcal{M}_{\mathbb{C}}(N_{2})$ the
smooth variety of complex points of this scheme, and that it is the quotient
of $\mathfrak{H}_{\mathfrak{K}}$ by the left action of $\Gamma _{1}^{1}(%
\mathfrak{\mathfrak{a}},N_{2})$. This action factors through the
projectivized group

\begin{equation*}
\mathbb{P}\Gamma_{1}^{1}(\mathfrak{\mathfrak{a}},N_{2})=\Gamma_{1}^{1}(%
\mathfrak{\mathfrak{a}},N_{2})/\{\pm1\}\text{.}
\end{equation*}

It is proved in prop. 1.1 of chap V in \cite{VdG} the existence of a morphism

\begin{equation*}
\mathcal{M}_{\mathbb{C}}(\Gamma _{0}(N_{1}))\longrightarrow \mathcal{M}_{%
\mathbb{C}}(\mathfrak{\mathfrak{a}},1)
\end{equation*}%
inducing in smooth compactifications

\begin{equation}
j_{X_{\mathbb{C}}(\Gamma _{0}(N_{1}))}:X_{\mathbb{C}}(\Gamma
_{0}(N_{1}))\longrightarrow Y_{\mathbb{C}}  \label{La inmersion}
\end{equation}%
This morphism is in fact generically injective. With a similar proof, we can
see that there is a morphism%
\begin{equation*}
\mathcal{M}_{\mathbb{C}}(\Gamma _{0}(N_{1}))\longrightarrow \mathcal{M}_{%
\mathbb{C}}(\mathfrak{\mathfrak{a}},N_{2})
\end{equation*}%
inducing%
\begin{equation}
j_{X_{\mathbb{C}}(\Gamma (N))}:X_{\mathbb{C}}(\Gamma (N))\longrightarrow Y_{%
\mathbb{C}}(N_{2})\text{ \ }  \label{complex immersion}
\end{equation}%
for the congruence group%
\begin{equation*}
\Gamma _{1}(N)\subseteq \Gamma (N):=\Gamma _{0}(d_{\mathfrak{K}}N_{1})\cap
\Gamma _{1}(N_{2})\subseteq \Gamma _{0}(N)
\end{equation*}%
because this is the subgroup of $\Gamma _{1}^{1}(\mathfrak{a,}N_{2})$
stabilizing the diagonal\ in its action on $\mathfrak{H}_{\mathfrak{K}}$.
Indeed, a matrix $\gamma =\left( 
\begin{array}{cc}
a & b \\ 
c & d%
\end{array}%
\right) $ of $\Gamma _{1}^{1}(\mathfrak{a,}N_{2})$ stabilizes the diagonal $%
\mathfrak{H}$ if and only if $\gamma =\gamma ^{\prime }$ or $\gamma =-\gamma
^{\prime }$. It cannot be $\gamma =-\gamma $ because then $d$ \ belongs to $%
\mathbb{Z[}\sqrt{d_{\mathfrak{K}}}]$ so it cannot be $1$ $\func{mod}.$ $N_{2}
$. Therefore $\gamma $ is selfconjugate, i.e. a matrix in 
\begin{equation*}
\Gamma _{1}^{1}(\mathfrak{a,}N_{2})\cap Sl(2,\mathbb{Q})=\Gamma (N)
\end{equation*}%
Both quotients \ $\Gamma _{1}^{1}(\mathfrak{a,}N_{2})$ $\setminus \mathfrak{H%
}_{\mathfrak{K}}$ \ and $\Gamma _{0}(N)\setminus \mathfrak{H}_{\mathfrak{K}}$
by the left actions of $\Gamma _{1}^{1}(\mathfrak{a,}N_{2})$ and of the
stabilizer $\Gamma _{0}(N)$ of the diagonal are defined over $\mathbb{Z}[%
\frac{1}{N}]$, thus over $R$.

The morphism described between the smooth varieties of complex points, of
both curve and surface, is in fact induced by a morphism of schemes. Indeed,
in \cite{Kudla y Rapop} it is justified that Hirzebruch-Zagier cycles in
Hilbert modular surfaces are in fact given by morphisms of schemes over $%
\mathbb{Z[}\frac{1}{N}\mathbb{]}_{\text{ }}$, thus over our ring $R\supseteq 
\mathfrak{o}$. The construction is made more generally for\ Hilbert-Zagier
cycles in "twisted" Hilbert modular surfaces, which in the particular case
of a Hilbert modular surface are just the usual Hirzebruch-Zagier cycles.
These are the images of suitable morphisms from Shimura curves associated to
quaternion algebras, and their transforms by the Hecke correspondences. In
our particular case, the morphism $X_{R}(N)\longrightarrow Y_{R}(N_{2})$ is
from a Shimura curve associated to a quaternion algebra 
\begin{equation*}
(\frac{a,b}{\mathbb{Q}})\subseteq M_{2}(\mathbb{Q}(\sqrt{a}))
\end{equation*}%
which is in fact the matrix algebra $M_{2}(\mathbb{Q})$ (so that $a=1$),
because it is the quaternion algebra $Q_{B}$ associated by V. 1. (5) \cite%
{VdG} to the skew-hermitian matrix $B=\left( 
\begin{array}{cc}
0 & 1 \\ 
-1 & 0%
\end{array}%
\right) $ whose $B-$diagonal $\mathcal{H}_{B}\subseteq \mathcal{H}^{2}$ of
equation V.1. (1.3) \cite{VdG} is the true diagonal $z_{1}=z_{2}$ of $%
\mathcal{H}^{2}$ .\ 

If $E_{X_{R}(N)}^{U}$ is the universal family on $X_{R}(N)$ as modular
curve, then\ the pullback\textbf{\ }of the universal family $%
A_{Y_{R}(N_{2})}^{U}$ by $j_{X_{R}(N)}$ is 
\begin{equation}
A_{X_{R}(N_{2})}^{U}:=E_{X_{R}(N)}^{U}\otimes \mathfrak{d}_{\mathfrak{K}%
}^{-1}\approx E_{X_{R}(N)}^{U}\times _{X_{R}(N)}E_{X_{R}(N)}^{U}
\label{A=E x E}
\end{equation}%
where the isomorphism is the relative version of\ the isomorphism $E\times _{%
\mathbb{Z}}\mathfrak{d}_{\mathfrak{K}}^{-1}\mathfrak{\approx }E\times E$ for
any elliptic curve $E$ over $R$ , as $M\otimes _{\mathbb{Z}}\mathfrak{j}%
\approx M\otimes _{\mathbb{Z}}\mathfrak{o\approx }M\oplus M$ for any
fractional ideal $\mathfrak{j}$ as we assume $R\supseteq $ $\mathfrak{o}$
(cf. (5.4) of \ and \cite{Split Katz} 2.1.1 ). The diagonal $%
E_{X_{R}(N)}^{U}\longrightarrow $ $A_{X_{R}(N)}^{U}$ \ induces the
"restriction" morphism between bundles of relative differentials $pr_{\ast
}\Omega _{A_{X_{R}(N)}^{U}/X_{R}(N)}\twoheadrightarrow pr_{\ast }\Omega
_{E_{X_{R}(N)}^{U}/X_{R}(N)}$, i.e. $j_{X_{R}(N)}^{\ast }\mathcal{\mathcal{R}%
}\twoheadrightarrow $ $L_{0,R}$ restricting isomorphisms%
\begin{equation*}
j_{X_{R}(N)}^{\ast }L_{Y_{R}(N_{2})}^{(1,0)}\approx L_{0,R}\text{ \ \ and \ }%
j_{X_{R}(N)}^{\ast }L_{Y_{R}(N_{2})}^{(0,1)}\approx L_{0,R}\text{ , }
\end{equation*}%
so that 
\begin{equation*}
j_{X_{R}(N)}^{\ast }L_{Y_{R}}^{(k,k^{\prime })}=L_{0,R}^{k+k^{\prime }}\text{
.}
\end{equation*}%
\footnote{%
Just to put this in relation with more general contexts, we observe that the
diagonal can be recovered in terms of the $\mathfrak{\mathfrak{o}}_{B}-$%
multiplication in $A_{X_{R}(N)}^{U}$ as kernel of $\left( 
\begin{array}{cc}
-1 & 0 \\ 
0 & 1%
\end{array}%
\right) \in \mathfrak{o}_{B}$. The quotient line bundle $L_{0,R}$ on this
particular Shimura curve still has sense on a Shimura curve on a non-split
quaternionic algebra $Q_{B}$, the global sections of its powers being the
Hilbert modular forms, according to the Atkin-Lehner correspondence. The
notation $L_{0,R}$ for this line bundle on the modular curve $X_{R}(N)$
usually denoted $\underline{\omega }_{X_{R}(N))}$ recalls that $\Gamma _{0}$
is the essential part of the level $\Gamma (N)=\Gamma _{0}(d_{\mathfrak{K}%
}N_{1})\cap \Gamma _{1}(N_{2})$ , as the additional $N_{2}-$level is just to
have universal families;\ our notation avoids confusions, as $\underline{%
\omega }_{X_{R}(N)}^{\otimes 2}(-cusps)$ has a Kodaira-Spencer isomorphism
with the cotangent line bundle $\omega _{X_{R}(N)}$of the modular curve; and
avoids an excessive use of the omega symbol in this article}

\subsection{\protect\bigskip$\ p-$adic Hilbert modular forms}

\bigskip Let $K\supseteq $ $\mathfrak{K}$ be a number field, having ring of
integers $R$ $\supseteq $ $\mathfrak{o}$. Choose $p$\ splitting in $K$\ , so
that it splits in $\mathfrak{K}$ as $p=\pi \pi ^{\prime }$. Then\textit{\ }$%
p $ does not divide the discriminant $d_{\mathfrak{K}}$\textit{\ }, and we
choose it prime with\textit{\ }$N=d_{\mathfrak{K}}N_{1}N_{2}$ .The embedding 
$\mathfrak{K\hookrightarrow \mathfrak{K}}_{\pi }\approx \mathbb{Q}_{p}$
extends to an embedding $K\hookrightarrow \overline{\mathbb{Q}}_{p}$ which
restricts to $R_{\mathfrak{m}}\approx \mathbb{Z}_{p}$ for a maximal ideal $%
\mathfrak{m}$ of $R$ lying on $\pi $, so that $\mathbb{Z}\subseteq \mathfrak{%
o\subseteq }R$ induces isomorphisms%
\begin{equation*}
\text{ }\mathbb{Z}_{p}\approx \mathfrak{o}_{\pi }\approx R_{\mathfrak{m}}%
\text{ and }\mathbb{Q}_{p}\approx \mathfrak{K}_{\pi }\approx K_{\mathfrak{m}}
\end{equation*}%
We get, analogously for $\pi ^{\prime }$,%
\begin{equation*}
\mathbb{Z}_{p}\approx \mathfrak{o}_{\pi ^{\prime }}\approx R_{\mathfrak{m}%
^{\prime }}\text{ and }\mathbb{Q}_{p}\approx \mathfrak{K}_{\pi ^{\prime
}}\approx K_{\mathfrak{m}^{\prime }}
\end{equation*}

Denote $\kappa $ the residual field $\mathbb{F}_{p}$ . Recall from section 2
of \cite{Goren on Hasse}, for instance, that the Hasse invariant $H_{%
\mathfrak{K,}p}$ is a $\mathfrak{a-}$Hilbert modular form over $\kappa $ of
level $N_{2}$ of parallel weight $(p-1,p-1)$ , i.e. a section of $L_{%
\overline{\mathcal{M}}_{\kappa }(N_{2})}^{(p-1,p-1)}$, whose $q-$expansion
is $1$. The open subscheme $\mathcal{M}_{\kappa }(N_{2})^{ord}$ of the
ordinary points is the complement of the support of the ample divisor $D^{h}$
of $\overline{\mathcal{M}_{\kappa }}(N_{2})$ where $H_{\mathfrak{K},p}$
vanishes\footnote{%
Points in the support of $D^{h}$ , i.e. nonordinary for the reduction $%
\mathcal{M}_{R}(N_{2})\times _{\mathbb{Z}}\mathcal{M}_{\mathbb{Z}/p}(N_{2})$
are union of points in the support of two divisors $D_{\pi }$ and $D_{\pi
^{\prime }}$ nonordinary for the reductions $\mathcal{M}_{R}(N_{2})\times _{%
\mathfrak{o}}\mathcal{M}_{\mathfrak{o}/\pi }(N_{2})$ and $\mathcal{M}%
_{R}(N_{2})\times _{\mathfrak{o}}\mathcal{M}_{\mathfrak{o}/\pi ^{\prime
}}(N_{2})$ , the two larger strata of the stratification of $D^{h}$ \ \cite%
{Goren generalizing Katz}considered in def. 2.1 of}. We still denote by $%
D^{h}$ the non-ample divisor which is strict transform of it in $Y_{\kappa
}(N_{2})$

\textbf{Proposition}. There is an integer $n_{0}\geq 0$ , which we take
minimal, such that an $\mathfrak{a-}$Hilbert modular form $H_{\mathbb{Z}%
_{p}} $ section of $L_{\mathcal{M}_{\mathbb{Z}%
_{p}}(N_{2})}^{(n_{0}(p-1),n_{0}(p-1))}$ reduces modulo $p$ to the section $%
H_{\mathfrak{K,}p}^{n_{0}}$ of $L_{\mathcal{M}_{\kappa
}(N_{2})}^{(n_{0}(p-1),n_{0}(p-1))}.$

\textbf{Proof}: From the ampleness of $L_{\overline{\mathcal{M}}_{\mathbb{Z}%
_{p}}}^{(1,1)},$ it follows the exactenes of the sequence, for $n>>0$,%
\begin{equation*}
0\longrightarrow L_{\overline{\mathcal{M}}_{\mathbb{Z}_{p}}(N_{2})}^{(n,n)}%
\overset{p}{\longrightarrow }L_{\overline{\mathcal{M}}_{\mathbb{Z}%
_{p}}(N_{2})}^{(n,n)}\longrightarrow L_{\overline{\mathcal{M}}_{\mathbb{F}%
_{p}}(N_{2})}^{(n,n)}\longrightarrow 0
\end{equation*}

and, from this, the exactness of the sequence%
\begin{eqnarray*}
0 &\longrightarrow &H^{0}(L_{\overline{\mathcal{M}}_{\mathbb{Z}%
_{p}}(N_{2})}^{(n,n)})\overset{p}{\longrightarrow }H^{0}(L_{\overline{%
\mathcal{M}}_{\mathbb{Z}_{p}}(N_{2})}^{(n,n)}) \\
&\longrightarrow &H^{0}(L_{\overline{\mathcal{M}}_{\mathbb{F}%
_{p}}(N_{2})}^{(n,n)})\longrightarrow H^{1}(L_{\overline{\mathcal{M}}_{%
\mathbb{Z}_{p}}(N_{2})}^{(n,n)})=0
\end{eqnarray*}%
$\square $

We still denote by $D^{h}$ the ample divisor of $\overline{\mathcal{M}}_{%
\mathbb{Z}_{p}}(N_{2})$ where $H_{\mathbb{Z}_{p}}$\ vanishes, which is a
normal crossing divisor, and denote the same its strict transform in $Y_{%
\mathbb{Z}_{p}}(N_{2})$, although it is not longer ample as it is disjoint
with the strict transforms of the cusps (The subscripts of the divisors will
not be needed, as is always clear where they belong) The schemes%
\begin{equation*}
Y_{\mathbb{Q}_{p}}^{\prime }(N_{2})=Y_{\mathbb{Q}_{p}}(N_{2})-D^{h}\text{
and }Y_{\kappa }^{\prime }(N_{2})=Y_{\kappa }(N_{2})-D^{h}\text{ ,}
\end{equation*}%
smooth and projective over $\mathbb{Q}_{p}$ and $\kappa $, are not affine
but lie thus on affine varieties 
\begin{equation*}
\overline{\mathcal{M}_{\mathbb{Q}_{p}}(N_{2})}^{\prime }=\overline{\mathcal{M%
}_{\mathbb{Q}_{p}}(N_{2})}^{\prime }-D^{h}\text{ \ and }\overline{\mathcal{M}%
_{\kappa }(N_{2})}^{\prime }=\overline{\mathcal{M}_{\kappa }(N_{2})}-D^{h}
\end{equation*}%
being $\mathcal{M}_{\kappa }^{ord}(N_{2})=$\ $\mathcal{M}_{\kappa
}(N_{2})^{\prime }$ the intersection of this last with $\mathcal{M}_{\kappa
}(N_{2})$ .

The scheme $Y_{\mathbb{Q}_{p}}(N_{2})$ has $Y_{\mathbb{Z}_{p}}(N_{2})$ as a
flat, smooth and projective model over $\mathbb{Z}_{p}$, with special fibre $%
Y_{\kappa }(N_{2})$ (this uses that $N_{2}$ is invertible in $\mathbb{Z}_{p}$%
). We will denote $\mathcal{Y}_{\mathbb{Z}_{p}}(N_{2})$ the formal
completion of $Y_{\mathbb{Z}_{p}}(N_{2})$ along this special fibre, i.e. its 
$p-$adic completion. Because of properness, the rigid points of the formal
scheme $\mathcal{Y}_{\mathbb{Z}_{p}}(N_{2})$ correspond bijectively with the
points of the rigidification $Y_{\mathbb{Q}_{p}}(N_{2})^{rig}$ of the $%
\mathbb{Q}_{p}$ $-$scheme $Y_{\mathbb{Q}_{p}}(N_{2})$ ( cf. 2.4 prop. 3, and
2.7 def.1, prop. 7 and 1.13 prop.4 of \cite{Bosch}). When rigid is clear by
the context -for instance while taking differentials on a rigid analytic
open subset- we just write $Y_{\mathbb{Q}_{p}}(N_{2})$;\ and we will also
denote the same the rigid analytic coherent sheaf $\mathcal{F}^{rig}$ on $Y_{%
\mathbb{Q}_{p}}(N_{2})^{rig}$ corresponding in a unique way to a coherent
sheaf $\mathcal{F}$ on $Y_{\mathbb{Q}_{p}}(N_{2})$ by the rigid version of
the GAGA principle (1.16 theory. 11-13 of Loc.cit), if clear by the context,
for instance while taking sections of it over a rigid open subset.

\bigskip The specialization map 
\begin{equation*}
sp:Y_{\mathbb{Q}_{p}}(N_{2})^{rig}\longrightarrow Y_{\kappa }(N_{2})(\kappa )
\end{equation*}%
induces%
\begin{equation*}
red:Y_{\mathbb{Q}_{p}}(N_{2})(\mathbb{Q}_{p})=Y_{\mathbb{Z}_{p}}(N_{2})(%
\mathbb{Z}_{p})\longrightarrow Y_{\kappa }(N_{2})(\kappa )\text{ \ }
\end{equation*}%
extending to the equally denoted reduction%
\begin{equation*}
\text{ }red:Y_{\mathbb{Q}_{p}}(N_{2})(\mathbb{C}_{p})=Y_{\mathbb{Z}%
_{p}}(N_{2})(\mathcal{O}_{\mathbb{C}_{p}})\longrightarrow Y_{\kappa }(N_{2})(%
\overline{\kappa }).
\end{equation*}

For a locally closed subset $S$ of $Y_{\kappa }(N_{2})(\kappa )$, the tube $%
]S[\subseteq Y_{\mathbb{Q}_{p}}(N_{2})^{rig}$ is the rigid $\mathbb{Q}_{p}$%
-space (denoted equally $]S[$ the corresponding rigid $\mathbb{C}_{p}$%
-space) which is counterimage of $S$ by the specialization map, so that $%
]Y_{\kappa }(N_{2})[=Y_{\mathbb{Q}_{p}}(N_{2})^{rig}$.

\bigskip Since $D^{h}$ is defined over $\mathbb{Z}_{p}$, we obtain a flat
and smooth scheme $Y_{\mathbb{Z}_{p}}^{\prime }(N_{2})=Y_{\mathbb{Z}%
_{p}}(N_{2})-D^{h}$ over $\mathbb{Z}_{p}$ having special fibre $Y_{\kappa
}^{\prime }(N_{2})=Y_{\kappa }(N_{2})-D^{h}$ . Let $\mathcal{A}\subseteq Y_{%
\mathbb{Q}_{p}}^{\prime }(N_{2})^{rig}$ $\subseteq Y_{\mathbb{C}%
_{p}}(N_{2})^{rig}$ be the rigid $\mathbb{Q}_{p}-$space $]Y_{\kappa
}^{\prime }(N_{2})[$ , and denote equally\ by$\mathcal{A}$ the corresponding
rigid $\mathbb{C}_{p}$-space.

\bigskip The Zariski open subset $Y_{\mathbb{Z}_{p}}^{\prime }(N_{2})(%
\mathbb{C}_{p})\subseteq Y_{\mathbb{Z}_{p}}(N_{2})(\mathbb{C}_{p})$ can be
seen, near the divisor $D^{h}$, as the set of $\mathbb{C}_{p}$-points where
a $\mathbb{C}_{p}$-valued function, still denoted $H_{\mathbb{C}_{p}}$ ,
does not vanish (This alternative way of looking at the $\mathfrak{a}$-
Hilbert modular forms is explained for instance in Katz \cite{Split Katz}
1.2.1 to 1.2.8). A local trivialization of the Rapoport rank one locally
free sheaf $\mathcal{R}$ induces a local trivialization of any associated
line bundle, and this allows us to see its local sections as just
functions;\ we use this to trivialize $L_{Y_{\mathbb{Z}%
_{p}}(N_{2})}^{(n_{0}(p-1),n_{0}(p-1))}$near $D^{h}$ . Since the section $H_{%
\mathbb{C}_{p}}$ is defined over $\mathbb{Z}_{p}$ , the equally denoted
function is defined over $\mathcal{O}_{\mathbb{C}_{p}}$). The analogous
holds for $Y_{\kappa }^{\prime }(N_{2})(\overline{\kappa })$. The rigid open
set $\mathcal{A\subseteq }Y_{\mathbb{C}_{p}}(N_{2})^{rig}$ appears as the
set of $\mathbb{C}_{p}-$points whose value by $H_{\mathbb{C}_{p}}$ has $p-$%
adic valuation $0$. It can be enlarged 
\begin{equation}
\mathcal{A}\mathcal{\subseteq \mathcal{W}}_{\varepsilon }\subseteq Y_{%
\mathbb{Z}_{p}}^{\prime }(N_{2})(\mathcal{\mathbb{C}}_{\substack{ p  \\ }})%
\text{ }  \label{wide opens}
\end{equation}%
by the wide open sets $\mathcal{\mathcal{W}}_{\varepsilon }$ (in the rigid
topology of $Y_{\mathbb{Z}_{p}}^{\prime }(N_{2})(\mathbb{C}_{p})$) where the 
$\mathbb{C}_{p}$-value of $H_{\mathbb{C}_{p}}$ has $p-$ adic order strictly
smaller than $\varepsilon >0$. Clearly, 
\begin{equation}
\mathcal{\mathcal{W}}_{\varepsilon }\subseteq \mathcal{W}_{\varepsilon
^{\prime }}\text{ \ for }\varepsilon <\varepsilon ^{\prime }\text{ , and }%
\dbigcup\limits_{\varepsilon >0}\mathcal{W}_{\varepsilon }=Y_{\mathbb{Z}%
_{p}}^{\prime }(N_{2})(\mathbb{C}_{p})\text{ .}
\label{inclusion of wide opens}
\end{equation}%
Denote\ $\overline{\mathcal{M}_{\mathbb{Q}_{p}}(N_{2})}^{ord}$ the
projection of $\mathcal{A}$ to $\overline{\mathcal{M}_{\mathbb{Q}_{p}}(N_{2})%
}$ , and denote $\mathcal{M}_{\mathbb{Q}_{p}}(N_{2})^{ord}$ its intersection
with $\mathcal{M}_{\mathbb{Q}_{p}}(N_{2})$. Analogously, denote by $%
\overline{\mathcal{M}_{\mathbb{Q}_{p}}(N_{2})}(\varepsilon )$ the projection
of $\mathbb{\mathcal{W}}_{\varepsilon }$ to $\overline{\mathcal{M}_{\mathbb{Q%
}_{p}}(N_{2})}$ , and $\mathcal{M}_{\mathbb{Q}_{p}}(N_{2})(\varepsilon )$ ,
or just $\mathcal{M}_{\varepsilon }$ , its intersection with $\mathcal{M}_{%
\mathbb{Q}_{p}}(N_{2})$\textbf{. }The $\mathbb{Z}_{p}$-module of $p-$adic $%
\mathfrak{a}$-Hilbert modular forms of weight $(k,k^{\prime })\in \mathbb{%
Z\times \mathbb{Z}}$ 
\begin{equation}
M_{k,k^{\prime }}^{p-\text{adic}}(N_{2},\mathbb{Z}_{p}):=\Gamma _{rig}(L_{%
\mathcal{A}}^{(k,k^{\prime })}/\mathbb{Z}_{p})\supseteq \Gamma L_{\mathcal{M}%
(N_{2})}^{(k,k^{\prime })}=M_{k,k^{\prime }}(N_{2},\mathbb{Z}_{p})
\label{inclusion}
\end{equation}%
is the one of rigid analytic sections on $\mathcal{A}$ , defined over $%
\mathbb{Z}_{p}$, of $L_{Y_{\mathbb{Z}_{p}}(N_{2})}^{(k,k^{\prime })}$ as a
rigid analytic line bundle. These are $p-$adic in the sense of Katz in \cite%
{Split Katz} 1.9.1 to 1.9.3 , recalled in 11 of \cite{And Go} \footnote{%
This is because abelian schemes over $R_{\mathfrak{m}}/\mathfrak{m}^{n}$
have ordinary $p-$reduction if and only if they admit a $\mu _{p\alpha }$
-level structure $\varepsilon _{p^{\alpha }}$, thus a $\mu _{N_{2}p^{\alpha
}}$-level structure $\varepsilon _{N_{2}p^{\alpha }}$, for some $\alpha $
(cf. 1.11 of \cite{Split Katz}). Since the tangent map of \ $\varepsilon
_{N_{2}p^{\alpha }}:$ $\mu _{N_{2}p^{\alpha }}\otimes _{\mathbb{Z}}\mathfrak{%
d}^{-1}\longrightarrow A_{R_{\mathfrak{m}}/\mathfrak{m}^{n}}$ is an
isomorphism we can push to $A_{R_{\mathfrak{m}}/\mathfrak{m}^{n}}$ the
canonical differential $\omega _{can}=\frac{dt}{t}$of $\mu _{N_{2}p^{\alpha
}}$ , so to obtain a nowhere vanishing section $\omega _{can}$ trivializing
the Rapoport bundle $\mathcal{R}_{\mathcal{M}_{R_{\mathfrak{m}}/\mathfrak{m}%
^{n}}(N_{2}p^{\alpha })}$ and thus the associated modular line bundles. For
those of parallel weight this extends to a trivialization at the cusps, so
that their sections become functions, this giving, as projective limit over $%
\alpha $, a section of the structure sheaf of the formal scheme $\mathcal{M}%
_{R_{\mathfrak{m}}}\mathcal{(}N_{2}p^{\infty })$ , i.e. a $p-$adic Hilbert
modular form in the sense of Katz, and the association is in fact
biunivocal. In the unparallel case, just add the fact that all forms are
then cuspidal.}, and thus have a $q-$expansion 
\begin{equation*}
\dsum\limits_{\nu \in (\mathfrak{a}^{-1})^{+}\cup \{0\}}a_{\nu }q^{\nu }%
\text{ ,}
\end{equation*}%
with all $a_{\nu }\in \mathbb{Z}_{p}$ , which agrees with the $q-$expansion
formerly defined, in case it is classical.

\bigskip

Denote $S_{k,k^{\prime }}^{p-\text{adic}}(N_{2},\mathbb{Z}_{p})$ the $%
\mathbb{Z}_{p}$ -module of the $p-$adic\textbf{\ }Frechet subspace of the
cuspidal ones, i.e. those having $a_{\nu }=0$ for this expansion at the
standard cusp and at all cusps. The $\mathbb{Z}_{p}$-module $M_{k,k^{\prime
}}^{p-\text{adic}}(N_{2},\mathbb{Z}_{p})$ contains the submodule of
overconvergent $\mathfrak{a}$-Hilbert modular forms , i.e. rigid analytic
sections on \textit{some} wide open set. This is the filtered union%
\begin{equation}
M_{k,k^{\prime }}^{oc}(N_{2},\mathbb{Z}_{p},):=\dbigcup\limits_{\varepsilon
>0}\Gamma _{rig}(L_{\mathcal{W}_{\varepsilon }}^{(k,k^{\prime })}/\mathbb{Z}%
_{p})\supseteq S_{k,k^{\prime }}^{oc}(N_{2},\mathbb{Z}_{p})\text{ }
\label{modular contains cuspidal}
\end{equation}%
made by the sections of the sheaf $j^{\dag }L_{\mathcal{M}_{\mathbb{Z}%
_{p}}^{\prime }(N_{2})}^{(k,k^{\prime })}$ where we denote%
\begin{equation}
j^{\dag }\mathcal{F}:=\lim_{\mathcal{\varepsilon }}j_{\ast }(\mathcal{F}%
\lfloor _{\mathbb{\mathcal{M}}_{\varepsilon }})  \label{definicion jcruz}
\end{equation}%
for any coherent sheaf $\mathcal{F}$ on $\mathcal{M}_{\mathbb{Z}%
_{p}}^{\prime }(N_{2})$ (and the notation holds for a coherent sheaf on $Y_{R%
\mathfrak{m}}^{\prime }(N_{2})$ ). The inclusion (\ref{modular contains
cuspidal}) is an equality in the unparallel case. We\textbf{\ }define in an
analogous way $p-$adic and overconvergent $\mathfrak{a}$-Hilbert modular
forms, and cuspidal forms, with coefficients in any complete extension of $%
\mathbb{Z}_{p}$, as for instance $\mathbb{Q}_{p}$ so to obtain $p-$adic\
Frechet spaces $M_{k,k^{\prime }}^{p-adic}(N_{2},\mathbb{Q}_{p},),$ $%
M_{k,k^{\prime }}^{oc}(N_{2},\mathbb{Q}_{p},)$, etc.

An $\mathfrak{a}$-Hilbert modular form $f$ \ of $q-$expansion (\ref%
{expansion}), classical, $p-$adic, or overconvergent, of weight $%
(k,k^{\prime })$ defined over $\mathbb{Z}_{p}$ "restricts" to a classic, $p-$%
adic, or overconvergent modular form $f\lfloor $ of weight $k+k^{\prime }$,
defined over $\mathbb{Z}_{p}$, of $q-$expansion 
\begin{equation}
f\lfloor (q)=\dsum\limits_{\nu \in (\mathfrak{a}^{-1})^{+}\cup
\{0\}}^{{}}a_{\nu }q^{Tr(\nu )}\text{ .}  \label{restriccion de expansion}
\end{equation}%
\qquad\ This is because $\mathcal{A}_{0}:=$ $j_{X(N)_{\mathbb{C}_{p}}}^{-1}%
\mathcal{A}$ and $\mathcal{W}_{\varepsilon ,0}:=j_{X(N)_{\mathbb{C}%
_{p}}}^{-1}\mathcal{W}_{\varepsilon }$ are the affinoid and the wide open
sets of the modular curve $X(N)_{\mathbb{C}_{p}}$ considered in the usual
definition of $p-$adic and overconvergent modular forms.

The $k-$linear Frobenius $\phi :\mathcal{M}_{\kappa }(N_{2})\longrightarrow 
\mathcal{M}_{\kappa }(N_{2})$ is realized on geometric points $x$ of $%
\mathcal{M}_{\kappa }(N_{2})^{ord}$ in the following way:\ \ For an abelian
surface $A_{x}$ with $\mathfrak{o}$-multiplication, ordered $\mathfrak{a-}$%
polarization and $N_{2}-$level, corresponding to a point $x$ , i.e.
admitting a level $\varepsilon _{p^{\alpha }}:\mu _{p^{\alpha }}\otimes _{%
\mathbb{Z}}\mathfrak{d}^{-1}$ $\longrightarrow $ $A_{x}$ , the kernel of the 
$\kappa -$linear Frobenius (cf. lemma 4.23 \ and cor. 4.22 of\textbf{\ \cite%
{Bess detallado}}) $\phi :A_{x}\longrightarrow A_{x}$ is the torsion $p$
group subscheme $H_{x}\approx \mu _{p}\times \mu _{p}$ which is\ $\mathfrak{o%
}$- invariant because the product by $\mathfrak{o}$ in $A_{x}$ preserves the 
$p-$torsion. The ordered $\mathfrak{a-}$polarization on $A_{x}$ is a
symmetric bilinear form which descends to $A_{x}/H_{x}$ as, being the target
torsion free, it kills any torsion subgroup. Furthermore, $\mu
_{N_{2}}\otimes _{\mathbb{Z}}\mathfrak{d}_{\mathfrak{K}}^{-1}\hookrightarrow
A_{x}\twoheadrightarrow A_{x}/H_{x}$ is still a monomorphism as $N_{2}$ is
prime to $p$. We obtain in this way $\phi (x)$\ (cf. \cite{Split Katz}, 1.
11). Among the zero characteristic liftings of $\phi $, all of them inducing
homotopic maps on the complex of differentials (cf. prop. 1 of \cite{Bess
pedag}), we will choose and denote equally a fixed morphism 
\begin{equation*}
\phi :\mathcal{M}_{\mathbb{Q}_{p}}(N_{2})^{ord}\longrightarrow \mathcal{M}_{%
\mathbb{Q}_{p}}(N_{2})^{ord}
\end{equation*}%
extending to 
\begin{equation*}
\phi :\mathcal{M}_{\mathbb{Q}_{p}}(N_{2})(\varepsilon )\longrightarrow 
\mathcal{M}_{\mathbb{Q}_{p}}(N_{2})(\varepsilon ^{\prime })
\end{equation*}%
for some $\varepsilon ^{\prime }>\varepsilon >0$ , which is provided by
lemma 3.1.1 and 3.1.7 of \cite{Kisin and Lai}, or alternatively by theorem.
4.3.1 ii) of \cite{Conrad}, where a finite etale subgroup $G_{1}\subseteq A_{%
\mathcal{M}_{\mathbb{Q}_{p}}(N_{2})(\varepsilon )}^{U}[p]$ of rank $p^{d}$
is found, with $d$ the relative dimension, in our case $d=2$. The morphism $%
\phi $ is then covered by a homomorphism $\Phi $ of abelian surfaces on
rigid spaces

\begin{equation}
\begin{array}{ccc}
A_{\varepsilon }:=A_{\mathcal{M}_{\mathbb{Q}_{p}}(N_{2})(\varepsilon )}^{U}
& \overset{\Phi }{\longrightarrow } & A_{\mathcal{M}_{\mathbb{Q}%
_{p}}(N_{2})(\varepsilon ^{\prime })}^{U}=:A_{\varepsilon ^{\prime }} \\ 
\downarrow &  & \downarrow \\ 
\mathcal{M}_{\varepsilon }:=\mathcal{M}_{\mathbb{Q}_{p}}(N_{2})^{ord}(%
\varepsilon )\text{\ } & \overset{\phi }{\longrightarrow } & \mathcal{M}_{%
\mathbb{Q}_{p}}(N_{2})(\varepsilon ^{\prime })=:\mathcal{M}_{\varepsilon
^{\prime }}\text{\ }%
\end{array}%
\text{ }  \label{diagrama analogo Coleman}
\end{equation}%
factoring by the quotient $A_{\mathcal{M}_{\mathbb{Q}_{p}}(N_{2})(%
\varepsilon )}^{U}/G_{1}$ isomorphic to the pullback $\phi ^{\ast }(A_{%
\mathcal{M}_{\mathbb{Q}_{p}}(N_{2})(\varepsilon ^{\prime })}^{U})$ of \ $A_{%
\mathcal{M}_{\mathbb{Q}_{p}}(N_{2})(\varepsilon ^{\prime })}^{U}$ by $\phi $
, followed by the natural projection from this pullback, just as quoted in
(2.1) of Coleman 's \cite{Coleman95} from sect. 3.10 of Katz' \cite{Katz LNM}%
. The construction of the lifting $\phi $ uses the additional fact,
mentioned in the comment preceding lemma 3.1.1 of \cite{Kisin and Lai}, that
this quotient inherits a natural $\mathfrak{o-}$multiplication, $\mathfrak{a-%
}$polarization and $\mu _{N_{2}}$-level, a fact which can be proved as above.

The diagram (\ref{diagrama analogo Coleman}) analogous to (2.1) in \cite%
{Coleman95} , induces as in Loc. cit. morphisms%
\begin{equation}
\Phi ^{\ast }\Omega _{A_{\varepsilon ^{\prime }}/\mathcal{M}_{\varepsilon
^{\prime }}}\longrightarrow \Omega _{A_{\varepsilon }/\mathcal{M}%
_{\varepsilon }}\text{ \ and \ }\Phi _{\ast }\Omega _{A_{\varepsilon }/%
\mathcal{M}_{\varepsilon }}\longrightarrow \Omega _{A_{\varepsilon ^{\prime
}}/\mathcal{M}_{\varepsilon ^{\prime }}}\text{ \ }
\label{Par morfismo en abeliana}
\end{equation}%
inducing%
\begin{equation}
\phi ^{\ast }\mathcal{R}_{\mathcal{M}_{\varepsilon ^{\prime
}}}\longrightarrow \mathcal{R}_{\mathcal{M}_{\varepsilon }}\text{ \ and \ }%
\phi _{\ast }\mathcal{R}_{\mathcal{M}_{\varepsilon }}\longrightarrow 
\mathcal{R}_{\mathcal{M}_{\varepsilon ^{\prime }}}
\label{Par morfismos en Rapoport}
\end{equation}%
thus%
\begin{equation*}
\phi ^{\ast }L_{\mathcal{M}_{\varepsilon ^{\prime }}}^{(k,k^{\prime
})}\longrightarrow L_{\mathcal{M}_{\varepsilon }}^{(k,k^{\prime })}\text{
and }\phi _{\ast }L_{\mathcal{M}_{\varepsilon }}^{(k,k^{\prime
})}\longrightarrow L_{\mathcal{M}_{\varepsilon ^{\prime }}}^{(k,k^{\prime })}
\end{equation*}%
and so inducing on sections%
\begin{equation*}
\Gamma L_{\mathcal{M}_{\varepsilon ^{\prime }}}^{(k,k^{\prime
})}\longrightarrow \Gamma \phi ^{\ast }L_{\mathcal{M}_{\varepsilon ^{\prime
}}}^{(k,k^{\prime })}\longrightarrow \Gamma L_{\mathcal{M}_{\varepsilon
}}^{(k,k^{\prime })}\text{ and }\Gamma L_{\mathcal{M}_{\varepsilon
}}^{(k,k^{\prime })}\longrightarrow \Gamma L_{\mathcal{M}_{\varepsilon
^{\prime }}}^{(k,k^{\prime })}
\end{equation*}%
or, by the Koecher' principle also valid for rigid analytic sections (cf.
lemma 4.14 of \cite{Kisin and Lai}),%
\begin{equation}
\phi :\Gamma L_{\mathcal{W}_{\varepsilon ^{\prime }}}^{(k,k^{\prime
})}\longrightarrow \Gamma L_{\mathcal{W}_{\varepsilon }}^{(k,k^{\prime })}%
\text{ and }Ver:\Gamma L_{\mathcal{W}_{\varepsilon }}^{(k,k^{\prime
})}\longrightarrow \Gamma L_{\mathcal{W}_{\varepsilon ^{\prime
}}}^{(k,k^{\prime })}  \label{Frob y Ver on overconver}
\end{equation}%
These induce equally denoted endomorphisms of the sheaves $j^{\dag }L_{%
\mathcal{M}_{\mathbb{Z}_{p}}^{\prime }(N_{2})}^{(k,k^{\prime })}$ and $%
j^{\dag }L_{Y_{\mathbb{Z}_{p}}^{\prime }(N_{2})}^{(k,k^{\prime })}$ ,
defined as in \ref{definicion jcruz}, so they induce equally denoted
operators $\phi $ and $Ver$\ on the $\mathbb{Z}_{p}-$module $M_{k,k^{\prime
}}^{oc}(N_{2},\mathbb{Z}_{p})$ of its sections, and on the submodule $%
S_{k,k^{\prime }}^{oc}(N_{2},\mathbb{Z}_{p})$ of the cuspidal (called
Vershiebung operator because it is induced by a morphism which on fibres of
the universal family and mod. $p$ is the dual of the Frobenius isogeny).
Just as in theorem. 5.1 of \cite{Coleman95}, where the arguments using
diagram (2.2)\cite{Coleman95}, analogous to our diagram (\ref{diagrama
analogo Coleman}), do not really depend on the dimension $d$, their
composition is 
\begin{equation*}
Ver\circ \phi =p^{k+k^{\prime }+d}=p^{k+k^{\prime }+2}
\end{equation*}%
and, based in 1.11.21 and lemma 1.11.22 of \cite{Split Katz}), $\phi $ acts
on $q-$expansions as $p^{k+k^{\prime }+2}V_{p}$ where 
\begin{equation*}
V_{p}:\dsum\limits_{\nu \in (\mathfrak{a}^{-1})^{+}\cup \{0\}}^{{}}a_{\nu
}q^{\nu }\mapsto \dsum\limits_{\nu \in (\mathfrak{a}^{-1})^{+}\cup
\{0\}}^{{}}a_{\nu }q^{p\nu }=\dsum\limits_{\substack{ \nu \in (\mathfrak{a}%
^{-1})^{+}\cup \{0\}  \\ p\mid \nu }}^{{}}a_{\nu /p}q^{\nu }
\end{equation*}%
(restriction of the equally denoted "Katz operator" of \cite{And Go}, 13.10,
which shows that this Katz operator preserves overconvergent forms).
Analogously, the Verschiebung operator $Ver$ acts as 
\begin{equation*}
U_{p}:\dsum\limits_{\nu \in (\mathfrak{a}^{-1})^{+}\cup \{0\}}^{{}}a_{\nu
}q^{\nu }\mapsto \dsum\limits_{\substack{ \nu \in (\mathfrak{a}%
^{-1})^{+}\cup \{0\}  \\ p/\nu }}^{{}}a_{\nu }q^{\nu /p}=\dsum\limits_{\nu
\in (\mathfrak{a}^{-1})^{+}\cup \{0\}}^{{}}a_{p\nu }q^{\nu }
\end{equation*}%
being clearly $U_{p}V_{p}=1$.

Since $Ver$ $\ $and $\ \phi $ act on $j^{\dag }L_{Y_{\mathbb{Q}%
_{p}}(N_{2})}^{(k,k^{\prime })}$ preserving cuspidality, so in fact on $%
j^{\dag }\Omega _{Y_{\mathbb{Q}_{p}}(N_{2})}^{i}$ for $i=0,1,2$, and this is
compatible with the exterior derivative%
\begin{equation}
j^{\dag }\mathcal{O}_{]Y_{\kappa }[}\longrightarrow j^{\dag }\Omega
_{]Y_{\kappa }[}^{1}\longrightarrow j^{\dag }\omega _{]Y_{\kappa }[}\text{ }
\label{short complex for rigid}
\end{equation}%
they also act on the hypercohomology of this complex, which is the rigid
cohomology $H_{rig}^{i}(Y_{\kappa _{\mathfrak{m}}}^{\prime }(N_{2}),\mathbb{Q%
}_{p}).$ \ As this $\phi $ -action reduces to the $k-$linear Frobenius of $%
Y_{\kappa }^{\prime }(N_{2}),$ it computes the natural Frobenius action on
this rigid cohomology, because it can be computed equivalently by any zero
characteristic lifting of the $k-$linear Frobenius. The natural restriction
gives a homomorphism of complexes%
\begin{equation*}
\begin{array}{ccccc}
\mathcal{O}_{]Y_{\kappa }[} & \longrightarrow & \Omega _{]Y_{\kappa }[}^{1}
& \longrightarrow & \omega _{]Y_{\kappa }[} \\ 
\downarrow &  & \downarrow &  & \downarrow \\ 
j^{\dag }\mathcal{O}_{]Y_{\kappa }[} & \longrightarrow & j^{\dag }\Omega
_{]Y_{\kappa }[}^{1} & \longrightarrow & j^{\dag }\omega _{]Y_{\kappa }[}%
\end{array}%
\end{equation*}%
inducing a $\mathbb{Q}_{p}-$linear map 
\begin{equation*}
H_{dR}^{i}(Y_{\mathbb{Q}_{p}}(N_{2}),\mathbb{Q}_{p})\approx
H_{rig}^{i}(Y_{\kappa _{\mathfrak{m}}}(N_{2}),\mathbb{Q}_{p})\longrightarrow
H_{rig}^{i}(Y_{\kappa _{\mathfrak{m}}}^{\prime }(N_{2}),\mathbb{Q}_{p})
\end{equation*}%
compatible with the $\phi -$action.

\bigskip Associated to the splitting $p=\pi \pi ^{\prime }$ in $\mathfrak{o}$%
, the $\mathfrak{o-}$morphisms (\ref{Par morfismo en abeliana}) and (\ref%
{Par morfismos en Rapoport}) give morphisms 
\begin{equation*}
\phi ^{\ast }L_{\mathcal{M}_{\varepsilon ^{\prime }}}^{(1,0)}\oplus \phi
^{\ast }L_{\mathcal{M}_{\varepsilon ^{\prime }}}^{(,0,1)}\longrightarrow L_{%
\mathcal{M}_{\varepsilon }}^{(1,0)}\text{ }\oplus L_{\mathcal{M}%
_{\varepsilon }}^{(0,1)}\text{\ and \ }\phi _{\ast }L_{\mathcal{M}%
_{\varepsilon }}^{(1,0)}\oplus \phi _{\ast }L_{\mathcal{M}_{\varepsilon
}}^{(0,1)}\longrightarrow L_{\mathcal{M}_{\varepsilon ^{\prime
}}}^{(1,0)}\oplus L_{\mathcal{M}_{\varepsilon ^{\prime }}}^{(0,1)}
\end{equation*}%
each of them direct sum of the $\mathbb{Z}_{p}$ -morphisms obtained by
extension of (\ref{Par morfismos en Rapoport}) via the two embeddings $%
\mathfrak{\mathfrak{o\hookrightarrow }}R$ and $\mathfrak{o}\overset{\sigma }{%
\mathfrak{\approx }}\mathfrak{\mathfrak{o\hookrightarrow }}R$ which induce
isomorphisms $\mathfrak{o}_{\pi }\approx R_{\mathfrak{m}}$ and\ $\mathfrak{o}%
_{\pi ^{\prime }}\approx R_{\mathfrak{m}^{\prime }}$ \textbf{, }all of them
isomorphic to $\mathbb{Z}_{p}$\textbf{\ . }The direct factors\textbf{\ }$%
\phi ^{\ast }L_{\mathcal{M}_{\varepsilon ^{\prime }}}^{(1,0)}\longrightarrow
L_{\mathcal{M}_{\varepsilon }}^{(1,0)}$ and \ $\phi _{\ast }L_{\mathcal{M}%
_{\varepsilon }}^{(1,0)}\longrightarrow $ $L_{\mathcal{M}_{\varepsilon
^{\prime }}}^{(1,0)}$ induce $\phi _{\pi }$ and $Ver_{\pi }$ on $j^{\dag }L_{%
\mathcal{M}_{\mathbb{Z}_{p}}^{\prime }(N_{2})}^{(k,0)}$ analogous to (\ref%
{Frob y Ver on overconver}). The endomorphism $\phi _{\pi }$ acts on $q-$%
expansions as $p^{k}V_{\pi }$ for%
\begin{equation}
V_{\pi }:\dsum\limits_{\nu \in (\mathfrak{a}^{-1})^{+}\cup \{0\}}^{{}}a_{\nu
}q^{\nu }\mapsto \dsum\limits_{\nu \in (\mathfrak{a}^{-1})^{+}\cup
\{0\}}^{{}}a_{\nu }q^{\pi \nu }=\dsum\limits_{\substack{ \nu \in (\mathfrak{a%
}^{-1})^{+}\cup \{0\}  \\ \pi \mid \nu }}^{{}}a_{\nu /\pi }q^{\nu }\text{ ,}
\label{Vpi}
\end{equation}%
and $Ver_{\pi }$ acts as 
\begin{equation}
U_{\pi }:\dsum\limits_{\nu \in (\mathfrak{a}^{-1})^{+}\cup \{0\}}^{{}}a_{\nu
}q^{\nu }\mapsto \dsum\limits_{\substack{ \nu \in (\mathfrak{a}%
^{-1})^{+}\cup \{0\}  \\ \pi /\nu }}^{{}}a_{\nu }q^{\nu /\pi
}=\dsum\limits_{\nu \in (\mathfrak{a}^{-1})^{+}\cup \{0\}}^{{}}a_{\pi \nu
}q^{\nu }  \label{Upi}
\end{equation}

Endomorphisms $\phi _{\pi ^{\prime }}$ and $Ver_{\pi ^{\prime }}$ of $%
j^{\dag }L_{\mathcal{M}_{\mathbb{Z}_{p}}^{\prime }(N_{2})}^{(0,k^{\prime })}$%
are analogously induced by the second direct factors, and they act on $q-$%
expansions via operators $p^{k^{\prime }}V_{\pi ^{\prime }}$ and $U_{\pi
^{\prime }}$ . We can still denote $\phi _{\pi }$ the endomorphism of 
\begin{equation*}
j^{\dag }L_{\mathcal{M}_{\mathbb{Z}_{p}}^{\prime }(N_{2})}^{(k,k^{\prime
})}=j^{\dag }L_{\mathcal{M}_{\mathbb{Z}_{p}}^{\prime
}(N_{2})}^{(k,0)}\otimes j^{\dag }L_{\mathcal{M}_{\mathbb{Z}_{p}}^{\prime
}(N_{2})}^{(0,k^{\prime })}
\end{equation*}%
obtained as limit of the behavior on sections of 
\begin{equation*}
(\phi ^{\ast }L_{\mathcal{M}_{\varepsilon ^{\prime }}}^{(k,0)})\otimes L_{%
\mathcal{M}_{\varepsilon ^{\prime }}}^{(0,k^{\prime })}\overset{\text{rest. }%
\otimes \text{ id.}}{\longrightarrow }(\phi ^{\ast }L_{\mathcal{M}%
_{\varepsilon }}^{(k,0)})\otimes L_{\mathcal{M}_{\varepsilon ^{\prime
}}}^{(0,k^{\prime })}\overset{\phi _{\pi }\otimes id}{\longrightarrow }L_{%
\mathcal{M}_{\varepsilon ^{\prime }}}^{(k,0)}\otimes L_{\mathcal{M}%
_{\varepsilon ^{\prime }}}^{(0,k^{\prime })}\text{ }
\end{equation*}%
It induces an operator on $M_{k,k^{\prime }}^{oc}(N_{2},\mathbb{Z}_{p})$ and 
$S_{k,k^{\prime }}^{oc}(N_{2},\mathbb{Z}_{p}),$ \footnote{%
Since, according a former footnote, $D^{h}$ is union of the two larger
strata $D_{\pi }$ and $D_{\pi ^{\prime }}$ it can be naturally defined
overconvergence respect $D_{\pi }$ , a property preserved by $U_{\pi }$ and $%
V_{\pi }$, and overconvergence respect $D_{\pi ^{\prime }}$, but we will
deal here with just overconvergence in the whole sense.}whose behavior on $q$%
-expansions is again $p^{k}V_{\pi }$ as in \ref{Vpi} ;\ and analogously with 
$\phi _{\pi ^{\prime }}$acting on $q-$expansions as $\ p^{k^{\prime }}V_{\pi
^{\prime }}$. Ditto for the Verschiebung operators $Ver_{\pi }$ and $%
Ver_{\pi ^{\prime }}$ which act on $q-$expansions as $U_{\pi }$ and $U_{\pi
^{\prime }}$.

Clearly, 
\begin{align*}
U_{\pi }U_{\pi ^{\prime }}& =U_{\pi ^{\prime }}U_{\pi }=U_{p}, \\
V_{\pi }V_{\pi ^{\prime }}& =V_{\pi ^{\prime }}V_{\pi }=V_{p}\text{ so that }%
\phi _{\pi }\phi _{\pi ^{\prime }}=\phi _{\pi ^{\prime }}\phi _{\pi }=\phi
_{p}, \\
U_{\pi }V_{\pi }& =U_{\pi ^{\prime }}V_{\pi ^{\prime }}=U_{p}V_{p}=1.
\end{align*}

\section{\textbf{The }$\mathbf{p-}$\textbf{adic Abel-Jacobi map}}

\subsection{\textbf{Besser theory}}

\bigskip

We recall first the essentials of Besser theory.$\ $There is a$\ p-$adic
version of the Abel-Jacobi map in terms of the syntomic cohomology, but
since it is not possible to mimic integration in this cohomology, this
framework does not allow us to compute this map in terms of $p-$adic
integration. Besser has solved this problem in \cite{Besser inventi} by
embedding the syntomic cohomology into the finite-polynomial cohomology,
where the integration is possible. As well known, $p-$adic integration was
already available for $1$-differentials in Coleman's theory by using
Frobenius as monodromy path to integrate, and Besser has generalized (in
cohomological guise) this use of Frobenius to integrate differentials of
higher order. Let $X_{\mathbb{Q}_{p}}$ be a smooth, irreducible scheme of
dimension $d$ over $\mathbb{Q}_{p}$, having a smooth, flat model $X_{\mathbb{%
Z}_{p}}$over $\mathbb{Z}_{p}$. Denote by $X_{\kappa }$ its reduction to $%
\kappa =\mathbb{F}_{p}.$ In \cite{Besser inventi}, Amnon Besser embeds the
syntomic cohomology of $X_{\mathbb{Z}_{p}}$ into a $\mathbb{Q}_{p}-$space
which he calls "finite polynomial cohomology"%
\begin{equation*}
H_{syn}^{i}(X_{\mathbb{Z}_{p}},n)\subseteq H_{fp}^{i}(X_{\mathbb{Z}_{p}},n)
\end{equation*}%
in a way which depends functorially on $X_{\mathbb{Z}_{p}}$ but is not a
true cohomology theory, in the sense of providing long exact sequences as
expected.\textbf{\ }To construct this wider space, Amnon Besser considers in
(2.1 to 2.4 of \cite{Besser inventi}), for a positive integer $u$ and $P(x)$
in the multiplicative set $\mathfrak{P}_{u}$ of\ rational polynomials\textbf{%
\ }$P(x)=\dprod (1-\alpha _{j}x)$ with first coefficient $1$and roots $%
\alpha _{j}^{-1}$of complex norm $p^{\frac{u}{2}},$ and for any integer $%
n\geq 0$, the syntomic $P-$complex $\mathbb{R}\Gamma _{P}(X_{\mathbb{Z}%
_{p}},n)$ defined as the mapping fibre complex \footnote{%
We recall that the mapping fibre of a homomorphism of complexes $\psi
:A^{\bullet }\longrightarrow B^{\bullet }$ is the complex given in degree $i$
by $A^{i}\oplus B^{i-1}$, with coboundary map 
\begin{equation*}
d(a^{i},b^{i-1})=(da^{i},\psi (a^{i})-db^{i-1})\text{ ,}
\end{equation*}%
so its cocycles are given by a cocycle $a^{i}$ in $A^{i}$ and a preimage $%
b^{i-1}$ in $B^{i-1}$ of $\psi (a)$ .} of the homomorphism of complexes 
\begin{equation}
Fil^{n}\mathbb{R}\Gamma _{dR}(X_{\mathbb{Q}_{p}}/\mathbb{Q}_{p})\overset{%
P(\phi )}{\longrightarrow }\mathbb{R}\Gamma _{rig}(X_{\kappa },\mathbb{Q}%
_{p})
\end{equation}%
composition of 
\begin{equation}
Fil^{n}\mathbb{R}\Gamma _{dR}(X_{\mathbb{Q}_{p}}/\mathbb{Q}%
_{p})\longrightarrow \mathbb{R}\Gamma _{dR}(X_{\mathbb{Q}_{p}},\mathbb{Q}%
_{p})\longrightarrow \mathbb{R}\Gamma _{rig}(X_{\kappa }/\mathbb{Q}_{p})
\label{De Rham to rig}
\end{equation}%
(the last being a map only in the derived category, as it is composition of
a morphism and the inverse of a quasi-isomorphism) with the evaluation of
the polynomial $P(x)$ in the endomorphism 
\begin{equation}
\mathbb{R}\Gamma _{rig}(X_{\kappa }/\mathbb{Q}_{p})\overset{\phi }{%
\longrightarrow }\mathbb{R}\Gamma _{rig}(X_{\kappa }/\mathbb{Q}_{p})\text{ }
\end{equation}%
of the complex of sheaves of $\mathbb{Q}_{p}-$spaces (def. 4.13 of \cite%
{Bess detallado}) induced, according cor. 4.22 and lemma 4.23 of Loc.cit.,
by the $\kappa -$ linear Frobenius $\phi :X_{\kappa }\longrightarrow
X_{\kappa }$. This we equally denote $\phi $ , as well as the $\mathbb{Q}%
_{p}-$linear map%
\begin{equation*}
\phi :H_{rig}^{i}(X_{\kappa }/\mathbb{Q}_{p})\longrightarrow
H_{rig}^{i}(X_{\kappa }/\mathbb{Q}_{p})
\end{equation*}%
which it induces on the finite dimensional $\mathbb{Q}_{p}-$space of rigid
cohomology of $X_{\kappa }$. If $X_{\mathbb{Z}_{p}}$ is proper over $\mathbb{%
Z}_{p}$ so that $X_{\mathbb{Q}_{p}}$ and $X_{\kappa }$ are proper over $%
\mathbb{Q}_{p}$and $\kappa $ resp.,\ the $\mathbb{Q}_{p}$-linear map in
cohomologies%
\begin{equation*}
H_{rig}^{i}(X_{\kappa }/\mathbb{Q}_{p})\longrightarrow H_{dR}^{i}(X_{\mathbb{%
Q}_{p}}/\mathbb{Q}_{p})
\end{equation*}%
induced by the map in the derived category is an isomorphism.

The finite polynomial complex $\mathbb{R}\Gamma _{fp,m}(X_{\mathbb{Z}%
_{p}},n) $ of weight $n$ is defined by Besser (def. 2.4 of \cite{Besser
inventi}) as the direct limit of $\mathbb{R}\Gamma _{P}(X_{\mathbb{Q}%
_{p}},n) $ for all polynomials $P\in \mathfrak{P}_{u}$ ordered by the
divisorial relation, and the "finite polynomial" cohomology $H_{fp,m}^{i}(X_{%
\mathbb{Z}_{p}},n)$ of weight $u$ is the cohomology of this complex. Since
direct limit is an exact functor, $H_{fp,m}^{i}(X_{\mathbb{Z}_{p}},n)$ is
the direct limit of the cohomologies $H_{P}^{i}(X_{\mathbb{Z}_{p}},n)$ of
the complexes $\mathbb{R}\Gamma _{P}(X_{\mathbb{Q}_{p}},n)$. Besser writes $%
H_{fp}^{i}(X_{\mathbb{Z}_{p}},n)$ for $H_{fp,i}^{i}(X_{\mathbb{Z}_{p}},n)$.
As a piece of the long exact sequence

\begin{align*}
... & \longrightarrow Fil^{n}H_{dR}^{i-1}(X_{\mathbb{Q}_{p}}/\mathbb{Q}_{p})%
\overset{P(\phi)}{\longrightarrow}H_{rig}^{i-1}(X_{\kappa}/\mathbb{Q}_{p}) \\
& \longrightarrow H_{P}^{i}(X_{\mathbb{Z}_{p}},n)\longrightarrow\text{$Fil$}%
^{n}H_{dR}^{i}(X/\mathbb{Q}_{p})\overset{P(\phi)}{\longrightarrow}%
H_{rig}^{i}(X_{\kappa}/\mathbb{Q}_{p})\longrightarrow...\text{ ,}
\end{align*}
we take the short exact sequence

\begin{align}
0& \longrightarrow H_{rig}^{i-1}(X_{\kappa }/\mathbb{Q}_{p})/P(\phi
)Fil^{n}H_{dR}^{i-1}(X_{\mathbb{Q}_{p}}/\mathbb{Q}_{p})\overset{}{%
\longrightarrow }H_{P}^{i}(X_{\mathbb{Z}_{p}},n)  \notag \\
& \longrightarrow \text{$Fil$}^{n}H_{dR}^{i}(X_{\mathbb{Q}_{p}}/\mathbb{Q}%
_{p})^{P(\phi )=0}\longrightarrow 0\ \text{\ }
\label{P-cohomology short sequence}
\end{align}%
as in (12) of \cite{Besser inventi}. Assume $X_{\mathbb{Z}_{n}}$ is proper
over $\mathbb{Z}_{p}$. The polynomials $P\in \mathfrak{P}_{i}\subseteq 
\mathbb{Q[}x]$ which are a multiple in $\mathbb{Q}_{p}\mathbb{[}x]$ \ of a
given polynomial $Q(x)=\dprod (1-\alpha _{j}x)\in \mathbb{Q}_{p}\mathbb{[}x]$
with first coefficient $1$ and roots $\alpha _{j}^{-1}$of complex norm $p^{%
\frac{m}{2}}$make a multiplicative subset $\mathfrak{P}_{i}^{Q}$of $%
\mathfrak{P}_{i}$, which is cofinal in the sense that any polynomial in $%
\mathfrak{P}_{i}$ divides one of them, so that $H_{fp}^{i}$ can be seen as $%
\lim_{P\in \mathfrak{P}_{i}^{Q}}H_{P}^{i}$ . The sequence (\ref{P-cohomology
short sequence}) holds in particular for any $P\in \mathfrak{P}_{i}^{Q_{\phi
}}$, where $Q_{\phi }$ is the characteristic polynomial of $\phi $ acting on 
$H_{dR}^{i}(X_{\mathbb{Q}_{p}}/\mathbb{Q}_{p})$ , so we obtain an exact
sequence

\begin{equation}
\begin{array}{c}
0\longrightarrow H_{dR}^{i-1}(X_{\mathbb{Q}_{p}}/\mathbb{Q}%
_{p})/Fil^{n}H_{dR}^{i-1}(X/\mathbb{Q}_{p})\overset{\mathit{\mathbf{i}}}{%
\longrightarrow } \\ 
H_{fp}^{i}(X_{\mathbb{Z}_{p}},n)\overset{\mathbf{p}}{\longrightarrow }\text{$%
Fil$}^{n}H_{dR}^{i}(X_{\mathbb{Q}_{p}}/\mathbb{Q}_{p})\longrightarrow 0%
\end{array}
\label{Besser exact sequence}
\end{equation}%
because this sequence holds for any $P\in \mathfrak{P}_{i}^{Q_{\phi }}$ with 
$\mathbf{i=}P(\phi )$ . Indeed, properness implies $H_{rig}^{i-1}\approx
H_{dR}^{i-1}$ and $P(\phi )$ acts invertible on it because $P\in \mathfrak{P}%
_{i}$ has eigenvalues of complex norm $p^{\frac{i}{2}}$while $\phi $ acts on
this space with eigenvalues of complex norm $p^{\frac{i-1}{2}}$ , thus $%
P(\phi )Fil^{n}H_{dR}^{i-1}(X_{\mathbb{Q}_{p}}/\mathbb{Q}_{p})\approx
Fil^{n}H_{dR}^{i-1}(X_{\mathbb{Q}_{p}}/\mathbb{Q}_{p})$. As a consequence%
\textbf{,} $H_{fp}^{i}(X_{\mathbb{Z}_{p}},n)=0$ if $i>2d+1$,\ and there is
an isomorphism 
\begin{equation}
H_{dR}^{2d}(X_{\mathbb{Q}_{p}}/\mathbb{Q}_{p})\overset{\mathbf{i}}{%
\longrightarrow }H_{fp}^{2d+1}(X_{\mathbb{Z}_{p}},d+1)
\label{isom dR and fp}
\end{equation}%
since both $Fil^{d+1}H_{dR}^{2d+1}(X_{\mathbb{Q}_{p}}/\mathbb{Q}_{p})$ and $%
Fil^{d+1}H_{dR}^{2d}(X_{\mathbb{Q}_{p}}/\mathbb{Q}_{p})$ are null. Composing
its inverse with the trace map in de Rham cohomology it is obtained in prop.
2.5. 4 of \cite{Besser inventi} the "trace map" in finite polynomial
cohomology%
\begin{equation*}
tr_{X_{\mathbb{Z}_{p}}}:H_{fp}^{2d+1}(X_{\mathbb{Z}_{p}},d+1)\overset{%
\mathbf{i}^{-1}}{\longrightarrow }H_{dR}^{2d}(X_{\mathbb{Q}_{p}}/\mathbb{Q}%
_{p})\overset{tr_{X}}{\longrightarrow }\mathbb{Q}_{p}\text{ .}
\end{equation*}

In various respects, finite polynomial cohomology behaves as de Rham
cohomology: by prop. 2.5. 3 of \cite{Besser inventi} there is a bilinear
cup-product 
\begin{equation*}
H_{fp}^{i}(X_{\mathbb{Z}_{p}},n)\otimes H_{fp}^{j}(X_{\mathbb{Z}_{p}},m)%
\overset{\cup }{\longrightarrow }H_{fp}^{i+j}(X_{\mathbb{Z}_{p}},n+m)
\end{equation*}%
which, composed with the trace map, gives a "Poincar\'{e} type" perfect
pairing%
\begin{equation}
<\_,\_>_{fp}:H_{fp}^{i}(X_{\mathbb{Z}_{p}},n)\otimes H_{fp}^{2d-i+1}(X_{%
\mathbb{Z}_{p}},d-n+1)\overset{tr_{X_{\mathbb{Z}_{p}}}}{\longrightarrow }%
\mathbb{Q}_{p}  \label{cup product in fp cohomology}
\end{equation}%
The commutativity of the diagram

\begin{tabular}{lllll}
$H_{fp}^{i}(X_{\mathbb{Z}_{p}},n)\otimes H_{dR}^{2d-i}(X_{\mathbb{Q}%
_{p}})/Fil^{d-n+1}$ & $\overset{1\otimes \mathbf{i}}{\longrightarrow }$ & $%
H_{fp}^{i}(X_{\mathbb{Z}_{p}},n)\otimes H_{fp}^{2d-i+1}(X_{\mathbb{Z}%
_{p}},d-n+1)$ &  &  \\ 
$\mathbf{p}\otimes 1\downarrow $ &  & $\downarrow <\_,\_>_{fp}$ &  &  \\ 
$Fil^{n}H_{dR}^{i}(X_{\mathbb{Q}_{p}})\otimes H_{dR}^{2d-i}(X_{\mathbb{Q}%
_{p}})/Fil^{d-n+1}$ & $\overset{<\_,\_>_{dR}}{\longrightarrow }$ & $\mathbb{Q%
}_{p}$ &  &  \\ 
&  &  &  & 
\end{tabular}

\bigskip

in \cite{Besser inventi} (14) provides an equivalent definition of the
perfect pairing (\ref{cup product in fp cohomology}) in terms of the
analogous pairing $<\_,\_>_{dR}$ in de Rham cohomology, since $1\times 
\mathbf{i}$ is an isomorphism, so the finite polynomial pairing inherits the
usual properties of it, as the fact that of being equivalently defined by
reduction to the diagonal. This means that $<\_,\_>_{fp}$ coincides with the
composition 
\begin{align*}
H_{fp}^{i}(X_{\mathbb{Z}_{p}},n)\otimes H_{fp}^{2d-i+1}(X_{\mathbb{Z}%
_{p}},d-n+1)& \hookrightarrow H_{fp}^{2d+1}(X_{\mathbb{Z}_{p}}\times X_{%
\mathbb{Z}_{p}},d+1) \\
& \overset{}{\longrightarrow }H_{fp}^{2d+1}(X_{\mathbb{Z}_{p}},d+1)\overset{%
tr_{X_{\mathbb{Z}_{p}}}}{\longrightarrow }\mathbb{Q}_{p}
\end{align*}

\bigskip

i.e. for $\widetilde{\alpha}\in H_{fp}^{i}(X_{\mathbb{Z}_{p}},n)$ and $%
\widetilde{\beta}\in H_{fp}^{2d-i+1}(X_{\mathbb{Z}_{p}},d+1-n)$, it is 
\begin{equation}
<\widetilde{\alpha},\widetilde{\beta}>_{fp}=tr_{X_{\mathbb{Z}_{p}}}(j_{X_{%
\mathbb{Z}_{p}},X_{\mathbb{Z}_{p}}\times X_{\mathbb{Z}_{p}}}^{\ast }(%
\widetilde{\alpha}\otimes\widetilde{\beta}))  \label{casi mi formula}
\end{equation}

\bigskip where $j_{X_{\mathbb{Z}_{p}},X_{\mathbb{Z}_{p}}\times X_{\mathbb{Z}%
_{p}}}:X_{\mathbb{Z}_{p}}\hookrightarrow X_{\mathbb{Z}_{p}}\times X_{\mathbb{%
Z}_{p}}$ is the diagonal map.

\bigskip\textbf{Lemma}.\label{El primer lema} For $j_{X_{\mathbb{Z}_{p}},Y_{%
\mathbb{Z}_{p}}}:X_{\mathbb{Z}_{p}}\hookrightarrow Y_{\mathbb{Z}_{p}}$, and
finite-polynomial classes $\widetilde{\alpha}\in H_{fp}^{i}(X_{\mathbb{Z}%
_{p}},n)$ and $\widetilde{\beta}\in H_{fp}^{2d-i+1}(Y_{\mathbb{Z}%
_{p}},d+1-n) $ on $X_{\mathbb{Z}_{p}}$ and $Y_{\mathbb{Z}_{p}}$, it is 
\begin{equation}
<\widetilde{\alpha},j_{X_{\mathbb{Z}_{p}},Y_{\mathbb{Z}_{p}}}^{\ast }%
\widetilde{\beta}>_{fp}=tr_{X_{\mathbb{Z}_{p}}}\text{ }j_{X_{\mathbb{Z}%
_{p}},X_{\mathbb{Z}_{p}}\times Y_{\mathbb{Z}_{p}}}^{\ast}(\widetilde{\alpha }%
\otimes\widetilde{\beta})\text{ }  \label{mi formula}
\end{equation}

\bigskip \textbf{Proof}:\ both are equal to $tr_{X_{\mathbb{Z}_{p}}}j_{X_{%
\mathbb{Z}_{p},X_{\mathbb{Z}_{p}}\times X_{\mathbb{Z}_{p}}}}^{\ast }(%
\widetilde{\alpha }\otimes j_{X_{\mathbb{Z}_{p}},Y_{\mathbb{Z}_{p}}}^{\ast }%
\widetilde{\beta }),$ the left hand side because of (\ref{casi mi formula}),
and the other because $j_{X_{\mathbb{Z}_{p}},X_{\mathbb{Z}_{p}}\times Y_{%
\mathbb{Z}_{p}}}^{\ast }(\widetilde{\alpha }\otimes \widetilde{\beta })=$ $%
j_{X_{\mathbb{Z}_{p},X_{\mathbb{Z}_{p}}\times X_{\mathbb{Z}_{p}}}}^{\ast }(%
\widetilde{\alpha }\otimes j_{X_{\mathbb{Z}_{p}},Y_{\mathbb{Z}_{p}}}^{\ast }%
\widetilde{\beta })$.$\ \square $

\textbf{Lemma}.\label{segundo lema} For $X_{\mathbb{Z}_{p}}$ of relative
dimension 1, finite-polynomial classes $\widetilde{\alpha},\widetilde{\beta }%
\in$ $H_{fp}^{2}(X_{\mathbb{Z}_{p}},2)$, and a point $o$ of $X$ defined over 
$\mathbb{Z}_{p}$, it is 
\begin{equation*}
j_{(X\times o)_{\mathbb{Z}_{p}}}^{\ast}(\widetilde{\alpha}\otimes \widetilde{%
\beta})=j_{(X\times o)_{\mathbb{Z}_{p}}}^{\ast}(\widetilde{\alpha }\otimes%
\widetilde{\beta})=0
\end{equation*}

(in the notation of a restriction morphism $j$ we drop the bigger scheme,
when clear).

\bigskip\textbf{Proof}:\ It is, for instance, $j_{(X\times o)_{\mathbb{Z}%
_{p}}}^{\ast}(\widetilde{\alpha}\otimes\widetilde{\beta})=$ $j_{X_{\mathbb{Z}%
_{p}}}^{\ast}\widetilde{\alpha}\otimes$ $j_{o_{\mathbb{Z}_{p}}}^{\ast }%
\widetilde{\beta}$ with $j_{o_{\mathbb{Z}_{p}}}^{\ast}\widetilde{\beta}$
belonging to $H_{fp}^{2}(o_{\mathbb{Z}_{p}},2)=0$ .$\square$

Besser theorem gives a "class map" from Chow groups to finite polynomial
cohomology $cl_{fp\text{ }}$compatible with the familiar class map $cl_{dR}$
to de Rham cohomology: 
\begin{equation*}
\begin{array}{ccc}
CH^{i}(X_{\mathbb{Z}_{p}}) & \overset{cl_{fp}}{\longrightarrow } & 
H_{fp}^{2i}(X_{\mathbb{Z}_{p}},i) \\ 
\downarrow  &  & \downarrow  \\ 
CH^{i}(X_{\mathbb{Q}_{p}}) & \overset{cl_{dR}}{\longrightarrow } & \text{$Fil
$}^{2i}H_{dR}^{2i}(X_{\mathbb{Q}_{p}})%
\end{array}%
\end{equation*}

Let us denote by $CH^{i}(X_{\mathbb{Z}_{p}})_{0}\longrightarrow CH^{i}(X_{%
\mathbb{Q}_{p}})_{0}$ the subgroups of classes applying to zero, say of
"null-homologous" classes. Let $Z_{\mathbb{Z}_{p}}=\dsum%
\limits^{{}}r_{l}Z_{l,\mathbb{Z}_{p}}$ represent a class in $CH^{i}(X_{%
\mathbb{Z}_{p}})_{0}$ , with the coefficients $r_{l}$ being rational and the
schemes $Z_{l,\mathbb{Z}_{p}}$ being proper and smooth of dimension $d-i$
over $\mathbb{Z}_{p}$. Let $Z=\dsum\limits^{{}}r_{l}Z_{l}$ be the
corresponding class in $CH^{i}(X_{\mathbb{Q}_{p}})_{0}$, obtained by taking
generic points. \ The $p-$adic or syntomic Abel Jacobi map evaluated at $Z$
is computed in \cite{Besser inventi} theory. 1.2 (and taken here as a
definition) as the functional $Fil^{d-i+1}H_{dR}^{2d-2i+1}(X_{\mathbb{Q}%
_{p}}/\mathbb{Q}_{p})\longrightarrow \mathbb{Q}_{p}$ applying the class of $%
\omega $ into%
\begin{equation}
AJ_{p}(Z)(\omega )=\dint_{Z_{\mathbb{Z}_{p}}}\omega :=\dsum\limits^{{}}r_{l}%
\text{ }tr_{Z_{l,\mathbb{Z}_{p}}}(j_{l}^{\ast }\widetilde{\omega })\in 
\mathbb{Q}_{p}  \label{Abel-Jacobi}
\end{equation}%
where $\widetilde{\omega }\in H_{fp}^{2d-2i+1}(X_{\mathbb{Z}_{p}},d-i+1)$ is
any lifting of $\omega $ by the Besser epimorphism (\ref{Besser exact
sequence}), and $j_{l}$ is the morphism $j_{Z_{l,\mathbb{Z}_{p}},X_{\mathbb{Z%
}_{p}}}:Z_{l,\mathbb{Z}_{p}}\longrightarrow X_{\mathbb{Z}_{p}}$ which Besser
assumes injective, but it is enough generically injective.

Finally, let us recall for the complement $X_{\mathbb{Z}_{p}}^{\prime }$ of
a subscheme $Z_{\mathbb{Z}_{p}}$ also proper over $\mathbb{Z}_{p}$, the
Mayer-Viatoris exact sequence

\begin{align}
...& \longrightarrow H_{rig,Z_{\kappa }}^{i}(X_{\kappa }/\mathbb{Q}%
_{p})\longrightarrow H_{rig}^{i}(X_{\kappa }/\mathbb{Q}_{p})\longrightarrow
H_{rig}^{i}(X_{\kappa }^{\prime }/\mathbb{Q}_{p})  \notag \\
& \longrightarrow H_{rig,Z_{\kappa }}^{i+1}(X_{\kappa }/\mathbb{Q}%
_{p})\longrightarrow ...\text{ \ \ \ \ }
\end{align}%
of $\mathbb{Q}_{p}-$spaces acted by the $\mathbb{Q}_{p}$-linear maps $\phi $
induced from the $\kappa $ $-$linear Frobenius $\phi :X_{\kappa
}\longrightarrow X_{\kappa }$ restricting $\phi :Z_{\kappa }\longrightarrow
Z_{\kappa }.$ Berthelot has proved in \cite{Berthelot} (cor. 5.7) a purity
theorem for $Z_{\mathbb{Z}_{\mathfrak{p}}}$ with $n$ irreducible components
of pure codimension $d$ 
\begin{equation}
H_{rig,Z_{\kappa }}^{2d-1}(X_{\kappa }/\mathbb{Q}_{p})=0\text{ \ and }%
H_{rig,Z_{\kappa }}^{2d}(X_{\kappa }/\mathbb{Q}_{p})=\mathbb{Q}_{p}(-d)^{n}%
\text{ ,}
\end{equation}%
where the Tate's twist notation means that the natural action of $\phi $ on
this copy of $\mathbb{Q}_{p}\approx \mathbb{Q}_{p}$ is given by
multiplication by $p^{d}$ (cf. lemma 7.2 in \ \cite{Besser inventi} and
preceding comment). As a consequence, we have, in case $d=1$, the following
inclusions $\imath ^{1},\imath ^{2}$ (followed of what we call residue maps) 
\begin{equation*}
\begin{array}{c}
0=H_{rig,Z_{\kappa }}^{1}(X_{\kappa }/\mathbb{Q}_{p})\longrightarrow
H_{rig}^{1}(X_{\kappa }/\mathbb{Q}_{p})\overset{\imath ^{1}}{\longrightarrow 
}H_{rig}^{1}(X_{\kappa }^{\prime }/\mathbb{Q}_{p}) \\ 
\overset{res^{1}}{\longrightarrow }H_{rig,Z_{\kappa }}^{2}(X_{\kappa }/%
\mathbb{Q}_{p})\approx \mathbb{Q}_{p}(-1)^{n}\overset{}{\longrightarrow }%
H_{rig}^{2}(X_{\kappa }/\mathbb{Q}_{p}) \\ 
\overset{\imath ^{2}}{\longrightarrow }H_{rig}^{2}(X_{\kappa }^{\prime }/%
\mathbb{Q}_{p})\text{ }\overset{res^{2}}{\longrightarrow }H_{rig,Z_{\kappa
}}^{3}(X_{\kappa }/\mathbb{Q}_{p}).%
\end{array}%
\end{equation*}%
This allows to present a class in $H_{rig}^{i}(X_{\kappa }/\mathbb{Q}_{p})$
as a class in $H_{rig}^{i}(X_{\kappa }^{\prime }/\mathbb{Q}_{p})$ having
null residue, i.e. null image in $H_{rig,Z_{\kappa }}^{i+1}(X_{\kappa }/%
\mathbb{Q}_{p})$, a presentation which is unique in case $i=1$.\ The
advantage of this presentation is that the cohomology $H_{rig}^{i}(X_{\kappa
}^{\prime }/\mathbb{Q}_{p})$ has a particularly handable description in
terms of a particular inclusion $X_{\kappa }^{\prime }\hookrightarrow
X_{\kappa }$ in the special fibre of a proper smooth scheme $X_{\mathbb{Z}%
_{p}}$ over $\mathbb{Z}_{p}$, as follows. First recall that wide open sets,
and thus overconvergence, can be mimicked in this general setting by the
notion of strict neighborhood: Write $X_{\mathbb{Q}_{p}}^{rig}$ for the
rigid analytic space associated to the formal completion of the generic
fibre of $X_{\mathbb{Z}_{p}}$ at the special fibre $X_{\kappa }$ (or $p-$%
completion), and consider the specialization map $sp:$ $X_{\mathbb{Q}%
_{p}}^{rig}\longrightarrow X_{\kappa }(\kappa )$ from this rigid analytic
space to the set of $\kappa -$points of the special fibre;\ for any locally
closed subset $S$ of this special fibre, associate a rigid analytic $\mathbb{%
Q}_{p}-$subspace $\left] S\right[ \subseteq X_{\mathbb{Q}_{p}}^{rig}$ whose
underlying set $sp^{-1}(S)$ consists of all points specializing to some
point in $S$ (cf.\cite{Bess detallado}, 4). A neighborhood $U$ of $\left]
X_{\kappa }^{\prime }\right[ $ in the rigid analytic\textbf{\ }topology of $%
\left] X_{\kappa }\right[ $ is said strict if $\{U,\left] Z_{\kappa }\right[
\}$ is a covering of $\left] X_{k}\right[ $ . Generalizing a notation we
have already used in our particular case, any rigid coherent sheaf $\mathcal{%
F}$ on the rigid analytic space $\left] X_{\kappa }\right[ $ has associated
the sheaf $j^{\dag }\mathcal{F}$ on $\left] X_{\kappa }\right[ $ of
overconvergent sections of $\mathcal{F}$ , i.e. sections over \textit{some}
strict neighborhood $U$ of $\left] X_{\kappa }^{\prime }\right[ $ , i.e. the
direct limit $j^{\dag }\mathcal{F=}\lim_{\longrightarrow }j_{U\ast }(%
\mathcal{F}\lfloor _{U})$ taken over all strict neighborhoods $U$ . In these
notations, the cohomology of the rigid analytic complex $\mathbb{R}\Gamma (%
\left] X_{\kappa }\right[ ,j^{\dag }\Omega _{\left] X_{\kappa }\right[
}^{\bullet })$ is $H_{rig}^{i}(X_{\kappa }^{\prime }/\mathbb{Q}_{p})$ (we
take advantage of the fact that the inclusion of $X_{\kappa }^{\prime }$ in
the special fibre $X_{\kappa }$of $X_{\mathbb{Z}_{p}},$ provides us, taking
formal completion along the special fibre, with a particular "rigid datum"
among those considered in sec. 4 of \cite{Bess detallado}, all the analogous
"rigid data" providing complexes which are mutually quasi-isomorphic, and
quasi-isomorphic to $\mathbb{R}\Gamma _{rig}(X_{\kappa }^{\prime }/\mathbb{Q}%
_{p}),$ i.e. have $H_{rig}^{i}(X_{\kappa }^{\prime }/\mathbb{Q}_{p})$ as
cohomology). If, furthermore, we are giving a morphism $\left] X_{\kappa
}^{\prime }\right[ \longrightarrow $ $\left] X_{\kappa }^{\prime }\right[ $
reducing to $\phi :$ $X_{\kappa }^{\prime }\longrightarrow X_{\kappa
}^{\prime }$ and extending to a morphism $U\longrightarrow V$ between strict
neighborhoods $U\subseteq $ $V$ \ (as will happen in the case we are
interested) then the natural action of this lifting on $H_{rig}^{i}(X_{%
\kappa }^{\prime }/\mathbb{Q}_{p})$ is just the $\phi $-action, so the
lifting provides an operative way to handle this action.

\bigskip

\subsection{\textbf{Application of Besser theory}}

We keep $N_{2},N_{1}$ and $N=d_{\mathfrak{K}}N_{1}N_{2}$ as in our
hypothesis. Recall that the modular curve $X(N)$ has a model smooth, flat
and projective over $\mathbb{Z}[\frac{1}{N}]$, thus over $\mathbb{Q}$ and
over $\mathbb{Z}_{p}$. Choose a point $o$\ of\textit{\ }$X(N)$ with this
definition. The fact that $Y(N_{2})$ has null odd homology (\cite{VdG}, IV,
6.1) implies that the modified diagonal\ $X(N)_{o}$ of

\textit{\ }%
\begin{equation*}
X(N)\times X(N)\overset{}{\longrightarrow }Y(N_{2})\times X(N)\text{ }
\end{equation*}%
is null-homologous null in $Y(N_{2})\times X(N)$ (we have adopted the
convention that omission of the subindex means $X_{\mathbb{Q}_{p}}(N)$ and $%
Y_{\mathbb{Q}_{p}}(N_{2})$ ). This is a consequence of the topological
observation made at the end of the introduction, applied to the complex
curve $X(N)(\mathbb{C})$ and obtained by base change (assumed the choice
made once and for all of an embedding of $\mathbb{Q}_{p}$ into $\mathbb{C}$%
). We will deal with the $p-$adic Abel-Jacobi map of $Y(N_{2})\times X(N)$
at the null-homologous null cycle $X(N)_{o}$, which we write sometimes $%
Y\times X$ and $X_{o}$ for short$.$

$\bigskip $The Hecke operators $T_{\mathfrak{n}\text{ }}$ act on $\mathfrak{%
a-}$Hilbert modular forms for all ideals $\mathfrak{n\subseteq o}$ because
of the hypothesis we have made that $\mathfrak{u}^{+}=(\mathfrak{u}^{+})^{2}$
, since Hecke operators act on the $\mathfrak{a-}$Hilbert modular forms of
group $\Gamma _{1}^{1}(\mathfrak{a,}N_{2})$ invariant by $\mathfrak{u}^{+}/(%
\mathfrak{u}^{+})^{2}$ , i.e. on those which are modular for the bigger
group $\Gamma _{1}(\mathfrak{a,}N_{2})$ ), and in this case both groups are
the same. In fact it acts null on $\mathfrak{a-}$Hilbert modular forms
unless $\mathfrak{n}$ is in the narrow class of $\mathfrak{a}$, i.e. $%
\mathfrak{n=\nu \mathfrak{a}}$ for some $\nu \in \mathfrak{o}^{+}$defined up
to totally positive units, and we denote $T_{\nu }$ this Hecke operator
acting on $\mathfrak{a-}$Hilbert modular forms $f$ . If $f$ is in fact an
eigenform, the coefficient $a_{\nu }$ of its\textbf{\ }$q-$expansion is
precisely the eigenvalue of $T_{\nu }$. We are going to make an explicit use
in this article of two of such operators, namely

\begin{equation}
T_{\pi }:=U_{\pi }+p^{k-1}V_{\pi }\text{ \ and }T_{\pi ^{\prime }}:=U_{\pi
}+p^{k-1}V_{\pi }  \label{Hecke operators}
\end{equation}%
(as follows from the computation of the Fourier coefficients of the
transforms for instance in VI. 1 of \cite{VdG}), so we can take here such
expressions \ref{Hecke operators}as definition of these particular
operators, as well as of the Hecke operator $T_{p}=U_{p}+p^{k-1}V_{p}$
acting on modular forms of weight $k$ , i.e. on sections of the line bundle $%
L_{0}^{k}$ on the modular curve $X_{\mathbb{Q}_{p}}(N)$.

\begin{description}
\item 
\begin{definition}
\label{definicion slope}\textbf{\ }In the notations above, for nonnegative
rationals $\sigma ,\sigma ^{\prime }$ , we denote $S_{k,k^{\prime }}(%
\mathfrak{a,}N_{2},R)^{\sigma ,\sigma ^{\prime }}\subseteq S_{k,k^{\prime }}(%
\mathfrak{a,}N_{2},R)$ the $R-$submodule of forms of slope $(\sigma ,\sigma
^{\prime })$, i.e. eigenforms of $T_{\pi }$ and $T_{\pi ^{\prime }}$ of
eigenvalues $a_{\pi }$ and $a_{\pi ^{\prime }}$ in $R$ such that $ord(a_{\pi
})=\sigma $ in $R_{\mathfrak{m}}\approx \mathfrak{o}_{\pi }\approx \mathbb{Z}%
_{p}$ and $ord(a_{\pi ^{\prime }})=\sigma ^{\prime }$ in $R_{\mathfrak{m}%
^{\prime }}\approx \mathfrak{o}_{\pi ^{\prime }}\approx \mathbb{Z}_{p}$ .
Denote

\begin{equation*}
S_{k,k^{\prime }}(\mathfrak{a,}N_{2},R)^{\sigma ,\ast
}=\dbigcup\limits_{\sigma ^{\prime }\geq 0}S_{k,k^{\prime }}(\mathfrak{a,}%
N_{2},R)^{\sigma ,\sigma ^{\prime }}\text{ }
\end{equation*}%
and analogously $S_{k,k^{\prime }}(\mathfrak{a,}N_{2},R)^{\ast ,\sigma
^{\prime }}$ . In the case of $\sigma =0$ we call such forms ordinary for $%
T_{\pi }$ , otherwise nonordinary, and analogously for $T_{\pi ^{\prime }}.$
\end{definition}
\end{description}

Let $f$\ be an $\mathfrak{a}$- Hilbert cuspidal eigenform of weight $\mathit{%
(2,2)}$ and level $N_{2}\geq 4$ defined over $K$, thus over $K_{\mathfrak{m}%
}\approx \mathbb{Q}_{p}$, to which corresponds a section 
\begin{align*}
\omega _{f}& \in H^{0}(L_{Y_{\mathbb{Q}_{p}}(N_{2})}^{(2,2)}(-D^{c}))=H^{0}(%
\omega _{Y}) \\
& \subseteq H_{dR}^{2}(Y_{\mathbb{Q}_{p}}/\mathbb{Q}_{p})\approx
H_{rig}^{2}(Y_{\kappa }/\mathbb{Q}_{p})\overset{i^{2}}{\hookrightarrow }%
H_{rig}^{2}(Y_{\kappa }^{\prime }/\mathbb{Q}_{p})\text{ }
\end{align*}%
in the smallest filter $Fil^{2}H_{dR}^{2}(Y_{\mathbb{Q}_{p}}(N_{2})/\mathbb{Q%
}_{p})$ of de Rham cohomology, since a section of the line bundle $L_{%
\mathcal{M}_{\mathbb{Q}_{p}}(N_{2})}^{(2,2)}$ extends, by normality, to a
section of the extended line bundle $L_{\overline{\mathcal{M}_{\mathbb{Q}%
_{p}}(N_{2})}}^{(2,2)}$ . Let

\begin{equation*}
H_{dR}^{2}(Y_{\mathbb{C}_{p}}(N_{2})/\mathbb{C}_{p})(f)\subseteq
H_{dR}^{2}(Y_{\mathbb{C}_{p}}(N_{2})/\mathbb{C}_{p})
\end{equation*}%
be the isotypic component of $f$ in cohomology, i.e. the subspace where $%
\omega _{f}$ lies, irreducible invariant for the natural action of $\phi $
on de Rham cohomology. \textbf{\ }Consider the algebraic values $\gamma _{j}$
inverse of the roots of the characteristic polynomial 
\begin{equation}
Q_{f}(x)=\dprod (1-\gamma _{j}x)=\det (1-\phi ^{-1}x\mid H_{dR}^{2}(Y_{%
\mathbb{Q}_{p}}(N_{2})/\mathbb{Q}_{p})(f))\in \mathbb{Q}_{p}[x]
\label{polinomio minimo}
\end{equation}%
with $Q_{f}(\phi )=0$ , all such roots $\gamma _{j}^{-1}$of complex norm $p$
. We will use Besser theory with the cofinal multiplicative system $%
\mathfrak{P}_{2}^{P_{f}}\subseteq \mathfrak{P}_{2}$\ \ for%
\begin{equation}
P_{f}(x):=(1-p^{-1}\phi )Q_{f}(x)\in \mathbb{Q}_{p}[x]
\label{Los dos polinomios de f}
\end{equation}

$\bigskip $Consider the embedding $j_{Y_{\kappa }^{\prime
}(N_{2})}:Y_{\kappa }^{\prime }(N_{2})\hookrightarrow Y_{\kappa }(N_{2})$,
with $Y_{\kappa }(N_{2})-Y_{\kappa }^{\prime }(N_{2})=D^{h_{\kappa }}$, in
which the rigid analytic $\mathbb{Q}_{p}-$space $\left] Y_{\kappa }^{\prime }%
\right[ $ is $\mathcal{A}$ . For the morphism solving singularities%
\begin{equation}
pr:Y^{\prime }\twoheadrightarrow \overline{\mathcal{M}_{\mathbb{Q}%
_{p}}(N_{2})}-D^{h}\text{,}  \label{projection}
\end{equation}%
it is $pr_{\ast }\mathcal{O}_{Y^{\prime }}=\mathcal{O}_{\overline{\mathcal{M}%
_{\mathbb{Q}_{p}}(N_{2})}-D^{h}}$, so that, in the spectral sequence
abutting to the rigid cohomology of $H_{rig}^{\bullet }(Y_{\kappa }^{\prime
}/\mathbb{Q}_{p})$, 
\begin{equation}
E_{1}^{1,0}=H_{rig}^{1}(Y^{\prime },\mathcal{O}_{Y^{\prime }})=H_{rig}^{1}(%
\overline{\mathcal{M}_{\mathbb{Q}_{p}}(N_{2})}-D^{h},\mathcal{O}_{\overline{%
\mathcal{M}_{\mathbb{Q}_{p}}(N_{2})}-D^{h}})=0\text{ ,}
\label{anulacion en sucesion espectral}
\end{equation}%
because the higher cohomology of a rigid analytic coherent sheaf on the
rigid analytification of an affine scheme is null (Kiehl's rigid analytic
version of Cartan's theorem" ). Therefore $E_{2}^{1,0}=0$, and thus the
elements of%
\begin{eqnarray}
Fil^{2}H_{rig}^{2}(Y_{\kappa }^{\prime }/\mathbb{Q}_{p}) &=&\frac{E_{2}^{2,0}%
}{\func{Im}[E_{2}^{1,0}\overset{d_{2}}{\longrightarrow }E_{2}^{0,2}]}%
=E_{2}^{2,0}=  \label{Fil^2 H^2} \\
\frac{E_{1}^{2,0}}{\func{Im}[E_{1}^{1,0}\overset{d_{1}}{\longrightarrow }%
E_{1}^{2,0}]} &=&\frac{\Gamma _{rig}(j_{Y_{\kappa }^{\prime }}^{\dag }\omega
_{\left] Y_{\kappa }\right[ })}{d(\Gamma _{rig}(j_{Y_{\kappa }^{\prime
}}^{\dag }\Omega _{\left] Y_{\kappa }\right[ }))}  \notag
\end{eqnarray}%
are just classes of overconvergent 2-differentials module exterior
derivatives of overconvergent differentials.

Let $P(x)\in \mathfrak{P}_{2}^{P_{f}}$ . We prefer to see, for what follows,
the de Rham class of $\omega _{f}$ as a class in $Fil^{2}H_{rig}^{2}(Y_{%
\kappa }^{\prime }/\mathbb{Q}_{p})$ thanks to the vertical monomorphisms%
\begin{equation}
\begin{array}{ccc}
Fil_{dR}^{2}H_{dR}^{2}(Y_{\mathbb{Q}_{p}}/\mathbb{Q}_{p}) & \approx & 
Fil^{2}H_{rig}^{2}(Y_{\kappa }/\mathbb{Q}_{p}) \\ 
\downarrow \imath _{dR}^{2} &  & \downarrow \imath _{rig}^{2} \\ 
Fil_{dR}^{2}H_{dR}^{2}(Y_{\mathbb{Q}_{p}}^{\prime }/\mathbb{Q}_{p}) & 
\longrightarrow & Fil^{2}H_{rig}^{2}(Y_{\kappa }^{\prime }/\mathbb{Q}_{p})%
\end{array}
\label{cuadrado magico}
\end{equation}%
The class of $P(\phi )($ $\omega _{f})\in \Gamma _{rig}(\left] Y_{\kappa }%
\right[ ,j_{Y_{\kappa }^{\prime }}^{\dag }\omega _{\left] Y_{\kappa }\right[
}^{\bullet })$ in this space (\ref{Fil^2 H^2}) vanishes, so by (\ref{Fil^2
H^2}) there is a one-form$\ \varrho $ on some wide open set $\mathcal{W}%
_{\varepsilon }$ such that 
\begin{equation}
d\varrho =P(\phi )\omega _{f}  \label{primitiva}
\end{equation}%
(we always equally denote a differential and its restriction to an open
set). Then 
\begin{equation}
\widetilde{\omega }_{f}=(\omega _{f},\varrho )\in H_{P}^{2}(Y_{\mathbb{Z}%
_{p}}^{\prime },2)  \label{lift de omega}
\end{equation}%
is a lift of $\omega _{f}$ by the epimorphism 
\begin{equation}
H_{P}^{2}(Y_{\mathbb{Z}_{p}}^{\prime },2)\longrightarrow \text{$Fil$}%
^{2}H_{dR}^{2}(Y^{\prime }/\mathbb{Q}_{p})^{P(\phi )=0}\longrightarrow 0%
\text{\ ,}  \label{Besser epimorphsim}
\end{equation}%
particular case of ($\ref{P-cohomology short sequence}$). Applying to the
lift $\widetilde{\omega }_{f}$ the "restriction" homomorphism%
\begin{equation*}
j_{X^{\prime }(N)}^{\ast }\text{: }H_{P}^{2}(Y_{\mathbb{Z}_{p}}^{\prime
},2)\longrightarrow H_{P}^{2}(X_{\mathbb{Z}_{p}}^{\prime },2)\text{,}
\end{equation*}%
we obtain an element in the $H_{P}^{2}$ cohomology, thus in the $H_{fp}^{2}$
cohomology, of

\begin{equation*}
X_{\mathbb{Z}_{p}}^{\prime }:=j_{X_{\mathbb{Z}_{p}}(N)}^{-1}(Y_{\mathbb{Z}%
_{p}}^{\prime })\text{ }
\end{equation*}%
which provides in turn an element $P(\phi )^{-1}j_{X^{\prime }(N)}^{\ast
}\varrho $ of $H_{rig}^{1}(X^{\prime }{}_{\kappa }/\mathbb{Q}_{p})$ via the
isomorphism

\begin{equation}
0\longrightarrow H_{rig}^{1}(X_{\kappa }^{\prime }/\mathbb{Q}_{p})\overset{%
P(\phi )}{\longrightarrow }H_{fp}^{2}(X^{\prime }{}_{\mathbb{Z}%
_{p}},2)\longrightarrow Fil^{2}H_{dR}^{2}(X_{\mathbb{Q}_{p}}^{\prime }/%
\mathbb{Q}_{p})=0\text{ ,}  \label{isom in X' of dR and fp}
\end{equation}%
Since $P(\phi )$ is a multiple of $(1-p^{-1}\phi ),$ this element is in fact
in the kernel%
\begin{equation*}
H_{dR}^{1}(X_{\mathbb{Q}_{p}}/\mathbb{Q}_{p})\approx H_{rig}^{1}(X_{\kappa }/%
\mathbb{Q}_{p})\overset{\imath ^{1}}{\hookrightarrow }H_{rig}^{1}(X_{\kappa
}^{\prime }/\mathbb{Q}_{p})
\end{equation*}%
of the Berthelot residue map, because $(1-p^{-1}\phi )$ vanishes on the
target $\mathbb{Q}_{p}(-1)$ of this map;\ i.e. it is an element in the first
space of the sequence (cf. $\ref{Besser exact sequence}$) 
\begin{equation}
0\longrightarrow H_{dR}^{1}(X_{\mathbb{Q}_{p}}/\mathbb{Q}_{p})\overset{%
P(\phi )}{\longrightarrow }H_{fp}^{2}(X_{\mathbb{Z}_{p}},2)\longrightarrow
Fil^{2}H_{dR}^{2}(X/\mathbb{Q}_{p})=0\text{ .}  \label{iso in X of dR and fp}
\end{equation}

$\bigskip$

\begin{lemma}
\label{lema computo por Poincare} In these notations, and for any $1-$form 
\begin{equation*}
\eta \in H_{dR}^{1}(X_{\mathbb{Q}_{p}}(N)/\mathbb{Q}_{p})
\end{equation*}%
and polynomial\textbf{\ }$P(x)\in \mathfrak{P}_{2}^{P_{f}}$ , the value on $%
\omega _{f}\otimes \eta $\textbf{\ }of the\ $p$-adic Abel-Jacobi map of $Y_{%
\mathbb{Q}_{p}}(N_{2})\times X_{\mathbb{Q}_{p}}(N)$ at the null-homologous
cycle $X(N)_{o}$ is given by the product in $H_{dR}^{1}(X_{\mathbb{Q}%
_{p}}(N)/\mathbb{Q}_{p})$%
\begin{equation}
AJ_{p}(X(N)_{o})(\omega _{f}\otimes \eta )=<P(\phi )^{-1}j_{X^{\prime
}(N)}^{\ast }\varrho ,\eta >  \label{the cup-product}
\end{equation}
\end{lemma}

\textbf{Proof:}\ \ First, observe that the expression on the right hand of ($%
\ref{the cup-product})$ is the same as for any other de Rham 1-cochain $%
\sigma $ in $Y_{\mathbb{Q}_{p}}^{\prime }(N_{2})$ applying to $P(\phi
)\omega _{f}$ . Indeed, the difference $\varrho -\sigma $ would be a de Rham
cocycle, thus defining a de Rham 1-class in $Y_{\mathbb{Q}_{p}}^{\prime
}(N_{2})$, in fact image of a de Rham class in $Y_{\mathbb{Q}_{p}}(N_{2})$
because $P(\phi )$ vanishes on the target $\mathbb{Q}_{p}(-1)$ of the
residue map . By a base change argument, the fact that $b_{1}(Y_{\mathbb{Q}%
_{p}}(N_{2})(\mathbb{C}))=0$ entails that this is in fact the null class, so
that $<P(\phi )^{-1}j_{X^{\prime }(N)}^{\ast }(\varrho -\sigma ),\eta >=0$.

Denote by $\widetilde{\eta }\in $ $H_{fp}^{2}(X(N)_{\mathbb{Z}_{p}},2)$ the
class in finite-polynomial cohomology corresponding to $\eta $ by the
isomorphism ($\ref{iso in X of dR and fp}$), given by multiplication by $%
P(\phi )$. Let $\widetilde{\omega }_{f}=(\omega _{f},\sigma )$ be a lifting
of $\omega _{f}$ by the Besser epimorphism $H_{fp}^{2}(Y(N_{2})_{\mathbb{Z}%
_{p}},2)\twoheadrightarrow $$H_{dR}^{2}(X_{\mathbb{Q}_{p}}/\mathbb{Q}_{p})$.

By (\ref{Abel-Jacobi}), and lemma (\ref{El primer lema}),

\begin{align}
AJ_{p}(X(N)_{o})(\omega _{f}\otimes \eta )& =\dint_{X(N)_{o}}\omega
_{f}\otimes \eta =  \label{sustracting terms} \\
& tr_{X(N)_{\mathbb{Z}_{p}}}j_{X(N)_{\mathbb{Z}_{p}},Y(N_{2})_{\mathbb{Z}%
_{p}}\times X(N)_{\mathbb{Z}_{p}}}^{\ast }(\widetilde{\omega }_{f}\otimes 
\widetilde{\eta })  \notag \\
& -tr_{(X(N)\times o)_{\mathbb{Z}_{p}}}j_{(X(N)\times o)_{\mathbb{Z}%
_{p}},Y(N_{2})_{\mathbb{Z}_{p}}\times X(N)_{\mathbb{Z}_{p}}}^{\ast }(%
\widetilde{\omega }_{f}\otimes \widetilde{\eta })  \notag \\
& -tr_{(o\times X(N))_{\mathbb{Z}_{p}}}j_{(o\times X(N))_{\mathbb{Z}%
_{p}},Y(N_{2})_{\mathbb{Z}_{p}}\times X(N)_{\mathbb{Z}_{p}}}^{\ast }(%
\widetilde{\omega }_{f}\otimes \widetilde{\eta })  \notag \\
& =tr_{X(N)_{\mathbb{Z}_{p}}}j_{X(N)_{\mathbb{Z}_{p}},Y(N_{2})_{\mathbb{Z}%
_{p}}\times X(N)_{\mathbb{Z}_{p}}}^{\ast }(\widetilde{\omega }_{f}\otimes 
\widetilde{\eta })
\end{align}

The vanishing of the first substracting terms is because 
\begin{align*}
& j_{(X(N)\times o)_{\mathbb{Z}_{p}},Y(N_{2})_{\mathbb{Z}_{p}}\times X(N)_{%
\mathbb{Z}_{p}}}^{\ast}(\widetilde{\omega}_{f}\otimes\widetilde{\eta}) \\
& =j_{(X(N)\times o)_{\mathbb{Z}_{p}},X(N)_{\mathbb{Z}_{p}}\times X(N)_{%
\mathbb{Z}_{p}}}^{\ast}(j_{X(N)_{\mathbb{Z}_{p}},Y(N_{2})_{\mathbb{Z}%
_{p}}}^{\ast}(\widetilde{\omega}_{f})\otimes\widetilde{\eta})
\end{align*}

so that we can apply lemma (\ref{segundo lema}), and analogously with the
vanishing of other substracting term.

By\ (\ref{mi formula}), this can be rewritten as the pairing in finite
polynomial cohomology 
\begin{equation*}
AJ_{p}(X(N)_{o})(\omega _{f}\otimes \eta )=<j_{X(N)}{}^{\ast }\widetilde{%
\omega }_{f},\widetilde{\eta }>_{fp}
\end{equation*}%
which by \cite{Besser inventi} (14) equals the paring of their
isomorphically corresponding classes in first de Rham cohomology (cf. \ref%
{iso in X of dR and fp}). 
\begin{equation*}
AJ_{p}(X(N)_{o})(\omega _{f}\otimes \eta )=<P(\phi )^{-1}j_{X(N)}^{\ast
}\sigma ,\eta >=<P(\phi )^{-1}j_{X^{\prime }(N)}^{\ast }\varrho ,\eta >\text{%
,}
\end{equation*}%
the last equality following from the remark made at the beginning of the
proof \footnote{%
The idea here is that, being the product independent of choice of cochain,
the advantage of $\sigma $ is that it is a cochain on the whole $Y_{K%
\mathfrak{m}}(N_{2})$ , so it serves to compute the $p$-adic Abel-Jacobi
map, and the cochain $\varrho $ has the advantage of being a differential.},
and from the observation that\textbf{\ } $<P(\phi )^{-1}j_{X^{\prime
}(N)}^{\ast }\varrho ,\eta >$ does not depend on the choice of $P(x)$ in $%
\mathfrak{P}_{2}^{P_{f}}$ . Indeed, for any other $Q(x)\in \mathfrak{P}%
_{2}^{P_{f}}$ , there are $P^{\prime },Q^{\prime }\in \mathbb{Q}_{p}\mathbb{[%
}x\mathbb{]}$ \ so that $P^{\prime }P=Q^{\prime }Q\in \mathfrak{P}%
_{2}^{P_{f}}$ , and the equality $d(P^{\prime }(\phi )\varrho )=P^{\prime
}P(\phi )\omega _{f}$ provides a lifting of the class of $\omega _{f}$ in
the Besser epimorphism ($\ref{Besser epimorphsim}$) relative to the
polynomial $P^{\prime }P$, "restricting" the de Rham class $P^{\prime }(\phi
)j_{X_{\mathbb{Q}_{p}}^{\prime }(N)}^{\ast }\varrho $ so that

\begin{equation*}
<(P^{\prime}P(\phi))^{-1}(P^{\prime}(\phi)j_{X_{\mathbb{Q}%
_{p}}^{\prime}(N)}^{\ast}\varrho),\eta>=<P(\phi)^{-1}j_{X_{\mathbb{Q}%
_{p}}^{\prime}(N)}^{\ast}\varrho,\eta>\text{.}
\end{equation*}

\bigskip and analogously for $Q$ and $Q^{\prime}$ $\square$

The evaluation in $\phi $ and $\omega _{f}$ 
\begin{equation*}
f^{\sharp }:=Q_{f}(p^{2}V_{p})(f)\text{ \ and \ }f^{\flat
}:=P_{f}(p^{2}V_{p})(f)
\end{equation*}%
of the two polynomials ($\ref{Los dos polinomios de f}$) provides
overconvergent 2-differentials%
\begin{equation*}
\omega _{f^{\sharp }}:=Q_{f}(\phi )\omega _{f}\text{ \ and }\omega
_{f^{\flat }}=P_{f}(\phi )\omega _{f}=(1-p^{-1}\phi )\omega _{f^{\sharp }}%
\text{ ,}
\end{equation*}%
Choose, once and for all, overconvergent $1-$differentials $\varrho
_{F^{\sharp }}$ and $\varrho _{F^{\flat }}=(1-p^{-1}\phi )\varrho ^{\sharp }$
, corresponding to pairs of $\mathfrak{a-}$Hilbert modular forms $F^{\sharp
}=(F_{(2,0)}^{\sharp },F_{(0,2)}^{\sharp })$ and $F^{\flat
}=(F_{(2,0)}^{\flat },F_{(0.2)}^{\flat })$, such that 
\begin{equation*}
d\varrho _{F^{\flat }}=\omega _{f^{\flat }}\text{ , i.e. }\theta ^{\prime
}F_{(2,0)}^{\flat }-\theta F_{(0,2)}^{\flat }=f^{\flat }
\end{equation*}%
or equivalently 
\begin{equation*}
d\varrho _{F^{\sharp }}=\omega _{f^{\sharp }}\text{ , i.e. }\theta ^{\prime
}F_{(2,0)}^{\sharp }-\theta F_{(0,2)}^{\sharp }=f^{\sharp }
\end{equation*}

We will deal, in the framework of the next section, with a primitive $%
F^{\sharp }=(F_{(2,0)}^{\sharp },0)$ which we will denote $F_{[2,0]}^{\sharp
}$, and analogous $F_{[0,2]}^{\sharp },F_{[2,0]}^{\flat },F_{[0,2]}^{\flat }$
. Lemma (\ref{lema computo por Poincare}) allows us to write 
\begin{equation}
AJ(X(N)_{o})(\omega _{f}\otimes \eta )=<P_{f}(\phi )^{-1}j_{X_{\mathbb{Q}%
_{p}}^{\prime }(N)}^{\ast }\varrho _{F^{\flat }},\eta >\text{,}
\label{Abel-Jacobi en dim 1}
\end{equation}%
a product in $H_{dR}^{1}(X_{\mathbb{Q}_{p}}(N)/\mathbb{Q}_{p})$.

\textbf{\bigskip}

\section{\protect\bigskip \textbf{Differentiation of }$\mathbf{p-}$\textbf{%
adic Hilbert Modular Forms}}

\subsection{\protect\bigskip\textbf{Gauss-Manin covariant derivative}}

Consider, in the variety $Y(N_{2})$, smooth and projective over $\mathbb{Q}%
_{p}$, the rank two locally free $\mathfrak{o\otimes \mathcal{O}}%
_{Y(N_{2})}- $sheaf 
\begin{equation*}
\mathcal{L}_{Y(N_{2})}=\mathcal{H}_{dR}^{1}(\widehat{A}^{U}):=\mathbb{R}%
^{1}pr_{Y(N_{2})}{}_{\ast }\Omega _{\widehat{A}^{U}/Y(N_{2})}^{\bullet
}(\log D_{\widehat{A}^{U}}^{c})
\end{equation*}%
\textbf{(}cf. (\ref{crisis}) and 1.0.15\textbf{\ \cite{Split Katz}}) of
first log-de Rham cohomology of a smooth compactification $\widehat{A}^{U}$%
of the universal family $A^{U}$ 
\begin{equation*}
pr_{Y(N_{2})}:\widehat{A}^{U}\longrightarrow Y(N_{2})\text{ .}
\end{equation*}%
where $D_{\widehat{A^{U}}}^{c}$ denotes the divisor which is counterimage of
the cuspidal divisor $D^{c}$ \ by this projection. The Frobenius action
preserves its natural filtration (cf. \cite{Mladen y Tilouine} 5.2 iii) by
rank one $\mathfrak{o\otimes \mathcal{O}}_{Y(N_{2})}$ -sheaves \qquad 
\begin{equation}
0\longrightarrow \mathcal{R}_{Y(N_{2})}\longrightarrow \mathcal{L}_{Y(N_{2})}%
\mathcal{\longrightarrow \mathcal{S}}_{Y(N_{2})}\longrightarrow 0\text{ .}
\label{extension buena de relative DR}
\end{equation}%
where $\mathcal{\mathcal{S}}_{Y(N_{2})}:=\mathcal{R}_{Y(N_{2})}^{\vee
}\otimes \mathcal{\mathfrak{d}}^{-1}\mathfrak{a}$ is isomorphic as $\mathcal{%
O}_{Y(N_{2})}-$sheaf to

\begin{equation*}
\text{ }\mathcal{R^{\vee }}_{Y(N_{2})}\approx L_{Y(N_{2})}^{(-1,0)}\oplus
L_{Y(N_{2})}^{(0,-1)}
\end{equation*}%
(cf. \cite{Mladen y Tilouine} 1.9 and \cite{Split Katz} 1.0.3 and 1.0.13).
For brevity, we adopt from now on the convention of omitting the subindex $%
Y(N_{2})$\ in the notation of vector bundles over\textit{\ }$Y(N_{2})$. Let 
\begin{equation*}
\text{ }\mathcal{L=\mathcal{L}}^{(1,0)}\oplus \mathcal{L}^{(0,1)}
\end{equation*}%
be the decomposition of\ the rank $2$ locally free $\mathfrak{o\otimes 
\mathcal{O}}_{Y(N_{2})}$- sheaf $\mathcal{L}$ as direct sum of two rank 2
locally free $\mathfrak{\mathcal{O}}_{Y(N_{2})}-$ sheaves. Using \cite{Split
Katz} 2.0.10, we can split (\ref{extension buena de relative DR}), viewed as
a short exact sequence of locally free $\mathfrak{\mathcal{O}}_{Y(N_{2})}$%
-sheaves, as a direct sum of the two exact sequences%
\begin{align*}
0& \longrightarrow L^{(1,0)}\longrightarrow \mathcal{\mathcal{L}}%
^{(1,0)}\longrightarrow L^{(-1,0)}\longrightarrow 0 \\
0& \longrightarrow L^{(0,1)}\longrightarrow \mathcal{\mathcal{L}}%
^{(0,1)}\longrightarrow L^{(0,-1)}\longrightarrow 0
\end{align*}%
Each of these two short exact sequences "restricts" via $j_{X(N)}:X(N)%
\longrightarrow Y(N_{2})$ to the well known sequence 
\begin{equation}
0\longrightarrow L_{0}\longrightarrow \mathcal{L}_{0}\longrightarrow
L_{0}^{-1}\longrightarrow 0  \label{main short sequence in dim 1}
\end{equation}%
on the compact modular curve $X(N)$, where, as in section 2, we denote $%
L_{0} $ the (weight $1$ ) modular line bundle on $X(N)$ and $\mathcal{L}_{0}$
the rank two bundle of (relative, log-cusps) first de Rham cohomology of the
universal generalized elliptic curve for the group $\Gamma (N)=\Gamma
_{0}(d_{\mathfrak{K}}N_{1})\cap \Gamma _{1}(N_{2})$ . Clearly $%
j_{X(N)}^{\ast }\mathcal{L\approx \mathcal{L}}_{0}\oplus \mathcal{L}_{0}$ ,
since, compactifying \ref{A=E x E}), the pull back of $\widehat{A}%
_{Y(N_{2})}^{U}$ by $j_{X(N)}$ is 
\begin{equation}
\widehat{A}_{Y(N_{2})}^{U}\lfloor _{X(N)}\approx \widehat{E}%
_{X(N)}^{U}\times _{X(N))}\widehat{E}_{X(N)}^{U}
\label{decomp of universal compact restr}
\end{equation}%
so this isomorphism is the relative version of the K\"{u}neth decomposition $%
H_{dR}^{1}(A)\approx H_{dR}^{1}(E)\oplus H_{dR}^{1}(E)$ for $A=E\times E,$
with $E$ a generalized elliptic curve .

We now take $n-$th symmetric powers%
\begin{equation*}
\mathcal{R}^{n}\subseteq \mathcal{R}^{\otimes n}\text{ \ \ \ and \ \ \ \ }%
\mathcal{L}^{n}\subseteq \mathcal{L}^{\otimes n}
\end{equation*}%
in the short sequence of locally free $\mathcal{\mathcal{O}}_{Y}-$sheaves 
\begin{equation*}
\text{ }0\longrightarrow \mathcal{R}\approx L^{(1,0)}\oplus
L^{(0,1)}\longrightarrow \mathcal{L\longrightarrow S}\approx
L^{(-1,0)}\oplus L^{(0,-1)}\longrightarrow 0,
\end{equation*}%
\textit{\ } namely (\cite{Split Katz}, 2.1.5)%
\begin{equation}
0\longrightarrow \mathcal{R}^{n}\approx \dbigoplus\limits_{\substack{ %
k+k^{\prime }=n  \\ k,k^{\prime }\geq 0}}^{{}}L^{(k,k^{\prime
})}\longrightarrow \mathcal{L}^{n}\longrightarrow \mathcal{S}%
^{(n)}\longrightarrow 0\text{ ,}  \label{nonsplit symmetric power}
\end{equation}%
where the notation $\mathcal{S}^{(n)}$ used for the quotient is just to
avoid confusion with the symmetric power $\mathcal{S}^{n}$ of $\mathcal{%
\mathcal{S}}$ . This short sequence (\ref{nonsplit symmetric power}) on $%
Y(N_{2})$ is the analogous of the filtration 
\begin{equation}
0\longrightarrow L_{0}^{n}\longrightarrow \mathcal{L}_{0}^{n}\longrightarrow
S_{0}^{(n)}\longrightarrow 0  \label{powered short sequence in dim 1}
\end{equation}%
\qquad of the rank 2 bundle $\mathcal{L}_{0}^{n}$ on the smooth compactified
modular curve $X(N)$, by the rank one bundle $L_{0}^{n}$ and a rank $n$
bundle denoted $S_{0}^{(n)}$ as quotient (and we warn the reader that $%
\mathcal{L}_{0}^{n}\subseteq \mathcal{L}_{0}^{\otimes n}$ is denoted $%
\mathcal{L}_{n}\subseteq \mathcal{L}^{n}$ in \cite{DR} ).

\bigskip

By Ex. II. 5.1.c of \cite{AG}, there is a filtration of locally free $%
\mathcal{\mathcal{O}}_{Y}-$sheaves 
\begin{equation*}
\mathcal{L}^{n}=\mathcal{F}^{0}\supseteq \mathcal{F}^{1}\supseteq
...\supseteq \mathcal{F}^{n}\supseteq \mathcal{F}^{n+1}=0
\end{equation*}%
by bundles $\mathcal{F}^{m}$ image of 
\begin{equation*}
\mathcal{R}^{m}\otimes \mathcal{L}^{n-m}\longrightarrow \mathcal{L}%
^{m}\otimes \mathcal{L}^{n-m}\twoheadrightarrow \mathcal{L}^{n}
\end{equation*}%
with quotients%
\begin{equation}
\mathcal{F}^{m}/\mathcal{F}^{m+1}\approx \mathcal{R}^{m}\otimes \mathcal{S}%
^{n-m}\approx \dbigoplus\limits_{\substack{ k+k^{\prime }=m  \\ k,k^{\prime
}\geq 0}}^{{}}L^{(k,k^{\prime })}\otimes \dbigoplus\limits_{\substack{ %
l+l^{\prime }=n-m  \\ l,l^{\prime }\geq 0}}^{{}}L^{(-l,-l^{\prime })}\approx
\dbigoplus\limits_{{}}^{{}}L^{(k-l,k^{\prime }-l^{\prime })}
\label{filtration}
\end{equation}%
so that $\mathcal{R}^{n}=\mathcal{F}^{n}$ .

\bigskip

\begin{proposition}
\textbf{\bigskip\ \label{rigid Koecher} }The bundle $\mathcal{L}^{n}$
satisfies the Koecher principle 
\begin{equation*}
\Gamma (U,\mathcal{L}^{n})=\Gamma (U\cap \mathcal{M}(N_{2}),\mathcal{L}^{n})
\end{equation*}%
for all open sets $U\subseteq Y(N_{2})$, i.e. $\mathcal{L}%
_{Y(N_{2})}^{n}=inc_{\ast }\mathcal{L}_{\mathcal{M}(N_{2})}^{n}$ for the
inclusion $inc:\mathcal{M}(N_{2})\longrightarrow Y(N_{2})$ .
\end{proposition}

\textbf{Proof:} We prove that all $\mathcal{F}^{m}$ satisfy this principle,
in particular $\mathcal{F}^{0}=\mathcal{L}^{n}$, proceeding by induction. If 
$\mathcal{F}^{m+1}$ does, then as $\mathcal{F}^{m}/\mathcal{F}^{m+1}$
clearly does because of (\ref{filtration}) , and as a consequence $\mathcal{F%
}^{m}$ satisfies also this principle. Indeed, in the diagram 
\begin{equation*}
\begin{array}{ccccc}
\Gamma (U,\mathcal{F}^{m+1}) & \hookrightarrow & \Gamma (U,\mathcal{F}^{m})
& \longrightarrow & \Gamma (U,\mathcal{F}^{m}/\mathcal{F}^{m+1}) \\ 
\downarrow &  & \downarrow &  & \downarrow \\ 
\Gamma (U\cap \mathcal{M}(N_{2}),\mathcal{F}^{m+1}) & \hookrightarrow & 
\Gamma (U\cap \mathcal{M}(N_{2}),\mathcal{F}^{m}) & \longrightarrow & \Gamma
(U\cap \mathcal{M}(N_{2}),\mathcal{F}^{m}/\mathcal{F}^{m+1})%
\end{array}%
\end{equation*}%
the left and the right vertical arrows are isomorphisms, and the middle
vertical arrow is a monomorphism, so it is also an isomorphism.$\square $

\bigskip

The rank four bundle $\mathcal{L}$ bears a connection 
\begin{equation}
\nabla ^{GM}=(\nabla ,\nabla ^{\prime }):\mathcal{\mathcal{L}}\overset{}{%
\longrightarrow }\mathcal{\mathcal{L}}\otimes _{\mathcal{O}%
_{Y(N_{2})}}\Omega _{Y(N_{2})}(\log D^{c})\overset{}{\approx }(\mathcal{%
\mathcal{L}}\otimes L_{Y(N_{2})}^{(2,0)})\oplus (\mathcal{\mathcal{L}}%
\otimes L_{Y(N_{2})}^{(0,2)})\text{ ,}
\end{equation}%
with $\log D^{c}$-poles \ (the Gauss-Manin connection) .The composition of
the second covariant derivative 
\begin{equation*}
\nabla ^{GM}:(\mathcal{L}^{n}\otimes L_{Y(N)}^{(2,0)})\oplus (\mathcal{L}%
^{n}\otimes L_{Y(N)}^{(0,2)})\longrightarrow \mathcal{L}^{n}\otimes
L_{Y(N)}^{(2,2)}
\end{equation*}%
with the inclusions of the two direct factors, are 
\begin{align*}
\nabla ^{\prime }& :\mathcal{L}^{n}\otimes L_{Y(N)}^{(2,0)}\longrightarrow 
\mathcal{L}^{n}\otimes L_{Y(N)}^{(2,2)} \\
-\nabla & :\mathcal{L}^{n}\otimes L_{Y(N)}^{(0,2)}\longrightarrow \mathcal{L}%
^{n}\otimes L_{Y(N)}^{(2,2)}
\end{align*}%
so we denote $\nabla ^{GM}=\nabla ^{\prime }-\nabla $. The Gauss-Manin
connection is integrable in the sense that 
\begin{equation}
\mathcal{L}^{n}\overset{\nabla ^{GM}=(\nabla ,\nabla ^{\prime })}{%
\longrightarrow }\mathcal{L}^{n}\otimes \Omega _{Y(N_{2})}(\log D^{c})%
\overset{\nabla ^{GM}=\nabla -\nabla ^{\prime }}{\longrightarrow }\mathcal{L}%
^{n}\otimes \omega _{Y(N_{2})}(D^{c})
\end{equation}%
is a complex of sheaves of $\mathbb{Q}_{p}$-spaces on $Y(N_{2})$ .

The hypercohomology of this complex is the $\mathcal{L}^{n}$ -valued $\log
D^{c}-$de Rham cohomology \ $H_{\log -dR}(Y(N_{2}),\mathcal{L}^{n}/\mathbb{Q}%
_{p})$of $Y(N_{2})$. There is again a $\mathbb{Q}_{p}$-linear map 
\begin{eqnarray*}
H_{\log -dR}^{i}(Y(N_{2}),\mathcal{L}^{n}/\mathbb{Q}_{p}) &\approx & \\
H_{\log -rig}^{i}(Y_{\kappa }(N_{2}),\mathcal{L}^{n}/\mathbb{Q}_{p})
&\longrightarrow &H_{\log -rig}^{i}(Y_{\kappa }^{\prime }(N_{2}),\mathcal{L}%
^{n}/\mathbb{Q}_{p})
\end{eqnarray*}%
to the analogously defined $\mathcal{L}^{n}-$valued $\log D^{c}-$rigid
cohomology of $Y_{\kappa }^{\prime }(N_{2})$ over $\mathbb{Q}_{p}$ , induced
by the restriction homomorphism%
\begin{equation*}
\begin{array}{ccccc}
\mathcal{L}_{\left] Y_{\kappa }\right[ }^{n} & \overset{}{\longrightarrow }
& \Omega _{\left] Y_{\kappa }\right[ }^{1}\otimes \mathcal{L}_{\left]
Y_{\kappa }\right[ }^{n} & \overset{}{\longrightarrow } & \omega _{\left]
Y_{\kappa }\right[ }\otimes \mathcal{L}_{\left] Y_{\kappa }\right[ }^{n} \\ 
\downarrow &  & \downarrow &  & \downarrow \\ 
j^{\dag }\mathcal{L}_{\left] Y_{\kappa }\right[ }^{n} & \overset{}{%
\longrightarrow } & j^{\dag }\Omega _{\left] Y_{\kappa }\right[ }^{1}\otimes 
\mathcal{L}_{\left] Y_{\kappa }\right[ }^{n} & \overset{}{\longrightarrow }
& j^{\dag }\omega _{\left] Y_{\kappa }\right[ }\otimes \mathcal{L}_{\left]
Y_{\kappa }\right[ }^{n}%
\end{array}%
\end{equation*}%
between the $\nabla ^{GM}-$complexes defining both $\log D^{c}-$rigid
cohomologies.

\bigskip Taking the first direct image $R^{1}pr_{\ast }$ in (\ref{Par
morfismo en abeliana}), we obtain%
\begin{equation}
\phi ^{\ast }\mathcal{L}_{\mathcal{M}_{\varepsilon ^{\prime
}}}\longrightarrow \mathcal{L}_{\mathcal{M}_{\varepsilon }}\text{ \ and \ }%
\phi _{\ast }\mathcal{L}_{\mathcal{M}_{\varepsilon }}\longrightarrow 
\mathcal{L}_{\mathcal{M}_{\varepsilon ^{\prime }}}
\end{equation}

\bigskip analogous to (\ref{Par morfismos en Rapoport}), and from this 
\begin{equation*}
\phi ^{\ast }(\Omega _{\mathbb{Q}_{p}}^{i}\otimes \mathcal{L}_{\mathcal{M}%
_{\varepsilon ^{\prime }}}^{n})\longrightarrow \Omega _{\mathbb{Q}%
_{p}}^{i}\otimes \mathcal{L}_{\mathcal{M}_{\varepsilon }}^{n}
\end{equation*}%
so that $Ver$ $\ $and $\ \phi $ act on $j^{\dag }\mathcal{L}_{\left] 
\mathcal{M}_{\kappa }\right[ }^{n}$ , thus also on $j^{\dag }\mathcal{L}_{%
\left] Y_{\kappa }\right[ }^{n}$, because of the Koecher principle for $%
\mathcal{L}^{n\text{ }},$ and, being these actions compatible with the
covariant derivative, it induces actions of $Ver$ and $\phi $ on $H_{\log
-rig}^{i}(Y_{\kappa }^{\prime }(N_{2}),\mathcal{L}^{n}/\mathbb{Q}_{p})$ .

\bigskip

These cohomologies are too big for our purposes, as in case $n=0$, they not
the de Rham and the rigid cohomology but their $\log D^{c}$\bigskip
-versions, so we switch to the parabolic-de Rham and parabolic-rigid
cohomology, which for $n=0$ gives the usual de Rham and the usual rigid
cohomologies.\bigskip\ That is to say, in the de Rham case, the cohomology $%
H_{par-dR}^{i}(Y(N_{2}),\mathcal{L}^{n}/\mathbb{Q}_{p})$ of the subcomplex
of sheaves of $\mathbb{Q}_{p}-$spaces 
\begin{equation*}
(\Omega _{Y(N_{2})}^{i}\otimes \mathcal{L}^{n})_{par}:=(\Omega
_{Y(N_{2})}^{i}\otimes \mathcal{L}^{n})+\nabla ^{GM}(\Omega
_{Y(N_{2})}^{i-1}\otimes \mathcal{L}^{n})\subseteq \Omega
_{Y(N_{2})}^{i}(\log D^{c})\otimes \mathcal{L}^{n}
\end{equation*}%
i.e. the cohomology of the de Rham "parabolic" complex%
\begin{equation*}
\mathcal{\mathcal{L}}_{Y(N_{2})}^{n}\longrightarrow (\Omega
_{Y(N_{2})}\otimes \mathcal{L}^{n})+\nabla ^{GM}\mathcal{L}%
_{Y(N_{2})}^{n}\longrightarrow (\omega _{Y(N_{2})}\otimes \mathcal{L}%
^{n})+\nabla ^{GM}(\Omega _{Y(N_{2})}\otimes \mathcal{L}^{n})\text{ \ \ \ \
\ \ \ \ \ \ \ \ \ \ \ \ \ \ \ \ \ \ \ }
\end{equation*}%
Note that $(\Omega _{Y(N_{2})}^{i}\otimes \mathcal{L}^{n})_{par}$ is a
locally free sheaf as kernel of the projection%
\begin{eqnarray*}
\Omega _{Y(N_{2})}^{i}(\log D^{c})\otimes \mathcal{L}^{n}
&\twoheadrightarrow &(\Omega _{Y(N_{2})}^{i}(\log D^{c})\otimes \mathcal{L}%
^{n})\otimes \mathcal{O}_{D^{c}} \\
&\twoheadrightarrow &(\Omega _{Y(N_{2})}^{i}(\log D^{c})\otimes \mathcal{L}%
^{n})\otimes \mathcal{O}_{D^{c}}/(\func{Im}\nabla ^{GM})_{D^{c}})\otimes 
\mathcal{O}_{D^{c}}
\end{eqnarray*}%
Note also that $\nabla ^{GM}(\mathcal{L}^{n}\otimes \Omega
_{Y(N_{2})})_{par}=\nabla ^{GM}(\mathcal{L}^{n}\otimes \Omega _{Y(N_{2})})$
because $\nabla ^{GM}$ is integrable.

The parabolic rigid cohomology $H_{par-rig}^{i}(Y_{\kappa }^{\prime }(N_{2}),%
\mathcal{L}^{n}/\mathbb{Q}_{p})$ is analogously defined as hypercohomology
of the $j^{\dag }(\Omega _{\left] Y_{\kappa }\right[ }^{i}\otimes \mathcal{L}%
_{\left] Y_{\kappa }\right[ }^{n})_{par}$ complex, i.e.

\begin{eqnarray}
j^{\dag }\mathcal{\mathcal{L}}_{\left] Y_{\kappa }\right[ }^{n}
&\longrightarrow &j^{\dag }(\Omega _{\left] Y_{\kappa }\right[ }\otimes 
\mathcal{L}_{\left] Y_{\kappa }\right[ }^{n})+j^{\dag }\nabla ^{GM}\mathcal{L%
}_{\left] Y_{\kappa }\right[ }^{n}  \label{complejo parabolico rigido} \\
&\longrightarrow &j^{\dag }(\omega _{\left] Y_{\kappa }\right[ }\otimes 
\mathcal{L}_{\left] Y_{\kappa }\right[ }^{n})+j^{\dag }\nabla ^{GM}(\Omega _{%
\left] Y_{\kappa }\right[ }\otimes \mathcal{L}_{\left] Y_{\kappa }\right[
}^{n}),\text{ \ \ \ \ \ }  \notag
\end{eqnarray}%
and it is acted by $\phi $ and $Ver$ because the parabolic complex inherits
such actions.

There is a $\mathbb{Q}_{p}$ -linear map, compatible with the $\phi $ action, 
\begin{equation*}
\begin{array}{c}
H_{par-dR}^{i}(Y(N_{2}),\mathcal{L}^{n}/\mathbb{Q}_{p})\approx
H_{par-rig}^{i}(Y_{\kappa }(N_{2}),\mathcal{L}^{n}/\mathbb{Q}_{p}) \\ 
\longrightarrow H_{par-rig}^{i}(Y_{\kappa }^{\prime }(N_{2}),\mathcal{L}^{n}/%
\mathbb{Q}_{p})%
\end{array}%
\end{equation*}%
induced by the restriction homomorphism%
\begin{equation}
\begin{array}{ccccc}
\mathcal{\mathcal{L}}_{\left] Y_{\kappa }\right[ }^{n} & \overset{\nabla
^{GM}}{\longrightarrow } & (\Omega _{\left] Y_{\kappa }\right[ }\otimes 
\mathcal{L}_{\left] Y_{\kappa }\right[ }^{n})+\nabla ^{GM}\mathcal{L}_{\left]
Y_{\kappa }\right[ }^{n} & \overset{\nabla ^{GM}}{\longrightarrow } & 
(\omega _{\left] Y_{\kappa }\right[ }\otimes \mathcal{L}_{\left] Y_{\kappa }%
\right[ }^{n})+\nabla ^{GM}(\Omega _{\left] Y_{\kappa }\right[ }\otimes 
\mathcal{L}_{\left] Y_{\kappa }\right[ }^{n}) \\ 
\downarrow &  & \downarrow &  & \downarrow \\ 
j^{\dag }\mathcal{\mathcal{L}}_{\left] Y_{\kappa }\right[ }^{n} & \overset{%
\nabla ^{GM}}{\longrightarrow } & j^{\dag }(\Omega _{\left] Y_{\kappa }%
\right[ }\otimes \mathcal{L}_{\left] Y_{\kappa }\right[ }^{n})+j^{\dag
}\nabla ^{GM}\mathcal{L}_{\left] Y_{\kappa }\right[ }^{n} & \overset{\nabla
^{GM}}{\longrightarrow } & j^{\dag }(\omega _{\left] Y_{\kappa }\right[
}\otimes \mathcal{L}_{\left] Y_{\kappa }\right[ }^{n})+j^{\dag }\nabla
^{GM}(\Omega _{\left] Y_{\kappa }\right[ }\otimes \mathcal{L}_{\left]
Y_{\kappa }\right[ }^{n})\text{ }%
\end{array}%
\end{equation}%
between the complexes defining these cohomologies.

Being just a formal homological tool, there is a Frobenius-invariant Mayor
-Vietoris sequence relating these parabolic cohomologies (de Rham and rigid)
for $Y(N_{2})$ and for $Y^{\prime }(N_{2})$ with their versions as parabolic
cohomologies for $Y(N_{2})$ supported in $D^{h}$(i.e. parabolic version of
the rigid supported cohomology defined in \cite{Besser inventi} 7.8 )

The natural inclusion $\mathcal{R\subseteq \mathcal{L}}$ has over the region 
$\mathcal{\mathcal{A}}$ in (\ref{wide opens}), a splitting 
\begin{equation*}
sp:\mathcal{\mathcal{L}}_{\mathcal{A}}\twoheadrightarrow \mathcal{R}_{%
\mathcal{A}}
\end{equation*}%
equivalent to a decomposition in Frobenius-invariant components%
\begin{equation}
\mathcal{\mathcal{L}}_{\mathcal{A}}=\mathcal{R}_{\mathcal{A}}\mathcal{\oplus 
\mathcal{S}}_{\mathcal{A}}=(L_{\mathcal{A}}^{(1,0)}\oplus L_{\mathcal{A}%
}^{(0,1)})\oplus (L_{\mathcal{A}}^{(-1,0)}\oplus L_{\mathcal{A}}^{(0,-1)})%
\text{ ,}  \label{ponerlo un numero}
\end{equation}%
defined in the noncuspidal points of $\mathcal{A}$ by the slope of the
Frobenius action, and extended by normality to the noncuspidal points. The
subbundle $\mathcal{S}_{\mathcal{A}}\subseteq \mathcal{L}_{\mathcal{A}}$ is
horizontal, i.e. $\nabla ^{GM}(\mathcal{S}_{\mathcal{A}})\subseteq \mathcal{S%
}_{\mathcal{A}}$ (cf. \cite{Split Katz} 1.11.26 and 27). This provides in
turn a split of the symmetric power of the rank four bundle $\mathcal{%
\mathcal{L}}_{\mathcal{A}}$ 
\begin{equation}
sp:\mathcal{L}_{\mathcal{A}}^{n}\twoheadrightarrow \mathcal{R}_{\mathcal{A}%
}^{n}=\dbigoplus\limits_{\substack{ k+k^{\prime }=n  \\ k,k^{\prime }\geq 0}}%
^{{}}L_{\mathcal{A}}^{(k,k^{\prime })}\text{ ,}  \label{split projection}
\end{equation}%
i.e. the exact sequence (\ref{nonsplit symmetric power}) splits in
Frobenius-stable direct factors%
\begin{equation*}
\mathcal{L}_{\mathcal{A}}^{n}=\mathcal{R}_{\mathcal{A}}^{n}\oplus \mathcal{S}%
_{\mathcal{A}}^{(n)}
\end{equation*}

\bigskip

\textbf{Lemma.\label{Complement is horizontal}} The complement $\mathcal{S}_{%
\mathcal{A}}^{(n)}$ is $\nabla ^{GM}-$horizontal, i.e. $\nabla ^{GM}(%
\mathcal{S}_{\mathcal{A}}^{(n)})\subseteq \mathcal{S}_{\mathcal{A}}^{(n)}$

\textbf{Proof}: The filtration $\mathcal{F}^{\bullet}$ splits on $\mathcal{A}
$, i.e.

\begin{equation*}
\mathcal{L}_{\mathcal{A}}^{n}=\dbigoplus \limits_{m=0}^{n}(\mathcal{R}_{%
\mathcal{A}}^{m}\otimes\mathcal{S}_{\mathcal{A}}^{n-m})\mathcal{=R}_{%
\mathcal{A}}^{n}\oplus\mathcal{S}_{\mathcal{A}}^{(n)}\text{ \ with }\mathcal{%
S}_{\mathcal{A}}^{(n)}=\text{ }\dbigoplus \limits_{m=0}^{n-1}(\mathcal{R}_{%
\mathcal{A}}^{m}\otimes\mathcal{S}_{\mathcal{A}}^{n-m})
\end{equation*}

For all $m\leq n-1$, 
\begin{equation*}
\nabla^{GM}(\mathcal{R}_{\mathcal{A}}^{m}\otimes\mathcal{S}_{\mathcal{A}%
}^{n-m})\subseteq(\nabla^{GM}(\mathcal{R}_{\mathcal{A}}^{m})\otimes \mathcal{%
S}_{\mathcal{A}}^{n-m})\oplus(\mathcal{R}_{\mathcal{A}}^{m}\otimes%
\nabla^{GM}(\mathcal{S}_{\mathcal{A}}^{n-m}))
\end{equation*}
The second direct factor is contained in $\mathcal{R}_{\mathcal{A}%
}^{m}\otimes\mathcal{S}_{\mathcal{A}}^{n-m}\subseteq$ $\mathcal{S}_{\mathcal{%
A}}^{(n)}$ and the first direct factor is contained in 
\begin{equation*}
\mathcal{\mathcal{L}}_{\mathcal{A}}^{m}\otimes\mathcal{S}_{\mathcal{A}%
}^{n-m}=(\dbigoplus \limits_{j=0}^{m}(\mathcal{R}_{\mathcal{A}}^{j}\otimes%
\mathcal{S}_{\mathcal{A}}^{m-j}))\otimes\mathcal{S}_{\mathcal{A}%
}^{n-m}=\dbigoplus \limits_{j=0}^{m}(\mathcal{R}_{\mathcal{A}}^{j}\otimes%
\mathcal{S}_{\mathcal{A}}^{n-j})\subseteq\mathcal{S}_{\mathcal{A}}^{(n)}
\end{equation*}

$\square$

\textbf{Remark. }Each tangent field $t$ on $\mathcal{M}=\mathcal{M}(N_{2})$
induces a derivation $\nabla _{t}^{GM}:\mathcal{L}_{\mathcal{M}}\overset{}{%
\longrightarrow }\mathcal{L}_{\mathcal{M}}$ , this providing, as explained
in \cite{Split Katz} 1.0.17 to 1.0.21, the isomorphism, as rank two $%
\mathfrak{\mathcal{O}}_{\mathcal{M}}$-sheaves, of the invertible $\mathfrak{%
o\otimes \mathcal{O}}_{\mathcal{M}}-$sheaves $\mathcal{T}_{\mathcal{M}}$ and 
$\mathcal{R}_{\mathcal{M}}^{\otimes (-2)}$ so that the first decomposes as

\begin{align*}
\mathcal{T}_{\mathcal{M}}& \approx L_{\mathcal{M}}^{(-2,0)}\oplus L_{%
\mathcal{M}}^{(0,-2)} \\
& =\mathcal{H}om_{\mathcal{O}_{\mathcal{M}}}(L_{\mathcal{M}}^{(1,0)},L_{%
\mathcal{M}}^{(-1,0)})\oplus \mathcal{H}om_{\mathcal{O}_{\mathcal{M}}}(L_{%
\mathcal{M}}^{(0,1)},L_{\mathcal{M}}^{(0,-1)})
\end{align*}%
\ For a tangent field $\tau $ in the first direct factor $L_{\mathcal{M}%
}^{(-2,0)}$ of $\mathcal{T}_{\mathcal{M}}$, and a local section $s$ of $%
\mathcal{L}$ , this says the following:\ if $s$ lies in $L_{\mathcal{M}%
}^{(1,0)}$ , then then the projection to $\mathcal{S}_{\mathcal{M}}$ of $%
\nabla _{\tau }^{GM}(s)$ lies in $L_{\mathcal{M}}^{(-1,0)}$ ;\ if $s$ lies
in $L_{\mathcal{M}}^{(0,1)}$ then $\nabla _{\tau }^{GM}(s)=0$ , and the
analogous holds for $\tau ^{\prime }$ in $L_{\mathcal{M}}^{(0,-2)}$.

Extending the remark to the whole of $Y(N_{2})$ (understood as subindex in
the notations), the dual $\mathcal{O}_{Y(N_{2})}-$isomorphism is $\Omega
_{Y(N_{2})}(\log D^{c})\approx \mathcal{R}^{\otimes 2}$ , and all three
compositions%
\begin{align*}
L^{(0,1)}& \hookrightarrow \mathcal{L}\overset{\nabla }{\longrightarrow }%
\mathcal{L}\otimes L^{(2,0)}\twoheadrightarrow L^{(-1,0)}\otimes L^{(2,0)} \\
L^{(0,1)}& \hookrightarrow \mathcal{L}\overset{\nabla }{\longrightarrow }%
\mathcal{L}\otimes L^{(2,0)}\twoheadrightarrow L^{(0,-1)}\otimes L^{(2,0)} \\
L^{(1,0)}& \hookrightarrow \mathcal{L}\overset{\nabla }{\longrightarrow }%
\mathcal{L}\otimes L^{(2,0)}\twoheadrightarrow L^{(0,-1)}\otimes L^{(2,0)}
\end{align*}%
are null, and 
\begin{equation*}
L^{(1,0)}\hookrightarrow \mathcal{L}\overset{\nabla }{\longrightarrow }%
\mathcal{L}\otimes L^{(2,0)}\twoheadrightarrow L^{(-1,0)}\otimes L^{(2,0)}
\end{equation*}%
is an isomorphism; and analogously $\nabla ^{\prime }$ provides the
isomorphism $L^{(0,1)}\approx L^{(0,-1)}\otimes L^{(0,2)}$ . In particular,
both $\nabla $ and $\nabla ^{\prime }$ preserve the decomposition $\mathcal{%
L\approx }$ $\mathcal{L}^{(1,0)}$ $\oplus \mathcal{L}^{(0,1)}$, a comment in 
\cite{Split Katz} to 2.3.18.

To understand the behavior of the covariant derivative on the symmetric
powers $\mathcal{L}^{n}$ of the rank four bundle $\mathcal{L}$ and its
filtrations, let us consider the canonical trivialization of the rank two
bundle $\mathcal{R}_{\mathcal{M}}\mathcal{\approx }L_{\mathcal{M}%
}^{(1,0)}\oplus L_{\mathcal{M}}^{(0,1)}$ in the punctured formal
neighborhood of the standard cusp by the canonical section $\omega _{can}$ \
(cf. comment of (\ref{expansion})), now just written $\omega $ for
simplicity, and its conjugate $\omega ^{\prime }$. This gives a
trivialization of the tangent bundle $\mathcal{\mathcal{T}}_{\mathcal{M}}%
\mathcal{\approx }L_{\mathcal{M}}^{(-2,0)}\oplus L_{\mathcal{M}}^{(0,-2)}$
by the square dual bases $\tau =\omega ^{-2}$ and $\tau ^{\prime }=\omega
^{\prime -2}$ . The rank four bundle $\mathcal{L}$ is trivialized, by adding
to $\omega $ and $\omega ^{\prime }$ the local sections 
\begin{equation*}
\eta :=\nabla _{\tau }^{GM}(\omega )\text{ and }\eta ^{\prime }:=\nabla
_{\tau ^{\prime }}^{GM}(\omega ^{\prime })
\end{equation*}%
and there are also the identities $\nabla _{\tau ^{\prime }}^{GM}(\omega )=0$
\ and $\nabla _{\tau }^{GM}(\omega ^{\prime })=0$. Therefore,

\begin{equation*}
\nabla ^{GM}\eta =0\text{ \ \ and }\nabla ^{GM}\eta ^{\prime }=0
\end{equation*}

The vector bundle $\mathcal{L}^{n}$ is consequently trivialized by sections%
\begin{equation*}
(\omega )^{k}(\omega ^{\prime })^{k^{\prime }}(\eta )^{l}(\eta ^{\prime
})^{l^{\prime }}\text{ with }k+k^{\prime }+l+l^{\prime }=n
\end{equation*}%
in terms of which

\begin{eqnarray}
&&\bigskip \nabla _{\tau }^{GM}f(\omega )^{k}(\omega ^{\prime })^{k^{\prime
}}(\eta )^{l}(\eta ^{\prime })^{l^{\prime }}  \label{derivada cov explicita}
\\
&=&\theta (f)(\omega )^{k}(\omega ^{\prime })^{k^{\prime }}(\eta )^{l}(\eta
^{\prime })^{l^{\prime }}+kf(\omega )^{k-1}(\omega ^{\prime })^{k^{\prime
}}(\eta )^{l+1}(\eta ^{\prime })^{l^{\prime }}  \notag
\end{eqnarray}%
and analogously for $\bigskip \nabla _{\tau ^{\prime }}^{GM}$ .

\bigskip We express the outcome of this remark in terms of the above
filtration $\mathcal{F}^{\bullet }$ and of a finer filtration $\mathcal{F}%
^{\bullet \bullet }$ defined by taking $\mathcal{F}^{(k,0)}\subseteq 
\mathcal{L}^{n}$ the image of

\begin{equation*}
L^{(k,0)}\otimes \mathcal{L}^{n-k}\hookrightarrow \mathcal{R}^{k}\otimes 
\mathcal{L}^{n-k}\hookrightarrow \mathcal{L}^{k}\otimes \mathcal{L}%
^{n-k}\twoheadrightarrow \mathcal{L}^{n}
\end{equation*}%
and analogous $\mathcal{F}^{(0,k^{\prime })}$ $\subseteq \mathcal{L}^{n}$,
and by being $\mathcal{F}^{(k,k^{\prime })}:=\mathcal{F}^{(k,0)}\cap 
\mathcal{F}^{(0,k^{\prime })}\subseteq \mathcal{L}^{n}$ . The outcome is
that these filtrations are compatible with the Gauss-Manin-connection in the
sense that

\begin{eqnarray*}
\nabla ^{GM}(\mathcal{F}^{m}) &\subseteq &\mathcal{F}^{m-1}\otimes \Omega
_{Y(N_{2})}(\log D^{c}) \\
\nabla (\mathcal{F}^{(k,k^{\prime })}) &\subseteq &\mathcal{F}%
^{(k-1,k^{\prime })}\otimes L^{(2,0)} \\
\nabla ^{\prime }(\mathcal{F}^{(k,k^{\prime })}) &\subseteq &\mathcal{F}%
^{(k,k^{\prime }-1)}\otimes L^{(0,2)}
\end{eqnarray*}%
so to induce isomorphisms

\begin{eqnarray}
\nabla &:&\text{$\mathcal{F}$}^{(k,k^{\prime })}/\text{$\mathcal{F}$}%
^{(k+k^{\prime }+1)}\overset{\approx }{\longrightarrow }(\text{$\mathcal{F}$}%
^{(k-1,k^{\prime })}/\text{$\mathcal{F}$}^{k+k^{\prime }})\otimes L^{(2,0)} 
\notag \\
\nabla ^{\prime } &:&\text{$\mathcal{F}$}^{(k,k^{\prime })}/\text{$\mathcal{F%
}$}^{(k+k^{\prime }+1)}\overset{\approx }{\longrightarrow }(\text{$\mathcal{F%
}$}^{(k,k^{\prime }-1)}/\text{$\mathcal{F}$}^{k+k^{\prime }})\otimes
L^{(0,2)}  \label{dos isomorfismos}
\end{eqnarray}

\begin{equation*}
(\mathcal{F}_{\mathbb{\mathcal{W}}_{\varepsilon }}^{m}/\mathcal{F}_{\mathbb{%
\mathcal{W}}_{\varepsilon }}^{m+1})\otimes \omega _{\mathbb{\mathcal{W}}%
_{\varepsilon }}(\log D^{c})\approx \dbigoplus\limits_{i+j=m}^{{}}\mathcal{F}%
_{\mathbb{\mathcal{W}}_{\varepsilon }}^{(i,j)}/\mathcal{F}_{\mathbb{\mathcal{%
W}}_{\varepsilon }}^{(i,j)}\cap \mathcal{F}_{\mathbb{\mathcal{W}}%
_{\varepsilon }}^{m+1}\otimes L_{\mathbb{\mathcal{W}}_{\varepsilon }}^{(2,2)}
\end{equation*}

\textbf{Definition\label{Defino near- oc y equivalencia}}. We say that a
rigid $L_{\mathcal{A}}^{(k,k^{\prime })}$-differential is\textit{\ nearly
overconvergent\ }$\varrho $ when it is\ in the image of \ 
\begin{equation*}
\Gamma _{rig}(\Omega _{\mathcal{W}_{\varepsilon }}^{i}\otimes \mathcal{L}_{%
\mathcal{W}_{\varepsilon }}^{n})_{par}\longrightarrow \Gamma _{rig}(\Omega _{%
\mathcal{A}}^{i}\otimes \mathcal{L}_{\mathcal{A}}^{n})_{par}\overset{sp^{n}}{%
\longrightarrow }\Gamma _{rig}(\Omega _{\mathcal{\mathcal{A}}}^{i}\otimes 
\mathcal{R}_{\mathcal{A}}^{n})
\end{equation*}%
where the last morphism is the restriction of the split morphism 
\begin{equation*}
sp^{n}:\Gamma _{rig}(\Omega _{\mathcal{A}}^{i}(\log D^{c})\otimes \mathcal{L}%
_{\mathcal{A}}^{n})\longrightarrow \Gamma _{rig}(\Omega _{\mathcal{A}%
}^{i}(\log D^{c})\otimes \mathcal{R}_{\mathcal{A}}^{n})\text{ ,}
\end{equation*}%
(The image of $\Gamma _{rig}(\Omega _{\mathcal{A}}^{i}\otimes \mathcal{L}_{%
\mathcal{A}}^{n})_{par}$ lies in $\Gamma _{rig}(\Omega _{\mathcal{\mathcal{A}%
}}^{i}\otimes \mathcal{R}_{\mathcal{A}}^{n})$ because of the description of
the $\nabla ^{GM}$ -derivative of local sections, in the former remark, in
terms of our choice of a local base of $\mathcal{L}$).

We say that two nearly overconvergent 1-differentials are \textit{equivalent}
$\varrho \approx \sigma $ when $e_{ord}j_{X^{\prime }(N)}^{\ast }(\varrho
-\sigma )=0$ (the reason is that they will have same contribution to the
cup-product computing the $p-$adic Abel Jacobi map); and say two
overconvergent $L^{(k,k^{\prime })}-$2-differentials are equivalent when its
difference has a $\nabla ^{GM}$ -primitive equivalent to $0$.

\bigskip

\ Considering now the generically injective morphism $j_{X(N)}:X(N)%
\longrightarrow Y(N_{2})$ , the splitting $sp$ induces, via $j_{X(N)}^{\ast
} $ , the equally denoted splitting $sp:\mathcal{\mathcal{L}}_{0,\mathcal{A}%
_{0}}\twoheadrightarrow L_{0,\mathcal{A}_{0}}$ of the restriction to $%
\mathcal{A}_{0}:=j_{X(N)}^{-1}\mathcal{A}$ of \ref{main short sequence in
dim 1} of each of the two direct factors of $j_{X(N)}^{\ast }\mathcal{%
L\approx \mathcal{L}}_{0}\oplus \mathcal{L}_{0}$ , hence an splitting 
\begin{equation*}
sp:\mathcal{\mathcal{L}}_{0,\mathcal{A}_{0}}^{n}\twoheadrightarrow L_{0,%
\mathcal{A}_{0}}^{n}
\end{equation*}%
of \ref{powered short sequence in dim 1}. Therefore, a nearly overconvergent 
$L^{(k,k^{\prime })}$ -differential in our sense "restricts", i.e. has image
by $j_{X(N)}^{\ast }$, to a nearly overconvergent $L_{0}^{k+k^{\prime }}$
-differential, in the sense of \cite{DR} def. 4.3, of which our definition
has been an analogue. We recall at this point that the Hida projector 
\begin{equation*}
e_{ord}=\lim_{n}U_{p}^{n!}\text{ , }
\end{equation*}%
\ from the space of overconvergent modular forms to the subspace of the
ordinary ones, extends to the space of nearly overconvergent forms, and has
the same image space, cf. lemma 2.7 in \cite{DR}.

\bigskip

\begin{proposition}
\textbf{\ \label{Clave en Lausanne}:}\ 
\end{proposition}

\textit{a) For a wide open set} $\mathbb{\mathcal{W}}_{\varepsilon }$ 
\textbf{, }\textit{any }$\varrho \mathit{\in \Gamma }_{rig}$\textit{\ }$(%
\mathcal{L}_{\mathbb{\mathcal{W}}_{\varepsilon }}^{n}\otimes \Omega _{%
\mathbb{\mathcal{W}}_{\varepsilon }})_{par}\subseteq \mathit{\Gamma }_{rig}%
\mathcal{L}_{\mathbb{\mathcal{W}}_{\varepsilon }}^{n}\otimes \Omega _{%
\mathbb{\mathcal{W}}_{\varepsilon }}(\log D^{c})$\textit{\ is sum }$\varrho
=\varrho _{n}+\nabla ^{GM}\alpha $\textit{\ for some }$\varrho _{n}\mathit{%
\in \Gamma }_{rig}\mathcal{R}_{\mathbb{\mathcal{W}}_{\varepsilon
}}^{n}\otimes \Omega _{\mathbb{\mathcal{W}}_{\varepsilon }}$ and $\alpha \in 
\mathit{\Gamma }_{rig}\mathcal{L}_{\mathbb{\mathcal{W}}_{\varepsilon }}^{n}$

b) $sp(\varrho )$ $\approx $ $\varrho _{n}$

\textbf{Proof:}\ 

Proof: a) This is particular case of the more general fact that any $\varrho
_{m}\in \Gamma _{rig}\mathcal{\mathcal{F}}_{\mathbb{\mathcal{W}}%
_{\varepsilon }}^{m}\otimes \Omega _{\mathbb{\mathcal{W}}_{\varepsilon
}}(\log D^{c})$ is sum $\varrho _{m}=\varrho _{m+1}+\nabla \alpha $ for some 
$\varrho _{m+1}\in \Gamma _{rig}\mathcal{\mathcal{F}}_{\mathbb{\mathcal{W}}%
_{\varepsilon }}^{m+1}\otimes \Omega _{\mathbb{\mathcal{W}}_{\varepsilon }}$
and $\alpha \in \mathit{\Gamma }_{rig}\mathcal{L}_{\mathbb{\mathcal{W}}%
_{\varepsilon }}^{n}$, which we prove by induction. For each summand $%
\overline{\varrho }_{i,j}$ of the class $\overline{\varrho }%
_{m}=\dsum\limits_{i+j=m}^{{}}\overline{\varrho }_{i,j}$ of $\varrho _{m}$
in the decomposition 
\begin{equation*}
(\mathcal{F}_{\mathbb{\mathcal{W}}_{\varepsilon }}^{m}/\mathcal{F}_{\mathbb{%
\mathcal{W}}_{\varepsilon }}^{m+1})\otimes \omega _{\mathbb{\mathcal{W}}%
_{\varepsilon }}(\log D^{c})\approx \dbigoplus\limits_{i+j=m}^{{}}\mathcal{F}%
_{\mathbb{\mathcal{W}}_{\varepsilon }}^{(i,j)}/\mathcal{F}_{\mathbb{\mathcal{%
W}}_{\varepsilon }}^{(i,j)}\cap \mathcal{F}_{\mathbb{\mathcal{W}}%
_{\varepsilon }}^{m+1}\otimes L_{\mathbb{\mathcal{W}}_{\varepsilon }}^{(2,2)}
\end{equation*}%
one of the two indexes $i,j$ is positive, for instance $i>0$ , i.e. $j<m$,
so it is the image $\nabla ^{\prime }(\overline{\varrho }_{i,j+1})=\nabla
^{GM}(\overline{\varrho }_{i,j+1})$ for some 
\begin{equation*}
\overline{\varrho }_{i,j+1}\in \Gamma _{rig}(\mathcal{F}_{\mathbb{\mathcal{W}%
}_{\varepsilon }}^{(i,j+1)}/\mathcal{F}_{\mathbb{\mathcal{W}}_{\varepsilon
}}^{(i,j+1)}\cap \mathcal{F}_{\mathbb{\mathcal{W}}_{\varepsilon
}}^{m+2})\otimes L_{\mathbb{\mathcal{W}}_{\varepsilon }}^{(2,0)}
\end{equation*}%
because of the above isomorphism \ref{dos isomorfismos} (tensored by $L_{%
\mathbb{\mathcal{W}}_{\varepsilon }}^{(2,0)}$) 
\begin{equation*}
\partial ^{\prime }:(\mathcal{F}_{\mathbb{\mathcal{W}}_{\varepsilon
}}^{(i,j+1)}/\mathcal{F}_{\mathbb{\mathcal{W}}_{\varepsilon }}^{(i,j+1)}\cap 
\mathcal{F}^{m+2})\otimes L_{\mathbb{\mathcal{W}}_{\varepsilon }}^{(2,0)}%
\overset{\approx }{\longrightarrow }(\mathcal{F}_{\mathbb{\mathcal{W}}%
_{\varepsilon }}^{(i.j)}/\mathcal{F}_{\mathbb{\mathcal{W}}_{\varepsilon
}}^{(i.j)}\cap \mathcal{F}_{\mathbb{\mathcal{W}}_{\varepsilon
}}^{m+1})\otimes L_{\mathbb{\mathcal{W}}_{\varepsilon }}^{(2,2)}
\end{equation*}

\bigskip b) From part a) we obtain $j_{X^{\prime }(N)}^{\ast }\varrho
=j_{X^{\prime }(N)}^{\ast }\varrho _{n}+\nabla (j_{X^{\prime }(N)}^{\ast }s)$
, so that, by lemma 2.7 of \cite{DR},%
\begin{equation*}
e_{ord}(j_{X^{\prime }(N)}^{\ast }\varrho _{n})=e_{ord}(sp(j_{X^{\prime
}(N)}^{\ast }\varrho ))=e_{ord}(j_{X^{\prime }(N)}^{\ast }sp(\varrho ))
\end{equation*}

$\square $

Composing the covariant derivative of a rigid analytic section $\alpha $ of $%
\mathcal{R}_{\mathcal{A}}^{n}\subseteq \mathcal{L}_{\mathcal{A}}^{n}$\ with
the split projection $sp$ in (\ref{split projection}), we obtain what may be
called its \textit{"split derivative"} 
\begin{equation*}
\nabla _{sp}^{GM}\alpha :=sp(\nabla ^{GM}\alpha )=(\nabla _{sp}\alpha
,\nabla _{sp}^{\prime }\alpha )\text{ ,}
\end{equation*}%
which is a rigid analytic section on $\mathcal{A}$, defined over $\mathbb{Q}%
_{p},$ of\textbf{\ }%
\begin{equation*}
\mathcal{R}^{n}\otimes \Omega _{\mathcal{A}}(\log D^{c})\approx (\mathcal{R}%
^{n}\otimes L_{\mathcal{A}}^{(2,0)})\oplus (\mathcal{R}^{n}\otimes L_{%
\mathcal{A}}^{(0,2)})
\end{equation*}%
\textbf{\ }(cf. \cite{Split Katz} 2.5.7, $p-$adic analogue of 2.3.2 and
2.3.12). A \textit{key fact} which follows of the behavior (\ref{derivada
cov explicita}) of the covariant derivative (cf. also \cite{Split Katz}
2.5.12) is that, if $\alpha $ is in fact a section of a modular line
subbundle $L_{\mathcal{A}}^{(k,k^{\prime })}\subseteq \mathcal{R}^{n}$, i.e.
corresponds to a $p-$adic form of weight $(k,k^{\prime })$, then $\nabla
_{sp}^{GM}\alpha =(\nabla _{sp}\alpha ,\nabla _{sp}^{\prime }\alpha )$ lies
in 
\begin{equation*}
L_{\mathcal{A}}^{(k,k^{\prime })}\otimes \Omega _{\mathcal{A}}(\log
D^{c})\approx L_{\mathcal{A}}^{(k+2,k^{\prime })}\oplus L_{\mathcal{A}%
}^{(k,k^{\prime }+2)}\text{,}
\end{equation*}%
i.e. it corresponds to a pair of $p-$adic $\mathfrak{a}$-Hilbert modular
forms of weights $(k+2,k^{\prime })$ and $(k,k^{\prime }+2)$. We can thus
write this split derivative, by abusing notation, as%
\begin{equation}
\nabla _{sp}:L_{\mathcal{A}}^{(k,k^{\prime })}\longrightarrow L_{\mathcal{A}%
}^{(k+2,k^{\prime })}\text{ \ and \ \ }\nabla _{sp}^{\prime }:L_{\mathcal{A}%
}^{(k,k^{\prime })}\longrightarrow L_{\mathcal{A}}^{(k,k^{\prime }+2)}\text{.%
}
\end{equation}

Generalizing (\ref{expansion first der. of diff}), if $\alpha $ has $q$%
-expansion $\dsum\limits_{\nu \in (\mathfrak{a}^{-1})^{+}\cup
\{0\}}^{{}}a_{\nu }q^{\nu }$, then the split derivatives $\nabla _{sp}\alpha 
$ and $\nabla _{sp}^{\prime }\alpha $ have $q$-expansions in the given
trivializations

\begin{align}
\theta (\dsum\limits_{\nu \in (\mathfrak{a}^{-1})^{+}\cup \{0\}}^{{}}a_{\nu
}q^{\nu })& =\dsum\limits_{\nu \in (\mathfrak{a}^{-1})^{+}\cup
\{0\}}^{{}}\nu a_{\nu }q^{\nu }\text{ \ }  \label{zeta on expansions} \\
\text{and \ }\theta ^{\prime }(\dsum\limits_{\nu \in (\mathfrak{a}%
^{-1})^{+}\cup \{0\}}^{{}}a_{\nu }q^{\nu }))& =\text{ }\dsum\limits_{\nu \in
(\mathfrak{a}^{-1})^{+}\cup \{0\}}^{{}}\nu ^{\prime }a_{\nu }q^{\nu }\text{ .%
}
\end{align}%
This follows from the remark of Andreatta-Goren in \cite{And Go} 12.26 that
the Gauss-Manin covariant derivative studied by Katz in 2.6 \cite{Split Katz}
acts on $p$-adic modular forms as the theta operators studied by
Andreatta-Goren in sec. 12 of \cite{And Go}, i.e. those acting on $q$
-expansions as (\ref{zeta on expansions}) (the two basic cases among those
considered in 12.23.2 of \cite{And Go}).

Obviously, if $\alpha $ is an overconvergent $\mathcal{L}^{n}-$differential,
then $\nabla ^{GM}\alpha $ is overconvergent, i.e. both $\nabla \alpha $ and 
$\nabla ^{\prime }\alpha $ are overconvergent, so that\ $\nabla
_{sp}^{GM}\alpha $ is nearly overconvergent. We can say more:

\begin{lemma}
\label{Conservacion near converg} If $\alpha$ is a nearly overconvergent $%
L^{(k,k^{\prime})}-$differential, then $\nabla_{sp}^{GM}\alpha$ is nearly
overconvergent, i.e. both $\nabla_{sp}\alpha$\textit{\ and }$\nabla
_{sp}^{\prime}\alpha$\textit{\ are nearly \ overconvergent.}
\end{lemma}

\textbf{Proof: } Write $n=k+k^{\prime}$ . It is enough to prove it\ for a
rigid section $\alpha$ of $\mathcal{R}_{\mathcal{A}}^{n}\otimes\Omega _{%
\mathcal{A}}^{i}$ which is restriction of some rigid section $\sigma$ of $%
\mathcal{L}_{\mathbb{\mathcal{W}}_{\varepsilon}}^{n}\otimes\Omega _{\mathbb{%
\mathcal{W}}_{\varepsilon}}^{i}$, i.e. $\sigma_{\mathcal{A}}=\alpha+\beta$
with $\alpha\in\mathcal{R}_{\mathcal{A}}^{n}$ and $\beta \in\mathcal{S}_{%
\mathcal{A}}^{(n)}$ . The rigid section $\nabla_{sp}^{GM}\alpha$ of $%
\mathcal{R}_{\mathcal{A}}^{n}\otimes\Omega_{\mathcal{A}}^{i+1}$ is the split
projection of the rigid section $\nabla^{GM}\sigma$ of $(\mathcal{L}_{%
\mathcal{\mathbb{\mathcal{W}}}_{\varepsilon}}^{n}\otimes \Omega_{\mathcal{W}%
_{\varepsilon}}^{i+1})_{par}$ because $\nabla ^{GM}(\mathcal{S}_{\mathcal{%
\mathbb{\mathcal{W}}}_{\varepsilon}}^{(n)}\otimes\Omega_{\mathcal{W}%
_{\varepsilon}}^{i})\subseteq(\mathcal{S}_{\mathcal{\mathbb{\mathcal{W}}}%
_{\varepsilon}}^{(n)}\otimes\Omega _{\mathcal{W}_{\varepsilon}}^{i+1})_{par}$
.$\square$

\bigskip

Observe that the "key fact" works also for the second log $D^{c}$-covariant
derivative $\nabla ^{GM}$ as it is induced by the first: any overconvergent
section $(\beta ,\gamma )$ of $L^{(k,k^{\prime })}\otimes \Omega (\log
D^{c})=L^{(k+2,k^{\prime })}\oplus L^{(k,k^{\prime }+2)}$ defined over $%
\mathbb{Q}_{p}$, corresponding to a pair of $p-$adic $\mathfrak{a}$-Hilbert
modular forms over $\mathbb{Q}_{p}$ of weights $(k+2,k^{\prime })$ and $%
(k,k^{\prime }+2)$ with $q-$expansions$\dsum\limits_{\nu \in (\mathfrak{a}%
^{-1})^{+}\cup \{0\}}^{{}}b_{\nu }q^{n}$ and$\dsum\limits_{\nu \in (%
\mathfrak{a}^{-1})^{+}\cup \{0\}}^{{}}c_{\nu }q^{n}$ , has as split
derivative $\nabla _{sp}^{GM}(\beta ,\gamma )=\nabla _{sp}^{\prime }\beta
-\nabla _{sp}\gamma $ $\ $a rigid analytic section of $\mathcal{L}_{\mathcal{%
A}}^{k+k^{\prime }}\otimes \omega _{\mathcal{A}}(D^{c})$, defined over $%
\mathbb{Q}_{p}$, lying in fact inside $L_{\mathcal{A}}^{(k,k^{\prime
})}\otimes \omega _{\mathcal{A}}(D^{c})\approx L_{\mathcal{A}%
}^{(k+2,k^{\prime }+2)}$, thus corresponding to a nearly overconvergent $%
\mathfrak{a}$-Hilbert modular form over $\mathbb{Q}_{p}$ of weight $%
(k+2,k^{\prime }+2)$ and $q$-expansion%
\begin{equation}
\theta ^{\prime }(\dsum\limits_{\nu \in (\mathfrak{a}^{-1})^{+}\cup
\{0\}}^{{}}b_{\nu }q^{\nu })-\theta (\dsum\limits_{\nu \in (\mathfrak{a}%
^{-1})^{+}\cup \{0\}}^{{}}c_{\nu }q^{\nu })=\dsum\limits_{\nu \in (\mathfrak{%
a}^{-1})^{+}\cup \{0\}}^{{}}(\nu ^{\prime }b_{\nu }-\nu c_{\nu })q^{n}\text{
.}  \label{second cov dev on expansions}
\end{equation}%
\textbf{\ }This generalizes (\ref{expansion second derivative of dif2}).

\bigskip

\bigskip\ Finally, it will be usefull also to recall that the second
exterior power%
\begin{equation*}
\dbigwedge\limits^{2}\mathcal{L}_{Y(N_{2})}=\mathcal{H}_{dR}^{2}(\widehat{%
A^{U}}):=\mathbb{R}^{2}pr_{Y(N_{2})}{}_{\ast }\Omega _{\widehat{A^{U}}%
/Y(N_{2})}^{\bullet }(\log D_{\widehat{A^{U}}}^{c})
\end{equation*}%
inherits an equally denoted Gauss-Manin connection, and that the smallest
filter $Fil^{2}$ of its Hodge filtration is the "Norm" modular line bundle \ 
$L_{Y(N_{2})}^{(1,1)}$ in the sequence%
\begin{equation*}
0\longrightarrow L_{Y(N_{2})}^{(1,1)}=\text{$Fil$}^{2}\dbigwedge\limits^{2}%
\mathcal{L}_{Y(N_{2})}\longrightarrow \dbigwedge\limits^{2}\mathcal{L}%
_{Y(N_{2})}\longrightarrow \mathcal{S}_{2}\longrightarrow 0
\end{equation*}%
analogous to (\ref{extension buena de relative DR}), from which we have that 
$L_{Y(N_{2})}^{(n,n)}$ is the smallest filter $Fil^{2n}$ of the Hodge
filtration of $\dbigwedge\limits^{2}\mathcal{L}_{Y(N_{2})}^{n}$ .

\subsection{\textbf{Null-homologous cycles of higher dimension. }}

On the $n-$th product of a generalized elliptic curve $E$ consider $\epsilon
_{E^{n}}=\epsilon _{E^{n}}^{sym}\circ \epsilon _{E^{n}}^{inv}$ with 
\begin{eqnarray}
\epsilon _{E^{n}}^{sym} &=&\frac{1}{n!}\dsum_{\sigma \in \Sigma
_{n}}sign(\sigma )\sigma \in \mathbb{Q}[AutE^{n}]  \label{Scholl operator} \\
\epsilon _{E^{n}}^{inv} &=&\text{ }(\frac{1-u_{1}}{2})...(\frac{1-u_{n}}{2}%
)\in \mathbb{Q}[AutE^{n}]
\end{eqnarray}%
where all $u_{i}$ are equal to the involution $u$ of $E$ , and denote
equally the idempotent operator it induces on rational cycles of $E^{n}$;\
this is alternatively presented in \cite{BDP} (1.4.4) in terms of the map $%
j: $ $(\mu _{2})^{n}\ltimes \Sigma _{n}\longrightarrow \pm 1$ product of the
identity and the sign character: 
\begin{equation}
\epsilon _{E^{n}}=\frac{1}{2^{n}n!}\dsum\limits_{\tau \in (\mu
_{2})^{n}\ltimes \Sigma _{n}}^{{}}j(\tau )\tau \in \mathbb{Q}[AutE^{n}]
\label{alternative Scholl operator}
\end{equation}
It is observed in\textbf{\ }lemma 1.8 of \cite{BDP} that $\epsilon _{E}$
vanishes at all the K\"{u}neth pieces of the de Rham homology $%
H_{dR}^{\bullet }(E^{n})$ of $E^{n}$ , except the one $%
\dbigotimes^{n}H^{1}(E)$ of the intermediate dimension $n$ , where it acts
as the projection to its symetric part%
\begin{equation*}
\epsilon _{E^{n}}H_{dR}^{n}(E^{n})=Sym^{n}H_{dR}^{1}(E)
\end{equation*}%
or, by duality, in terms of the singular homology $H_{\bullet }(E^{n})$, it
is null at all pieces $H_{i}(E^{n})$ of order $i\neq n$ , and

\begin{equation*}
\epsilon_{E^{n}}H_{n}(E^{n})=Sym^{n}H_{1}(E)
\end{equation*}

We will denote by $\mathcal{E}$ the smooth compactification $\widehat{E}%
_{X(N)}^{U}$ (universal generalized elliptic curve) of the universal abelian
scheme $E_{X(N)}^{U}$ over $X(N)$, and denote $\mathcal{E}^{n}$ its $n-$th
product fibred over $X(N)$ (but will keep the notation $\widehat{A}%
_{Y(N_{2})}^{U,n}\overset{pr}{\longrightarrow }Y(N_{2})$ for the proper,
smooth variety of dimension $2n+2$ smooth $n-$th fibred product of the
compactified universal family $\widehat{A}_{Y(N_{2})}^{U}$over $Y(N_{2})$).
The operator $\epsilon _{\mathcal{E}^{n}}$ is defined fibrewise, and as a
consequence of the above killing of nonintermediate cohomology, it is 
\begin{equation}
\epsilon _{\mathcal{E}^{n}}H_{dR}^{n+1}(\mathcal{E}^{n}\mathcal{)=}%
H_{dR}^{1}(\mathcal{L}_{0}^{n})\text{ and }\epsilon _{\mathcal{E}%
^{n}}H_{dR}^{i}(\mathcal{E}^{n}\mathcal{)=}0\text{ for }i\neq n+1
\label{symmetric are fixed}
\end{equation}%
(lemma 2.2 of \cite{BDP}), i.e. $\epsilon _{\mathcal{E}^{n}}H_{i}(\mathcal{E}%
^{n}\mathcal{)}=0$ for all $i\neq n+1$.

\bigskip For an integer $n\geq 0\mathfrak{\ }$and even integer $0\mathfrak{%
\leq }n_{0}\mathfrak{<}2n$, denote $m=n+\frac{n_{0}}{2}$ . We now construct
null-homologous $(m+1)$-algebraic cycles $\Delta _{n,n_{0}}$ of $\widehat{A}%
_{Y(N_{2})}^{U,n}\times \mathcal{E}^{n_{0}}$, \ on which to compute the $p-$%
adic Abel-Jacobi map, using the generically injective embedding 
\begin{equation}
\mathcal{E}^{2n}=\mathcal{E}^{n}\times _{X(N)}\mathcal{\mathcal{E}}^{n}=(%
\mathcal{E}\times _{X(N)}\mathcal{\mathcal{E}})^{n}\longrightarrow \widehat{A%
}_{Y(N_{2})}^{U,n}  \label{main embedding}
\end{equation}%
given by \ref{decomp of universal compact restr}. We will use the idempotent
operator $\epsilon =\epsilon _{\mathcal{E}^{n}}\times _{X(N)}\epsilon _{%
\mathcal{E}^{n}}$ in cohomology of $\mathcal{E}^{n}\times _{X(N)}\mathcal{%
\mathcal{E}}^{n}$ , so restricting $\epsilon _{\mathcal{E}^{n}}$ to the
diagonal $\mathcal{E}^{n}\overset{diag}{\mathcal{\longrightarrow }}\mathcal{%
\mathcal{E}}^{n}\mathcal{\mathcal{\times }}_{X(N)}\mathcal{\mathcal{\mathcal{%
E}}}^{n}$ . It is restricted from the equally denoted operator on the whole
of $\widehat{A}_{Y(N_{2})}^{U,n}$ with a definition on its fibres analogous
to \ref{Scholl operator} but which does not have the property of killing all
nonintermediate cohomology, as $\epsilon =\epsilon _{\mathcal{E}^{n}}\times
_{X(N)}\epsilon _{\mathcal{E}^{n}}$ does:\ this last can be proved just as
for lemma 2.2 of \cite{BDP} i.e. $\epsilon H^{i}(\mathcal{E}^{n}\times
_{X(N)}\mathcal{\mathcal{E}}^{n}\mathcal{)=}0$ for $i\neq 2n+1$ and $%
\epsilon H^{2n+1}(\mathcal{E}^{n}\times _{X(N)}\mathcal{\mathcal{E}}^{n}%
\mathcal{)=}H^{1}(\mathcal{L}_{0}^{n}\otimes \mathcal{L}_{0}^{n})$ so $%
\epsilon $ fixes the cohomology classes lying inside this subspace.

We consider first the case $n_{0}=0$, so that $m=n.$ In case $n=0$, the
cycle $\Delta _{0,0}$ of $Y(N_{2})\times X(N)$ is the one already studied,
so we can assume $n\geq 1$ . Consider $\mathcal{G}r_{n}\mathcal{\approx E}%
^{n}$ the graph of the projection $pr_{n}:\mathcal{E}^{n}\longrightarrow
X(N) $ of the diagonal cycle%
\begin{equation*}
\mathcal{E}^{n}\overset{diag}{\mathcal{\longrightarrow }}\mathcal{\mathcal{E}%
}^{n}\mathcal{\mathcal{\times }}_{X(N)}\mathcal{\mathcal{\mathcal{E}}}%
^{n}\longrightarrow \widehat{A}_{Y(N_{2})}^{U,n}
\end{equation*}%
and define the cycle of $\widehat{A}_{Y(N_{2})}^{U,n}\times X(N)$ 
\begin{equation*}
\Delta _{n,0}:=(\epsilon \times id_{X(N)})(\mathcal{G}r_{n})=(\epsilon _{%
\mathcal{\mathcal{E}}^{n}}\times id_{X(N)})(\mathcal{G}r_{n})
\end{equation*}%
Since $\epsilon _{\mathcal{\mathcal{E}}^{n}}$ kills the homology of $%
\mathcal{E}^{n}$ of dimension different of $n+1$, the operator induced by $%
\epsilon _{\mathcal{\mathcal{E}}^{n}}\times id_{X(N)}$ in the homology of $%
\mathcal{E}^{n}\times X(N)$ kills its homology different of $n+3,n+2,n+1$.
Since $\mathcal{G}r_{n}\subseteq \mathcal{E}^{n}\times X(N)$ has dimension $%
2(n+1)$, its homology class is killed by the operator, if $n\geq 2.$

In case $n=1,n_{0}=0$ , we choose a rational point $o\in X(N)$ and write $%
E_{o}$ the generalized elliptic curve $pr^{-1}(o)$ , so to have a
decomposition%
\begin{equation*}
\mathcal{G}r_{1}=\mathcal{E}\times \{o\}\text{ }+n_{ij}\text{ }\dsum
E_{i}\times F_{j}+(E_{o}\times X(N))
\end{equation*}%
analogous to\ (\textbf{\ref{key topological formula})}\ for cycles $F_{j}$
making a base of the 1-homology of $X(N)$, and cycles $E_{i}$ making a base
of the 3-homology of $\mathcal{E}$ . Since this last homology is killed by
the action of $\epsilon _{\mathcal{E}}$, we obtain that the homology of 
\begin{eqnarray*}
\Delta _{1,0} &:&=(\epsilon \times id_{X(N)})(\mathcal{G}r-\mathcal{E}\times
\{o\}-E_{o}\times X(N)) \\
&=&(\epsilon _{\mathcal{E}}\times id_{X(N)})(\mathcal{G}r-\mathcal{E}\times
\{o\}-E_{o}\times X(N))
\end{eqnarray*}%
is null.

\bigskip Finally, consider the case $0<n_{0}<2n$ . Write $\{1,...,m\}$ as
union $A\cup B$ of two subsets of cardinality $n$ , so with an intersection $%
T=A\cap B$ of cardinality $t=n-\frac{n_{0}}{2}$ $=m-n_{0}>0$ . The
projections $pr_{A}:\mathcal{E}^{m}\longrightarrow \mathcal{E}^{n}$ and $%
pr_{B}:\mathcal{E}^{m}\longrightarrow \mathcal{E}^{n}$ to the arguments in $%
A $ and $B$ , induce a generically injective embedding 
\begin{equation}
\varphi _{m;2n}:\mathcal{E}^{m}\hookrightarrow \mathcal{E}^{n}\times _{X(N)}%
\mathcal{E}^{n}\approx \mathcal{E}^{2n}\longrightarrow \widehat{A}%
_{Y(N_{2})}^{U,n}  \label{inmersion en A^U}
\end{equation}%
the first being essentially the diagonal embedding 
\begin{equation*}
\mathcal{E}^{m}\approx \mathcal{E}^{t}\times _{X(N)}\mathcal{E}^{n_{0}}%
\overset{(\text{diag)}^{t}}{\hookrightarrow }(\mathcal{E}\times _{X(N)}%
\mathcal{E})^{t}\times _{X(N)}\mathcal{E}^{n_{0}}
\end{equation*}%
in the $t$ arguments of $T$ . Observe that $\epsilon $ restricts $\epsilon _{%
\mathcal{E}^{t}}\times _{X(N)}\epsilon _{\mathcal{E}^{n_{0}}}$ to $\mathcal{E%
}^{m}$ , so it kills all the nonintermediate cohomology of $\mathcal{E}^{m}$
, i.e. $\epsilon H^{i}(\mathcal{E}^{m})=0$ for $i\neq m+1$ , and $\epsilon
H^{m+1}(\mathcal{E}^{m})=H^{1}(\mathcal{L}_{0}^{t}\otimes \mathcal{L}%
_{0}^{n_{0}})$ so that $\epsilon $ fixes the cohomology clases lying in this
subspace. We define the cycle of $\widehat{A}_{Y(N_{2})}^{U,n}\times 
\mathcal{E}^{n_{0}}$ transform 
\begin{equation*}
\Delta _{n,n_{0}}:=(\mathcal{\epsilon }\times \epsilon _{\mathcal{E}%
^{n_{0}}})(\mathcal{G}r_{m,n_{0}})
\end{equation*}%
of the graph $\mathcal{G}r_{m,n_{0}}\subseteq \mathcal{E}^{m}\times \mathcal{%
E}^{n_{0}}$ of the projection 
\begin{equation*}
\varphi _{m;n_{0}}:\mathcal{E}^{m}\longrightarrow \mathcal{E}^{n_{0}}
\end{equation*}
to the $n_{0}$ arguments in $(A-T)\sqcup $ $(B-T)$. It vanishes on the
homology of $\mathcal{E}^{m}\times \mathcal{E}^{n_{0}}$ of order different
of its intermediate dimension $(m+$ $1)+(n_{0}+1)$ . In particular, it
vanishes on the homology class of $\mathcal{G}r_{m,n_{0}}\approx \mathcal{E}%
^{m}$ as it has topological dimension $2(m+1)$, thus differing $t>0$ from
the intermediate dimension.

\bigskip The other ingredient for the computation of the $p-$adic
Abel-Jacobi map will be a cuspidal form $g$ of weight $k_{0}=n_{0}+2\geq 2$
and $\Gamma _{0}-$level $N$, thus modular for the group $\Gamma (N)$, and a $%
\mathfrak{a}$-Hilbert cuspidal form $f$ of level $N_{2}$ and weight $%
(k,k)=(n+2,n+2)\geq (2,2)$ \ to which we associate the $\mathcal{L}^{2n}-$%
differential, 
\begin{align}
\omega _{f}& \in \text{$Fil$}^{2+2n}H_{dR}^{2+2n}(\widehat{A}%
_{Y(N_{2})}^{U,n})=H^{0}(\omega _{\widehat{A}_{Y(N_{2})}^{U,n}})=  \notag \\
H^{0}(pr^{\ast }\omega _{Y(N_{2})}\otimes \omega _{\widehat{A}%
_{Y(N_{2})}^{U.n}/Y(N_{2})})& =H^{0}(\omega _{Y(N_{2})}\otimes pr_{\ast
}(\omega _{\widehat{A}_{Y(N_{2})}^{U}/Y(N_{2})})^{n})  \notag \\
& =H^{0}(\omega _{Y(N_{2})}\otimes L^{(n,n)}).
\end{align}

\textbf{Remark. }It will be used in the proof of the proposition below that
the making the product $\mathcal{L}_{0}\mathcal{\otimes \mathcal{L}}%
_{0}\longrightarrow \mathcal{O}_{X(N)}(-1)$ in each of the $t$ fibred\
diagonal factors of $\mathcal{E}^{m}\hookrightarrow \mathcal{E}^{2n}$, an $%
\mathcal{L}_{0}^{2n}-$differential in $X(N)$ , thus differential of order $%
2n+1$ in $\mathcal{E}^{2n}$, transforms, after restriction to $\mathcal{E}%
^{m}\hookrightarrow \mathcal{E}^{2n}$ an equally denoted $\mathcal{L}%
_{0}^{n_{0}}$ -differential in $X(N)$ ,so differential of order $n_{0}+1$ in 
$\mathcal{E}^{n_{0}}$, via%
\begin{equation}
\mathcal{L}_{0}^{n}\otimes \mathcal{L}_{0}^{n}\hookrightarrow \mathcal{L}%
_{0}^{n_{0}/2}\otimes (\mathcal{L}_{0}\mathcal{\otimes \mathcal{L}}_{0}%
\mathcal{\mathcal{)}}^{\otimes t}\otimes \mathcal{L}_{0}^{n_{0}/2}\overset{}{%
\longrightarrow \mathcal{L}_{0}^{n_{0}/2}\otimes \mathcal{L}%
_{0}^{n_{0}/2}(-t)\overset{sym}{\longrightarrow }}\mathcal{L}_{0}^{n_{0}}(-t)
\label{desapareciendo}
\end{equation}%
(as in proof of prop. 2.9 in \cite{DR}).

\begin{proposition}
\textbf{\ \label{Clave 2 en Lausana}}

a) The $p-$adic Abel-Jacobi map is computed in terms of the product 
\begin{equation*}
<\_.\_>:H_{dR-par}^{1}(X(N),\mathcal{L}_{0}^{n_{0}}(-t)\times
H_{dR-par}^{1}(X(N),\mathcal{L}_{0}^{n_{0}})\longrightarrow \mathbb{C}%
_{p}(-1-m)
\end{equation*}%
by 
\begin{equation}
AJ_{p}(\Delta _{n,n_{0}})(\omega _{f}\otimes \eta _{g}^{u-r})=<Q_{f}(\phi
)^{-1}e_{g}e_{ord}j_{X^{\prime }(N)}^{\ast }\varrho _{\mathcal{F^{\sharp }}%
},\eta _{g}^{u-r}>  \label{AJ 1}
\end{equation}%
for the $\mathcal{L}^{2n}$-differential $\varrho _{\mathcal{F^{\sharp }}}$
with $\nabla ^{GM}(\varrho _{\mathcal{F^{\sharp }}})=\omega _{f^{\sharp }}$.

b) It can also be computed as 
\begin{equation*}
AJ_{p}(\Delta _{n,n_{0}})(\omega _{f}\otimes \eta _{g}^{u-r})=<Q_{f}(\phi
)^{-1}e_{g}e_{ord}j_{X^{\prime }(N)}^{\ast }sp(\varrho _{\mathcal{F^{\sharp }%
}}),\eta _{g}^{u-r}>
\end{equation*}
\end{proposition}

\bigskip \textbf{Proof:} a) Besser theory is applied similarly as in section
3, generalizing the polynomial $P_{f}(x)$ used in that computation by 
\begin{equation*}
P_{f}(x):=(1-p^{1-k}\phi )Q_{f}(x)\in \mathbb{Q}_{p}[x]
\end{equation*}

By the same reason as in section 3, this polynomial kills the target of the
residue map, so that (\ref{AJ 1}) equals\ \ $<P_{f}(\phi
)^{-1}e_{g}e_{ord}j_{X^{\prime }(N)}^{\ast }\varrho _{\mathcal{\mathcal{%
F^{\flat }}}},\eta _{g}^{u-r}>$ for the overconvergent $\mathcal{L}^{2n}$%
-differential $\varrho _{\mathcal{F}^{\flat }}$with $\nabla ^{GM}(\varrho _{%
\mathcal{F^{\flat }}})=\omega _{f^{\flat }}$ , which is equivalent to $%
d\varrho _{\mathcal{F^{\flat }}}=\omega _{f^{\flat }}$ as $\varrho _{%
\mathcal{F^{\flat }}}$ and $\omega _{f^{\flat }}$ can be seen as
overconvergent differentials in $\widehat{A}_{Y(N_{2})}^{U,n}$ of order $%
2n+1 $ and $2n+2$.

To prove \ref{AJ 1} , we can assume $t>0$ because this has been already
proved for $n=n_{0}=2$. As analogous generalization of this case, and just
as in the proof of the analogous in theor. 3.8 of \cite{DR}, the $p-$adic
Abel Jacobi map at the null-homologous cycle $\Delta _{n,n_{0}}=(\epsilon
\times \epsilon _{\mathcal{E}^{n_{0}}})\mathcal{E}^{m}$ of components
isomorphic to $\mathcal{E}^{m}$ and defined over $\mathbb{Z}_{p}$ \ (by \ref%
{alternative Scholl operator}), is given by the integral \ref{Abel-Jacobi} 
\begin{equation*}
AJ_{p}(\Delta _{n,n_{0}})(\omega _{f}\otimes \eta
_{g}^{u-r})=<cl_{fp}(\Delta _{n,n_{0}}),\widetilde{\omega }_{f}\otimes 
\widetilde{\eta }_{g}^{u-r}>_{fp}
\end{equation*}%
in terms of liftings of the given differentials by the Besser epimorphism. A
lifting $\widetilde{\omega }_{f}$ is given by the primitive $\varrho _{%
\mathcal{F^{\flat }}}$ of $\omega _{f^{\flat }}=P(\phi )\omega _{f}$ , and
the product equals 
\begin{equation*}
tr_{\mathcal{E}_{\mathbb{Z}_{p}}^{m}}(\varphi _{m;2n,n_{0}}^{\ast }(%
\widetilde{\omega }_{f}\otimes \widetilde{\eta }_{g}^{u-r}))=<\varphi
_{m;2n}^{\ast }(\widetilde{\omega }_{f}),\varphi _{m;n_{0}}^{\ast }(%
\widetilde{\eta }_{g}^{u-r})>_{fp}
\end{equation*}%
by the same argument as decomposition (86) \cite{DR} in components of $%
\Delta _{n,n_{0}}$ (here $\varphi _{m;2n,n_{0}}:\mathcal{E}%
^{m}\longrightarrow \mathcal{E}^{2n}\times \mathcal{E}^{n_{0}}$ is the
product of $\varphi _{m;2n}$ and $\varphi _{m;n_{0}}$ , defined over $%
\mathbb{Z}_{p}$). Just as in the argument of section 3, a preimage of the
finite polynomial class $\varphi _{m;2n}^{\ast }(\widetilde{\omega }_{f})$
for the map $\mathbf{i}$ of \ref{Besser exact sequence} is the de Rham class
of $P_{f}(\phi )^{-1}j_{X^{\prime }(N)}^{\ast }\varrho _{\mathcal{\mathcal{%
F^{\flat }}}}$so that 
\begin{eqnarray*}
AJ_{p}(\Delta _{n,n_{0}})(\omega _{f}\otimes \eta _{g}^{u-r}) &=&<P_{f}(\phi
)^{-1}e_{g}e_{ord}j_{X^{\prime }(N)}^{\ast }\varrho _{\mathcal{\mathcal{%
F^{\flat }}}},\eta _{g}^{u-r}>_{dR} \\
&=&<Q_{f}(\phi )^{-1}e_{g}e_{ord}j_{X^{\prime }(N)}^{\ast }\varrho _{%
\mathcal{F^{\sharp }}},\eta _{g}^{u-r}>_{dR}
\end{eqnarray*}%
proving a) in this case. To make sense of this product, observe that
restriction to $\mathcal{E}^{m}$ of the primitive $\varrho _{\mathcal{%
F^{\sharp }}}$ of $\omega _{f^{\sharp }}=P(\phi )\omega _{f}$ appearing in
the its factor corresponds to an overconvergent section of $\mathcal{L}%
_{0}^{n}\otimes \mathcal{L}_{0}^{n}$ (as class acted by $\epsilon $ thus in
its image $\mathcal{L}_{0}^{n}\otimes \mathcal{L}_{0}^{n}$ ), and it is then
pulled-back by $pr_{m;2n}:$ $\mathcal{E}^{m}\longrightarrow \mathcal{E}^{2n}$
to a differential of $\mathcal{E}^{m}$ or section of $\mathcal{L}_{0}^{m}$;
and observe that the second factor is an overconvergent \ section of $%
\mathcal{L}_{0}^{n_{0}}$ or overconvergent differential of $\mathcal{E}%
^{n_{0}}$, pulled-back by $pr_{m;n_{0}}:$ $\mathcal{E}^{m}\longrightarrow 
\mathcal{E}^{n_{0}}$ to a section of $\mathcal{L}^{m}$ or differential of $%
\mathcal{E}^{m}$. The product of a de Rham class in $\mathcal{E}^{m}$ and
the de Rham class $\eta _{g}^{u-r}$ in $\mathcal{E}^{n_{0}}$ makes sense
because, by projection formula, it can be equivalently understood with
pull-back by $pr_{m;n_{0}}$ in the second argument or with push-forward in
the first argument.

As for the remaining case $n_{0}=0$, $n=1$, two analogous terms to those (%
\ref{sustracting terms}) in the proof of (\ref{lema computo por Poincare})
must be subtracted, and only one term if $n_{0}=0,n>1$, but in both cases,
and by analogous reason, they do not contribute to the computation of
product \ref{AJ 1}.

b) We decompose $\varrho _{\mathcal{F^{\sharp }}}=\varrho _{n}+\nabla s$
with $\varrho _{n}\approx sp(\varrho )$\ by (\ref{Clave en Lausanne}), and
remark that $\varrho _{n}$ is also a primitive of $\varrho $ so that it can
be replaced en (\ref{AJ 1}), so also $sp(\varrho )$ , instead of $\varrho _{%
\mathcal{F^{\sharp }}}$ $\square $

\bigskip

\bigskip \textbf{Remark.} Observe that we have not used as polynomial $%
P_{f}(x)$ the analogous to the one considered by \cite{DR} (to fulfill
condition ii) of subsection 3.2) as this would be $(1-p^{-\frac{m+2}{2}%
}x)Q_{f}(x)=(1-p^{n+1}x)Q_{f}(x)$. As some of his roots have complex norm $%
p^{\frac{m+2}{2}}=p^{\frac{t+n_{0}+2}{2}}$ and some have $p^{\frac{2n+2}{2}%
}=p^{k-1}$, so not all of them of the same complex norm, it would not be
allowed for our use in Besser theory. The first factor is used in order to
vanish, while evaluated in $\phi $, in the target of the residue map for the
first argument of the cup product $\ref{AJ 1}$ . Since this argument lies in 
$H_{dR-par}^{1}(X(N),\mathcal{L}_{0}^{n_{0}}(-t))$ the target of its residue
map is twisted by $\mathbb{Q}_{p}(-\frac{n_{0}+2}{2}-t)=\mathbb{Q}_{p}(-n-1)$
, so the factor $(1-p^{-\frac{2n+2}{2}}x)=(1-p^{1-k}x)$ we consider in our
definition of $P_{f}(x)$ does the job.

\bigskip

\subsection{\textbf{Primitives for depleted Hilbert modular forms}}

Keeping our notations, if an overconvergent $\mathfrak{a}$-Hilbert cuspidal
form $f$ of weight $(k,k^{\prime })$ , level $N_{2}$ , defined over $\mathbb{%
Q}_{p}$, has $q-$expansion (\ref{expansion})%
\begin{equation}
f(q)=\dsum\limits_{\nu \in (\mathfrak{a}^{-1})^{+}}^{{}}a_{\nu }q^{\nu }%
\text{ ,}  \label{q-expansion}
\end{equation}%
its $p-$depletion and $\pi -$depletion are the overconvergent $\mathfrak{a}$%
-Hilbert modular forms 
\begin{equation*}
f^{[p]}=(1-V_{p}U_{p})f\text{ \ \ and }f^{[\pi ]}=(1-V_{\pi }U_{\pi })f\text{
}
\end{equation*}%
of same weight and level, and $q-$expansion%
\begin{equation}
f^{[p]}(q)=\dsum\limits_{\nu \in (\mathfrak{a}^{-1})^{+}\text{, }p\nmid \nu
}^{{}}a_{\nu }q^{\nu }\text{ \ and }f^{[\pi ]}(q)=\dsum\limits_{\nu \in (%
\mathfrak{a}^{-1})^{+}\text{, }\pi \nmid \nu }^{{}}a_{\nu }q^{\nu }\text{.}
\label{depleted expansion}
\end{equation}%
An overconvergent $\mathfrak{a}$-Hilbert modular form $f$ is said $p$%
-depleted if $f=f^{[p]}$ and $\pi -$depleted if $f=f^{[\pi ]}$. Obviously
these are idempotent operators, and $\pi -$depleted implies $p-$depleted. We
assume from now on that $f$ is nonordinary at both $\pi $ and $\pi ^{\prime
} $.

\textbf{Proposition}. Under this nonordinarity condition for cuspidal $f$ of
weight $(k,k^{\prime })=(n+2,n^{\prime }+2)$,

a) there is a unique $\mathcal{L}^{n+n^{\prime }}-$valued $(1,0)-$%
differential $\varrho _{\mathcal{F}^{[\pi ^{\prime }]}}\in \Omega _{\mathcal{%
W}_{\varepsilon }}^{(1,0)}\otimes \mathcal{L}_{\mathcal{W}_{\varepsilon
}}^{n+n^{\prime }}\subseteq L_{\mathcal{W}_{\varepsilon }}^{(2,0)}\otimes 
\mathcal{L}_{\mathcal{W}_{\varepsilon }}^{n+n^{\prime }}$ such that%
\begin{equation*}
\nabla ^{GM}(\varrho _{\mathcal{F}^{[\pi ^{\prime }]}})=\nabla ^{\prime
}\varrho _{\mathcal{F}^{[\pi ^{\prime }]}}=\omega _{f^{[\pi ^{\prime }]}}
\end{equation*}%
and is rational sum of $n^{\prime }$ nearly overconvergent forms $(\theta
^{\prime })^{-j-1}f^{[\pi ^{\prime }]}$ in the sense that, near the standard
cusp, it trivializes (in the above base trivializing $\mathcal{L}%
^{n+n^{\prime }}$) as sum of $q-$expansions 
\begin{equation}
\mathcal{F}^{[\pi ^{\prime }]}(q)=\dsum_{j=0}^{n^{\prime }}(-1)^{j}j!\binom{%
n^{\prime }}{j}(\theta ^{\prime })^{-j-1}f^{[\pi ^{\prime }]}(q)\omega
^{n}(\omega ^{\prime })^{n^{\prime }-j}(\eta ^{\prime })^{j}\text{ .}
\label{expansion of scriptF primitive}
\end{equation}

b) There is a unique nearly convergent $L^{(n,n^{\prime })}-(1,0)-$%
differential $\varrho (F_{[k,k^{\prime }-2]}^{[\pi ^{\prime }]})$ associated
to overconvergent cuspidal form $F_{[k,k^{\prime }-2]}^{[\pi ^{\prime }]}$
of weight $(k,k^{\prime }-2)=(n+2,n^{\prime })$\ such that 
\begin{equation*}
\nabla _{sp}^{GM}\varrho (F_{[k,k^{\prime }-2]}^{[\pi ^{\prime }]})=\nabla
_{sp}^{\prime }\varrho (F_{[k,k^{\prime }-2]}^{[\pi ^{\prime }]})=\omega
_{f^{[\pi ^{\prime }]}}
\end{equation*}%
\ \ \ \ \ \ \ \ \ \ \ \ \ \ \ \ \ \ \ \ \ \ \ thus $F_{[k,k^{\prime
}-2]}^{[\pi ^{\prime }]}$ having $q-$expansion%
\begin{equation}
F_{[k,k^{\prime }-2]}^{[\pi ^{\prime }]}(q)=\dsum_{\nu \in (\mathfrak{a}%
^{-1})^{+}\text{, }\pi ^{\prime }\nmid \nu }(\nu ^{\prime })^{-1}a_{\nu
}q^{\nu }  \label{primitive depleted expansion}
\end{equation}

i.e. $\theta ^{\prime }(F_{[k,k^{\prime }-2]}^{[\pi ^{\prime }]}=f^{[\pi
^{\prime }]}$ \ , so that we write (because of the uniqueness) $%
F_{[k,k^{\prime }-2]}^{[\pi ^{\prime }]}$ $=(\theta ^{\prime })^{-1}f^{[\pi
^{\prime }]}$

\textbf{Proof: \ }a) The purely formal fact that \ref{expansion of scriptF
primitive} is the $q-$expansion of the $\nabla ^{GM}-$primitive if it
exists, thus the unicity, is an obvious computation which goes exactly as
for depleted modular forms, since we only derive one of the two variables
(this is made in detail\textbf{\ }in comment preceding lemma 9.2 of \cite%
{Coleman95}\ \textbf{,} in the particular case where $f_{i}=0$ in notations
of Loc. cit.) The essential fact of its existence is proved, under the
nonordinarity hypothesis of $f$ at both $\pi $ and $\pi ^{\prime }$ , in 
\cite{David et al.} (this nonordinarity condition was not imposed in a
previous version of our article, and is in fact unnecesary according a
conjecture in Loc. cit.). In particular, it is proved in Loc. cit. that $%
(\theta ^{\prime })^{-n^{\prime }-1}f^{[\pi ^{\prime }]}$ is overconvergent,
so it follows that its iterated $\theta ^{\prime }$ deriviatives $(\theta
^{\prime })^{-j-1}f^{[\pi ^{\prime }]}$ , for all $j=0,...,n^{\prime }-1$ ,
are all nearly-overconvergent, by lemma \ref{Conservacion near converg}.

b) It follows from a) by just taking $\varrho _{F_{[k,k^{\prime }-2]}^{[\pi
^{\prime }]}}=sp(\varrho _{\mathcal{F}^{[\pi ^{\prime }]}})$ , i.e. the term 
$j=0$ in the sum (\ref{expansion of scriptF primitive})\ of $q-$expansions.$%
\square $

Observe that the coefficients of the expansion (\ref{primitive depleted
expansion}) stay all in $\mathbb{Z}_{p}$ as $\pi \nmid \nu ^{\prime }$ in
all of them, so that $\nu ^{\prime }$ is invertible in $\mathfrak{o}_{\pi
}\approx \mathbb{Z}_{p}$ .

\bigskip \textbf{Remark}. Applying (\ref{desapareciendo}) to the $\mathcal{L}%
_{0}^{2n}$ -differential $j_{X^{\prime }(N)}^{\ast }\mathcal{\varrho }(%
\mathcal{F}^{[\pi ^{\prime }]})$ in (\ref{expansion of scriptF primitive}) ,
we obtain an equally denoted overconvergent differential in $\mathcal{L}%
_{0}^{n_{0}}(-t)$ whose split projection is the nearly overconvergent $%
(-1)^{t}t!j_{X^{\prime }(N)}^{\ast }\varrho (\theta ^{\prime -1-t}f^{[\pi
^{\prime }]})$, thus with an overconvergent $e_{ord}$ projection which we
will prove to compute the $p-$ adic Abel-Jacobi map for higher dimensional
null-homologous cycles.

\section{\textbf{Computation of the }$p-$ \textbf{adic Abel--Jacobi map}}

Let us recall first how the isotypic component and the stabilizations are
made for cuspidal forms, so to imitate it for $\mathfrak{a}$-Hilbert
cuspidal forms. Assume that $g$ \textit{is a cuspidal form over }$\mathbb{Q}$%
\textit{\ of weight }$k_{0}=n_{0}+2$ with $n_{0}\geq 0$\textit{\ for the
modular group} group $\Gamma (N)$, and consider it as defined over $\mathbb{Q%
}_{p}$, for $p\nmid N$ . Let%
\begin{equation*}
\dsum\limits_{n\geq 1}^{{}}b_{n}q^{n}
\end{equation*}%
be its $q-$expansion. \textit{Assume also that }$g$\textit{\ is an eigenform
of eigenvalue} $b_{n}$ \textit{for all the Hecke operators} $T_{n}$, \textit{%
and that it is normalized} , so that all $b_{n}\in \mathbb{Z}.$\textbf{\ }In
particular $b_{p}g=T_{p}g=U_{p}g+p^{k_{0}-1}V_{p}g$ . As a consequence, and
using $U_{p}V_{p}=1$, the operator\textbf{\ } 
\begin{equation}
1-V_{p}b_{p}+p^{k_{0}-1}V_{p}^{2}=(1-\beta _{0}V_{p})(1-\beta _{1}V_{p})%
\text{ with }\beta _{0}+\beta _{1}=b_{p}\text{ and }\beta _{0}\beta
_{1}=p^{k_{0}-1}  \label{polinomio de grado 2}
\end{equation}%
acts on the cuspidal form $g$ as the $p-$depletion operator $1-V_{p}U_{p}$,
where $\beta _{0},\beta _{1}\in \overline{\mathbb{Q}}$ are the reciprocal
roots of this polynomial in $V_{p}$, i.e. $\beta _{0}^{-1}=p^{1-k_{0}}\beta
_{1}$ and $\beta _{1}^{-1}=p^{1-k_{0}}\beta _{0}$ are its roots. Again we
will not work in $X(N)$ with all the $\mathcal{L}_{0}-$valued $\log D^{c}$
-de Rham cohomology, but with its parabolic piece 
\begin{equation*}
H_{dR-par}^{1}(X(N),\mathcal{L}_{0}^{n_{0}}/\mathbb{Q}_{p})\subseteq
H_{dR}^{1}(X(N),\mathcal{L}_{0}^{n_{0}}/\mathbb{Q}_{p})
\end{equation*}%
which inherits from de Rham an equally denoted product $<-,->$, and an
action of $V_{p}$ as $p^{1-k_{0}}\phi $ , having $U_{p}$ as inverse (as $%
U_{p}V_{p}=1$ in this finite-dimensional space). The polynomial \ref%
{polinomio de grado 2} in $V_{p}$ vanishes on the class of the $%
L_{0}^{n_{0}} $-valued differential $\omega _{g}$ , so $V_{p}$ leaves
invariant a $2$-dimensional subspace or isotypic component of $g$ 
\begin{equation*}
H_{dR-par}^{1}(X(N),\mathcal{L}_{0}^{n_{0}}/\mathbb{Q}_{p})(g)\subseteq
H_{dR-par}^{1}(X(N),\mathcal{L}_{0}^{n_{0}}/\mathbb{Q}_{p})
\end{equation*}%
the inclusion following from the fact that the parabolic cohomology is $%
V_{p} $ -invariant.

\bigskip

On the $\mathbb{\mathbb{C}}_{p}$-extended isotypic component $%
H_{dR-par}^{2}(X_{\mathbb{C}_{p}}(N),\mathcal{L}_{0}^{n_{0}}/\mathbb{C}%
_{p})(g)$ , 
\begin{align}
\text{ }V_{p}& =p^{1-k_{0}}\phi \text{ acts as }\left( 
\begin{array}{cc}
p^{1-k_{0}}\beta _{1} &  \\ 
& p^{1-k_{0}}\beta _{0}%
\end{array}%
\right)  \notag \\
\text{ }U_{p\text{ \ }}& =V_{p}^{-1}\text{ acts as }\left( 
\begin{array}{cc}
\beta _{0} &  \\ 
& \beta _{1}%
\end{array}%
\right) \text{ \ .}  \label{diagonalizacion debil de Up}
\end{align}

We define the two $p-$stabilizations $g_{i}$ of $g$ , $i\in \{0,1\}$, as 
\begin{equation}
g_{i}:=(1-\beta _{j}V_{p})g\in S_{k_{0}}^{oc}(\Gamma (N),\mathbb{C}_{p})%
\text{\ for }j\in \{0,1\}\text{ with }j\neq i
\end{equation}

and decompose $g$ as linear combination of these two stabilizations%
\begin{equation}
g=\frac{\beta _{0}}{\beta _{0}-\beta _{1}}g_{0}-\frac{\beta _{1}}{\beta
_{0}-\beta _{1}}g_{1}\text{ .}  \label{g suma de dos}
\end{equation}%
Observe that $g_{i}$ has still eigenvalue $b_{n}$ for all operators $T_{n}$
with $p\nmid n$ . Indeed, 
\begin{equation*}
T_{n}g_{i}=T_{n}g-\beta _{j}V_{p}T_{n}g=b_{n}(g-\beta _{j}V_{p}g)=b_{n}g_{i}%
\text{ ,}
\end{equation*}%
since $V_{p}$ commutes with $U_{n}$ , thus also with $T_{n}$. Observe also
that $g_{i}$ is an eigenform for $U_{p}$ of eigenvalue $\beta _{i}$. Indeed,
since $T_{p}g=b_{p}g=\beta _{i}g+\beta _{j}g$, it is 
\begin{equation*}
U_{p}g_{i}=U_{p}g-\beta _{j}g=T_{p}g-p^{k_{0}-1}V_{p}g-\beta _{j}g=\beta
_{i}g-p^{k_{0}-1}V_{p}g=\beta _{i}g_{i}\text{ ,}
\end{equation*}%
the last equality following from $\beta _{i}\beta _{j}=p^{k_{0}-1}$. The
decomposition \ref{g suma de dos} tells us that $g$ belongs to the
2-dimensional span 
\begin{equation*}
S_{k_{0}}^{oc}(\Gamma (N),\mathbb{\mathbb{C}}_{p})(g)\subseteq
S_{k_{0}}^{oc}(\Gamma (N),\mathbb{\mathbb{C}}_{p})\text{ }
\end{equation*}%
of $g_{0}$ and $g_{1}$, a subspace acted by $U_{p}$ as 
\begin{equation}
U_{p}=\left( 
\begin{array}{cc}
\beta _{0} &  \\ 
& \beta _{1}%
\end{array}%
\right)  \label{diagonalizacion fuerte de Up}
\end{equation}%
(in the strong sense, not only in the sense of cohomology classes, what we
already had in \ref{diagonalizacion debil de Up}). Clearly, the natural
homomorphism 
\begin{equation*}
S_{k_{0}}^{oc}(\Gamma (N),\mathbb{\mathbb{Q}}_{p})\longrightarrow
H_{rig-par}^{1}(X^{\prime }(\Gamma (N)),\mathcal{L}_{0}^{n_{0}}/\mathbb{Q}%
_{p})
\end{equation*}%
restricts to an isomorphism 
\begin{eqnarray}
S_{k_{0}}^{oc}(\Gamma (N),\mathbb{Q}_{p})(g) &\approx &
\label{two dimensional isotypic} \\
H_{dR-par}^{1}(X(N),\mathcal{L}_{0}^{n_{0}}/\mathbb{Q}_{p})(g) &\subseteq
&H_{rig-par}^{1}(X^{\prime }(N),\mathcal{L}_{0}^{n_{0}}/\mathbb{Q}_{p})\text{%
,}  \notag
\end{eqnarray}%
the inclusion (up to isomorphism) being analogous to \ \ref{cuadrado magico}%
. \textit{Assume now that }$g$\textit{\ is ordinary for} $p$, i.e. $U_{p}$
acts on it via a $p-$adic unit $\beta _{0}$. This is equivalent to saying
that $U_{p}$ acts via this unit on the cohomology class of its associated $%
L_{0}^{n_{0}}-$valued differential $\omega _{g}$ . Therefore, its inverse
operator $V_{p}$ in cohomology acts also via a $p-$adic unit $\beta
_{0}^{-1}=p^{1-k_{0}}\beta _{1}$ on this class, i.e. $\phi =p^{k_{0}-1}V_{p}$
acts with algebraic eigenvalue $\beta _{1}$ of $p-$adic value $k_{0}-1$, and
this is said "acting with slope $k_{0}-1$". In other words, $g$ lies in the
first eigenspace of the diagonalization \ref{diagonalizacion fuerte de Up},
with $\beta _{0}$ being a $p-$adic unit. Since both algebraic eigenvalues $%
\beta _{0}$ and $\beta _{1}$ can be distinguished, in the ordinary case, by
having $p-$adic value $0$ or $k_{0}-1$ (which explains our notation, taken
from the case $k_{0}=2$), we could properly write $\beta _{0}(g)$ and $\beta
_{1}(g)$, but spare this for convenience except when we want to emphasize
it. We then denote%
\begin{eqnarray}
&&H_{dR-par}^{1}(X_{\mathbb{C}_{p}}(N),\mathcal{L}_{0}^{n_{0}})^{\phi
,k_{0}-1}(g)\oplus H_{dR-par}^{1}(X_{\mathbb{C}_{p}}(N),\mathcal{L}%
_{0}^{n_{0}})^{\phi ,0}(g)  \notag \\
&=&H^{ord}(g)\oplus H^{u-r}(g)  \label{descompongo en dos monodim2}
\end{eqnarray}%
the decomposition of the $\mathbb{C}_{p}-$isotypic component of $g$ in the
two monodimensional\ ("ordinary" and "unit root) subspaces diagonalizing (%
\ref{diagonalizacion debil de Up}), as already recalled for the special case 
$k_{0}=2$ in (\ref{descompongo en dos monodim2}). The class of the
corresponding $L_{0}^{n_{0}}$-valued differential $\omega _{g}$ lies in the
first direct factor. Since the the $\mathbb{C}_{p}$-valued product
establishes a perfect pairing between the two factors, \textit{there is a
unique\ class} $\eta _{g}^{u-r}$ \textit{in the second direct factor such
that}%
\begin{equation}
<\omega _{g},\eta _{g}^{u-r}>=1\text{ ,}  \label{unit-root in higher weights}
\end{equation}%
(This class is presented in \cite{DR} as lift of the class in\textbf{\ }$%
H^{1}(X(N),L_{0}^{-n_{0}})$ represented by the $L_{0}^{-n_{0}}$-valued
antiholomorphic differential which is hermitian-dual to $\omega _{g}$ $\in
H^{0}(X(N),\omega _{X(N)}\otimes L_{0}^{n_{0}})$).

Observe that this works also if the ordinary cuspidal form $g$ of level $N$
is not classical but only overconvergent, thus classical of level $Np$
(Coleman classicality theorem), except that its isotypic component in $%
H_{dR-par}^{2}(X^{\prime }(N),\mathcal{L}_{0}^{n_{0}}/\mathbb{Q}_{p})$ does
not lie then inside $H_{dR-par}^{2}(X(N),\mathcal{L}_{0}^{n_{0}}/\mathbb{Q}%
_{p})$ as in (\ref{two dimensional isotypic}).

Denote still by $e_{ord}$ the action of the Hida projector induced on
cohomology. The obvious equation $e_{ord}(g_{0})=e_{ord}g-\beta
_{1}e_{ord}V_{p}g$ \ induces the\ equation in cohomology%
\begin{eqnarray}
e_{ord}(g_{0}) &=&e_{ord}g-\beta _{1}Ve_{ord}g=\mathcal{E}_{0}(g)e_{ord}g%
\text{ ,}  \label{Eps_0} \\
\text{for }\mathcal{E}_{0}(g) &:&=(1-\beta _{1}\beta _{0}^{-1})=(1-\beta
_{1}^{2}p^{1-k_{0}})\text{ ,}  \notag
\end{eqnarray}%
since $V_{p}$ commutes in cohomology with $U_{p}$, so also with $e_{ord}$ .

Since we are assuming that $g$ is ordinary at $p$, it makes sense to \textit{%
call }$g_{0}$ \textit{the ordinary }$p-$\textit{stabilization} $g^{(p)}$ 
\textit{of} $g$, so the equality (\ref{Eps_0}) becomes $e_{ord}g^{(p)}=%
\mathcal{E}_{0}(g)g$ .

Next, we describe next similar notions for an $\mathfrak{a}$-Hilbert modular
form $f$ \ defined over $R$, thus over $\mathbb{Q}_{p}$, of parallel weight $%
(k,k)\geq (2,2)$ and level $N_{2}\geq 4$ which is eigenform for all Hecke
operators $T_{\nu }$, in particular eigenform for $T_{\pi }=$ $U_{\pi
}+p^{k-1}V_{\pi }$, i.e.

\begin{equation*}
T_{\pi }f=U_{\pi }f+p^{k-1}V_{\pi }f\text{ }=a_{\pi }f
\end{equation*}%
This decomposition -kind of Eichler relation since $\phi _{\pi \text{ }}$%
acts as $p^{k-1}V_{\pi }$ on the corresponding $L^{(k-2,k-2)}$-valued
differential $\omega _{f}$ - follows from the computation of the
coefficients of the $q-$expansion of the $T_{\pi }$ transform (found for
instance in the discussion before (3) in VI. 1 of \cite{VdG}). Since $U_{\pi
}V_{\pi }=1$, the operator

\begin{equation}
1-a_{\pi }V_{\pi }+p^{k-1}V_{\pi }^{2}=(1-\alpha _{0}V_{\pi })(1-\alpha
_{1}V_{\pi })\text{ with }\alpha _{0}+\alpha _{1}=a_{\pi }\text{ and }\alpha
_{0}\alpha _{1}=p^{k-1}  \label{polynomial in Vpi}
\end{equation}%
acts on $f$ as $\pi -$depletion operator $1-V_{\pi }U_{\pi }$, where the
algebraic values $\alpha _{0},\alpha _{1}$ are the reciprocal roots of this
polynomial in $V_{\pi }$, i.e. $\alpha _{0}^{-1}=p^{1-k}\alpha _{1}$ and $%
\alpha _{1}^{-1}=p^{1-k}\alpha _{0}$ are the roots; and $\alpha _{0^{\prime
}},\alpha _{1^{\prime }}$ are defined analogously for the prime $\pi
^{\prime }$ . Recall we are under the nonordinarity hypothesis that all $%
\alpha _{0},\alpha _{1},\alpha _{0^{\prime }},\alpha _{1^{\prime }}$ are $p-$%
adic nonunits.

The two "$\pi $-stabilizations" $f_{i}\in $ $S_{k,k}^{oc}(N_{2},\mathbb{C}%
_{p})$ of $f$ are 
\begin{align}
f_{i}& :=(1-\alpha _{j}V_{\pi })f\text{ \ for }j\neq i\text{ ,} \\
\text{i.e. }f_{i}& =\dsum\limits_{\nu \in (\mathfrak{a}^{-1})^{+}}^{{}}c_{%
\nu }q^{\nu }\text{ with }c_{\nu }=a_{\nu }-\alpha _{j}a_{\nu /\pi }\text{ ,
if }\pi \mid \nu \text{, and }a_{\nu }\text{ otherwise.}  \notag
\end{align}%
Observe that $f_{i}$ has still eigenvalue $a_{\nu }$ for the Hecke operators 
$T_{\nu }$ with $\pi \nmid \nu $ (thus same eigenvalue $a_{\nu }$ for all
the Hecke operators $T_{\nu }$ with $\pi \nmid \nu $). Indeed, 
\begin{equation*}
T_{\nu }f_{i}=T_{\nu }g-\alpha _{j}V_{\pi }T_{\nu }f=a_{\nu }(g-\alpha
_{j}V_{\pi }f)
\end{equation*}%
since $V_{\pi }$ commutes with $U_{\nu }$ , thus also with $T_{\nu }$.
Observe furthermore that $f_{i}$ is an eigenform for $U_{\pi }$ of
eigenvalue $\alpha _{i}$. Indeed, since $T_{\pi }f=a_{\pi }f=\alpha
_{i}f+\alpha _{j}f$ , 
\begin{equation*}
U_{\pi }f_{i}=U_{\pi }f-\alpha _{j}f=T_{\pi }f-p^{k-1}V_{\pi }f-\alpha
_{j}f=\alpha _{i}f-p^{k-1}V_{\pi }f=\alpha _{i}f_{i}\text{ ,}
\end{equation*}%
the last equality because $\alpha _{i}\alpha _{j}=p^{k-1}$. We have then the
decomposition%
\begin{equation}
f=\frac{\alpha _{0}}{\alpha _{0}-\alpha _{1}}f_{0}-\frac{\alpha _{1}}{\alpha
_{0}-\alpha _{1}}f_{1}\text{ ,}
\end{equation}%
and the analogous 
\begin{equation*}
f=\frac{\alpha _{0^{\prime }}}{\alpha _{0^{\prime }}-\alpha _{1^{\prime }}}%
f_{0^{\prime }}-\frac{\alpha _{1^{\prime }}}{\alpha _{0^{\prime }}-\alpha
_{1^{\prime }}}f_{1^{\prime }}\text{ ,}
\end{equation*}%
in terms of the two $\pi ^{\prime }$-stabilizations $f_{i^{\prime }}$ having
eigenvalue $\alpha _{i^{\prime }}$ for $U_{\pi ^{\prime }}$ and eigenvalue $%
a_{\nu }$ for all $T_{\nu }$ $\pi ^{\prime }\nmid \nu $ .

We can define also the four "$\pi ,\pi ^{\prime }$-stabilizations" $%
f_{ii^{\prime }}\in $ $S_{k,k}^{oc}(N_{2},\mathbb{C}_{p})$ of $f$ as 
\begin{equation}
f_{ii^{\prime }}:=(1-\alpha _{i^{\prime }}V_{\pi ^{\prime }})(1-\alpha
_{i}V_{\pi })f=(1-\alpha _{i}V_{\pi })(1-\alpha _{i^{\prime }}V_{\pi
^{\prime }})f\text{ .}  \label{pi,pi'-estabilizaciones}
\end{equation}%
This is an eigenform for\ $U_{\pi }$ and $U_{\pi ^{\prime }}$ of eigenvalue $%
\alpha _{i}$ and $\alpha _{i^{\prime }}$, thus for $U_{p}=U_{\pi }U_{\pi
^{\prime }}$ of eigenvalue $\alpha _{i}\alpha _{i^{\prime }}$; and it is an
eigenform of $T_{\nu }$ of eigenvalue $a_{\nu }$ whenever $\pi ,\pi ^{\prime
}\nmid \nu $. In particular, the $\nu -$coefficients of the $q-$expansions
of both $f_{ii^{\prime }}$ and $f$ are the same, for $\pi ,\pi ^{\prime
}\nmid \nu $ .

From the two decompositions above, we obtain the four terms decomposition

\begin{equation}
f=\dsum\limits_{i,i^{\prime }}^{{}}(-1)^{i+i^{\prime }}\frac{\alpha
_{i}\alpha _{i^{\prime }}}{(\alpha _{0}-\alpha _{1})(\alpha _{0^{\prime
}}-\alpha _{1^{\prime }})}f_{ii^{\prime }}  \label{f como suma de cuatro}
\end{equation}%
in eigenforms of the $U_{p}$ operator, all having same eigenvalue $a_{\nu }$
for the operators$T_{\nu }$ with $\pi ,\pi ^{\prime }\nmid \nu $.

\bigskip Putting this together, $f$ belongs to the 2-dimensional span 
\begin{equation*}
S_{k,k}^{oc}(N_{2},\mathbb{C}_{p})(U_{\pi },f)\subseteq S_{k,k}^{oc}(N_{2},%
\mathbb{C}_{p})
\end{equation*}

of $f_{0}$and $f_{1}$ which is acted by $U_{\pi }$ as 
\begin{equation}
\left( 
\begin{array}{cc}
\alpha _{0} &  \\ 
& \alpha _{1}%
\end{array}%
\right) \text{ ,}
\end{equation}%
and to the span of $f_{0^{\prime }}$, $f_{1^{\prime }}$ acted by $U_{\pi
^{\prime }}$ as 
\begin{equation*}
\left( 
\begin{array}{cc}
\alpha _{0^{\prime }} &  \\ 
& \alpha _{1^{\prime }}%
\end{array}%
\right) \text{ ;}
\end{equation*}

\bigskip And on the other hand, $f$ belongs to the 4-dimensional span%
\begin{equation*}
S_{k,k}^{oc}(N_{2},\mathbb{C}_{p})(f)\subseteq S_{k,k}^{oc}(N_{2},\mathbb{C}%
_{p})
\end{equation*}

of $f_{00^{\prime}}$, $f_{10^{\prime}},$ $f_{10^{\prime}},$ $f_{11^{\prime}}$
, which is acted by $U_{p}=U_{\pi}U_{\pi^{\prime}}$ as 
\begin{equation}
\left( 
\begin{array}{cccc}
\alpha_{0}\alpha_{0^{\prime}} &  &  &  \\ 
& \alpha_{0}\alpha_{1^{\prime}} &  &  \\ 
&  & \alpha_{1}\alpha_{0^{\prime}} &  \\ 
&  &  & \alpha_{1}\alpha_{1^{\prime}}%
\end{array}
\right)
\end{equation}

\bigskip so that 
\begin{equation}
S_{k,k}^{oc}(N_{2},\mathbb{C}_{p})(f)\approx S_{k,k}^{oc}(N_{2},\mathbb{C}%
_{p})(U_{\pi},f)\otimes S_{k,k}^{oc}(N_{2},\mathbb{C}_{p})(U_{\pi^{%
\prime}},f)  \label{desc. tensorial de cuspidales}
\end{equation}

as $U_{p}=U_{\pi}U_{\pi^{\prime}}$ -acted space.

As a consequence, in $H_{dR-par}^{2}(Y(N_{2}),\mathcal{L}^{n}/\mathbb{C}%
_{p}) $, for $n=2k-2$, where $V_{\pi }$ acts as $p^{1-k}\phi _{\pi }$ and
has inverse $U_{\pi }$ (because $U_{\pi }V_{\pi }=1$ in this
finite-dimensional space), the polynomial ($\ref{polynomial in Vpi}$) in $%
V_{\pi }$ vanishes on the class of the $L^{(k-2,k-2)}$-valued differential $%
\omega _{f}$ , so $V_{\pi }$ leaves invariant a 2-dimensional subspace%
\begin{equation*}
H_{dR-par}^{2}(Y(N_{2}),\mathcal{L}^{n}/\mathbb{C}_{p})(\phi _{\pi
},f)\subseteq H_{dR-par}^{2}(Y(N_{2}),\mathcal{L}^{n}/\mathbb{C}_{p})\text{ ,%
}
\end{equation*}%
in which it acts with minimal polynomial (\ref{polynomial in Vpi}). On the $%
\mathbb{C}_{p}$ -extension $H_{dR-par}^{2}(Y_{\mathbb{C}_{p}}(N_{2}),%
\mathcal{L}^{n})(\phi _{\pi },f)$ of this subspace, 
\begin{equation}
V_{\pi }=p^{1-k}\phi _{\pi }=U_{p}^{-1}\text{ acts as }\left( 
\begin{array}{cc}
p^{1-k}\alpha _{1} &  \\ 
& p^{1-k}\alpha _{0}%
\end{array}%
\right)  \notag
\end{equation}%
and on the isotypic component $H_{dR-par}^{2}(Y(N_{2}),\mathcal{L}^{k-2}/%
\mathbb{C}_{p})(f)$ of $f$ \ in $H_{dR-par}^{2}(Y(N_{2}),\mathcal{L}^{n}/%
\mathbb{C}_{p})$, i.e. the minimal subspace invariant by $\phi =\phi _{\pi
}\phi _{\pi ^{\prime }}$ in which the class of $\omega _{g}$ lies,%
\begin{eqnarray*}
V_{p}\text{ } &=&V_{\pi }V_{\pi ^{\prime }}=U_{p}^{-1}=p^{2-2k}\phi \text{
acts as} \\
&&\text{ }\left( 
\begin{array}{cccc}
p^{2-2k}\alpha _{1}\alpha _{1^{\prime }} &  &  &  \\ 
& p^{2-2k}\alpha _{1}\alpha _{0^{\prime }} &  &  \\ 
&  & p^{2-2k}\alpha _{0}\alpha _{1^{\prime }} &  \\ 
&  &  & p^{2-2k}\alpha _{0}\alpha _{0^{\prime }}%
\end{array}%
\right)
\end{eqnarray*}

so that%
\begin{equation}
\begin{array}{c}
H_{dR-par}^{2}(Y(N_{2}),\mathcal{L}^{k-2}/\mathbb{C}_{p})(f)\approx \\ 
H_{dR-par}^{2}(Y(N_{2}),\mathcal{L}^{n}/\mathbb{C}_{p})(\phi _{\pi
},f)\otimes H_{dR-par}^{2}(Y(N_{2}),\mathcal{L}^{n}/\mathbb{C}_{p})(\phi
_{\pi },f)%
\end{array}
\label{desc. tensorial en cohomologia}
\end{equation}

\bigskip as $\phi _{p}=\phi _{\pi }\phi _{\pi ^{\prime }}$ -acted space.
Clearly, (\ref{desc. tensorial de cuspidales}) applies isomorphically to (%
\ref{desc. tensorial en cohomologia}) by the homomorphism 
\begin{equation*}
S_{k,k}^{oc}(N_{2},\mathbb{C}_{p})\longrightarrow H_{par-rig}^{2}(Y_{\mathbb{%
C}_{p}}^{\prime }(N_{2}),\mathcal{L}^{n}/\mathbb{C}_{p})(f)\text{ .}
\end{equation*}%
This works also if $f$ is only overconvergent, with $H_{par-dR}^{2}(Y_{%
\mathbb{C}_{p}}(N_{2}),\mathcal{L}^{n}/\mathbb{C}_{p})(f)$ instead of $%
H_{par-dR}^{2}(Y_{\mathbb{C}_{p}}^{\prime }(N_{2}),\mathcal{L}^{n}/\mathbb{C}%
_{p})(f)$.

Consider the degree four characteristic polynomial of $\phi $ on the
isotypic component of $\omega _{f}\in H^{0}(L_{Y(N_{2})}^{(k,k)}(-D^{c}))$ ,
with coefficients in $R$, thus in $\mathbb{Q}_{p}$, 
\begin{equation}
Q_{f}(x)=\dprod (1-\alpha _{i}^{-1}\alpha _{i^{\prime }}^{-1}x)=\dprod
(1-\alpha _{i}\alpha _{i^{\prime }}p^{2-2k}x)\in \mathbb{Q}_{p}[x]
\label{Q^f(t)}
\end{equation}%
generalizing the polynomial given in \ (\ref{polinomio minimo}) (Its
coefficients must be symmetric in $\alpha _{0}$ and $\alpha _{1}$ , and also
symmetric in $\alpha _{0^{\prime }}$ and $\alpha _{1^{\prime }}$, thus must
be rational expressions of $\alpha _{0}+\alpha _{1}=a_{\pi }$ and $\alpha
_{0}\alpha _{1}=p^{k-1}$and $\alpha _{0^{\prime }}+\alpha _{1^{\prime
}}=a_{\pi ^{\prime }}$and $\alpha _{0^{\prime }}\alpha _{1^{\prime
}}=p^{k-1} $, all values in $R$). As $Q_{f}(\phi )=0$ it divides the
characteristic polynomial of $\phi $ acting on $H_{dR}^{2+2n}(\widehat{A}%
_{Y(N_{2})}^{U,n})$ , so its roots $\alpha _{i}\alpha _{i^{\prime }}$ have
complex norm $p^{\frac{u}{2}}=p^{k-1}$ , for $u=2+2n$. As operative tool, we
will denote 
\begin{equation*}
Q_{ii^{\prime }}(x):=\frac{Q_{f}(x)}{(1-\alpha _{i}\alpha _{i^{\prime
}}p^{2-2k}x)}\text{ \ .}
\end{equation*}

.

Denote $\omega _{f^{\sharp }}$ the $L^{(k-2,k-2)}-$differential $\omega
_{f^{\sharp }}$ associated to the overconvergent $\mathfrak{a}$-Hilbert
cuspidal form

\begin{equation*}
f^{\sharp }:=Q_{f}(\phi )(f)=\dsum\limits^{{}}(-1)^{i+i^{\prime }}\frac{%
\alpha _{i}\alpha _{i^{\prime }}}{(\alpha _{0}-\alpha _{1})(\alpha
_{0^{\prime }}-\alpha _{1^{\prime }})}Q_{f}(\phi )f_{ii^{\prime }}\text{ }
\end{equation*}%
as follows from \ref{f como suma de cuatro} (where, abusing notation, we
write $\phi $ for the operator $p^{2k-2}V_{p}$ acting in $S_{k,k}^{oc}(N_{2},%
\mathbb{C}_{p})$). In the following proof of the main theorem \ref{Main
theorem copy(1)} we find, explicitely up to the equivalence $\approx $
defined in (\ref{Defino near- oc y equivalencia}), an $\mathcal{L}^{2n}$%
-differential $\varrho _{\mathcal{F^{\flat }}}$ with $\nabla ^{GM}(\varrho _{%
\mathcal{\mathcal{F^{\sharp }}}})=\omega _{f^{\sharp }}$.

\bigskip

\bigskip \textbf{Proof of theorem \ref{TEOREMA PRINCIPAL INTRODUCCION}: }To
compute 
\begin{align*}
Q_{f}(\phi )f& =\dsum\limits_{{}}^{4}(-1)^{i+i^{\prime }}\frac{\alpha
_{i}\alpha _{i^{\prime }}}{(\alpha _{0}-\alpha _{1})(\alpha _{0^{\prime
}}-\alpha _{1^{\prime }})}Q_{f}(\phi )f_{ii^{\prime }} \\
& =\dsum\limits_{{}}^{4}(-1)^{i+i^{\prime }}\frac{\alpha _{i}\alpha
_{i^{\prime }}}{(\alpha _{0}-\alpha _{1})(\alpha _{0^{\prime }}-\alpha
_{1^{\prime }})}Q_{ii^{\prime }}(\phi )(1-\alpha _{i}\alpha _{i^{\prime
}}V_{\pi }V_{\pi ^{\prime }})f_{ii^{\prime }}
\end{align*}%
decompose%
\begin{equation}
(1-\alpha _{i}\alpha _{i^{\prime }}V_{\pi }V_{\pi ^{\prime }})f_{ii^{\prime
}}=\frac{1}{2}(1-\alpha _{i}V_{\pi })(1+\alpha _{i^{\prime }}V_{\pi ^{\prime
}})f_{ii^{\prime }}+\frac{1}{2}(1-\alpha _{i^{\prime }}V_{\pi ^{\prime
}})(1+\alpha _{i}V_{\pi })f_{ii^{\prime }}\text{ .}  \label{Decompose}
\end{equation}%
The first summand equals%
\begin{equation}
\frac{1}{2}(1-\alpha _{i}V_{\pi })f_{ii^{\prime }}+\frac{1}{2}\alpha
_{i^{\prime }}V_{\pi ^{\prime }}(1-\alpha _{i}V_{\pi })f_{ii^{\prime }}=%
\frac{1}{2}(1-V_{\pi }U_{\pi })f_{ii^{\prime }}+\frac{1}{2}\alpha
_{i^{\prime }}V_{\pi ^{\prime }}(1-V_{\pi }U_{\pi })f_{ii^{\prime }}\text{ .}
\label{key step}
\end{equation}%
Observe that%
\begin{align*}
(1-V_{\pi }U_{\pi })f_{ii^{\prime }}& =(1-V_{\pi }U_{\pi })f_{i}-\alpha
_{j^{\prime }}(1-V_{\pi }U_{\pi })V_{\pi ^{\prime }}f_{i} \\
& =(1-V_{\pi }U_{\pi })f-\alpha _{j^{\prime }}(1-V_{\pi }U_{\pi })V_{\pi
^{\prime }}f\text{ ,}
\end{align*}%
because $(1-V_{\pi }U_{\pi })V_{\pi }=0$, as $U_{\pi }V_{\pi }=1$ .
Therefore, (\ref{key step}) equals 
\begin{equation}
\frac{1}{2}(1-V_{\pi }U_{\pi })f-\frac{\alpha _{j^{\prime }}}{2}(1-V_{\pi
}U_{\pi })V_{\pi ^{\prime }}f+\frac{\alpha _{i^{\prime }}}{2}(1-V_{\pi
}U_{\pi })V_{\pi ^{\prime }}f_{ii^{\prime }}  \label{All told}
\end{equation}%
Analogously, the second summand in \ref{Decompose} equals%
\begin{equation}
\frac{1}{2}(1-V_{\pi ^{\prime }}U_{\pi ^{\prime }})f-\frac{\alpha _{j}}{2}%
(1-V_{\pi ^{\prime }}U_{\pi ^{\prime }})V_{\pi }f+\frac{\alpha _{i}}{2}%
(1-V_{\pi ^{\prime }}U_{\pi ^{\prime }})V_{\pi }f_{ii^{\prime }}
\label{ultima}
\end{equation}

\bigskip

\textit{Lemma.} For any overconvergent Hilbert modular form $f$ of weight $%
(k,k)$, denote by $f^{\pi ^{\prime }[\pi ]}$ the $\pi -$depletion of $f^{\pi
^{\prime }}:=V_{\pi ^{\prime }}f$.

a) $\theta ^{-1}$ $f^{\pi ^{\prime }[\pi ]}\approx 0$ and $\theta ^{\prime
-1}$ $f^{\pi \lbrack \pi ^{\prime }]}\approx 0$

b) $\mathcal{F}^{\pi \lbrack \pi ^{\prime }]}\approx 0$, where $\nabla ^{GM}%
\mathcal{F}^{\pi \lbrack \pi ^{\prime }]}=\nabla ^{\prime }\mathcal{F}^{\pi
\lbrack \pi ^{\prime }]}=f^{\pi \lbrack \pi ^{\prime }]}$, i.e. $f^{\pi
\lbrack \pi ^{\prime }]}\approx 0$;\ analogously $\mathcal{F}^{\pi ^{\prime
}[\pi ]}\approx 0$, i.e. $f^{\pi ^{\prime }[\pi ]}\approx 0$

c) $f_{ii^{\prime }}^{[p]}\approx f^{[p]}$ $\approx \frac{1}{2}f^{[\pi ]}+%
\frac{1}{2}f^{[\pi ^{\prime }]}$

\bigskip\textit{Proof of the lemma}:

a) If $f(q)=\dsum\limits_{\nu \in (\mathfrak{a}^{-1})^{+}}^{{}}b_{\nu
}q^{\nu }$ , then%
\begin{equation*}
f^{\pi ^{\prime }[\pi ]}(q)=\dsum\limits_{\nu \in (\mathfrak{a}^{-1})^{+},%
\text{ }\pi \nmid \nu }^{{}}b_{\nu }q^{\pi ^{\prime }\nu }
\end{equation*}%
so that 
\begin{equation*}
\text{ }\theta ^{-1}f^{\pi ^{\prime }[\pi ]}(q)=\dsum\limits_{\nu \in (%
\mathfrak{a}^{-1})^{+},\text{ }\pi \nmid \nu }^{{}}(\pi ^{\prime }\nu
)^{-1}b_{\nu }q^{\pi ^{\prime }\nu }
\end{equation*}%
and 
\begin{equation*}
\theta ^{-1}f^{\pi ^{\prime }[\pi ]}\lfloor (q)=\dsum\limits_{\nu \in (%
\mathfrak{a}^{-1})^{+},\text{ }\pi \nmid \nu }^{{}}(\pi ^{\prime }\nu
)^{-1}b_{\nu }q^{Tr(\pi ^{\prime }\nu )}=\dsum\limits_{n\geq
1}^{{}}c_{n}q^{n}\text{ }
\end{equation*}%
with%
\begin{equation*}
\text{ }c_{n}=\dsum\limits_{\substack{ \nu \in (\mathfrak{a}^{-1})^{+},\text{
}\pi \nmid \nu  \\ Tr(\pi ^{\prime }\nu )=n}}^{{}}(\pi ^{\prime }\nu
)^{-1-t}b_{\nu }
\end{equation*}%
This is in the kernel of $U_{p}$ because $p\nmid n$ whenever $c_{n}\neq 0$.
Indeed, if $p$ divides $n=\pi ^{\prime }\nu +\pi x^{\prime }$, then also $%
\pi $ does, so it divides $\nu $, a contradiction. Therefore, it is also in
the kernel of $e_{ord}=\lim_{n}U_{p}^{n!}$.

\bigskip b) It is a consequence of a) and of the expression (\ref{expansion
of scriptF primitive}) of the $\nabla ^{\prime }-$primitive.

c) It is a consequence of b), and of (\ref{Decompose}), (\ref{All told}) and
(\ref{ultima}) $\square $

We use this lemma in combination with the fact that, for the overconvergent
primitive $(1,0)$-differential $\mathcal{\varrho }_{\mathcal{F}^{[\pi
^{\prime }]}}$ of $\omega _{f^{[\pi ^{\prime }]}}$ (i.e. with $\nabla ^{GM}%
\mathcal{\varrho }_{\mathcal{F}^{[\pi ^{\prime }]}}=\nabla ^{\prime }%
\mathcal{\varrho }_{\mathcal{F}^{[\pi ^{\prime }]}}=\omega _{f^{[\pi
^{\prime }]}}$ ), the first argument of the product proved in \ref{Clave 2
en Lausana} b) to compute the $p-$adic Abel Jacobi map is 
\begin{equation*}
e_{ord}j_{X^{\prime }(N)}^{\ast }sp(\mathcal{\varrho }_{\mathcal{F}^{[\pi
^{\prime }]}})=(-1)^{t}t!e_{ord}j_{X^{\prime }(N)}^{\ast }\varrho (\theta
^{\prime -1-t}f^{[\pi ^{\prime }]})
\end{equation*}%
(cf. last remark of the former section). Since $(1-\alpha _{i}\alpha
_{i^{\prime }}V_{p})f_{ii^{\prime }}=f_{ii^{\prime }}^{[p]}$, 
\begin{align*}
f^{^{\sharp }}& =Q_{f}(\phi )f=\dsum\limits_{i,i^{\prime
}}^{{}}(-1)^{i+i^{\prime }}\frac{\alpha _{i}\alpha _{i^{\prime }}}{(\alpha
_{0}-\alpha _{1})(\alpha _{0^{\prime }}-\alpha _{1^{\prime }})}Q_{ii^{\prime
}}(\phi )(1-\alpha _{i}\alpha _{i^{\prime }}V_{p})f_{ii^{\prime }} \\
& \approx \dsum\limits_{i,i^{\prime }}^{{}}(-1)^{i+i^{\prime }}\frac{\alpha
_{i}\alpha _{i^{\prime }}}{(\alpha _{0}-\alpha _{1})(\alpha _{0^{\prime
}}-\alpha _{1^{\prime }})}Q_{ii^{\prime }}(\phi )(\frac{f^{[\pi ]}}{2}+\frac{%
f^{[\pi ^{\prime }]}}{2})
\end{align*}%
so that $\omega _{f^{\sharp }}$ has a $\nabla ^{GM}-$primitive $\mathcal{%
\varrho }_{\mathcal{F^{\sharp }}}$ with the equivalence 
\begin{equation*}
\mathcal{\varrho }(\mathcal{F}^{\sharp })\approx
(-1)^{t}t!\dsum\limits_{i,i^{\prime }}^{{}}(-1)^{i+i^{\prime }}\frac{\alpha
_{i}\alpha _{i^{\prime }}}{(\alpha _{0}-\alpha _{1})(\alpha _{0^{\prime
}}-\alpha _{1^{\prime }})}Q_{ii^{\prime }}(\phi )\varrho (\frac{\theta
^{\prime -1-t}f^{[\pi ^{\prime }]}}{2}-\frac{\theta ^{-1-t}f^{[\pi ]}}{2})
\end{equation*}%
Therefore, the expression computing the $p-$adic Abel-Jacobi by \ref{AJ 1}
is 
\begin{align}
& (-1)^{t}t!\text{ }\dsum\limits_{i,i^{\prime }}(-1)^{i+i^{\prime }}\frac{%
\alpha _{i}\alpha _{i^{\prime }}}{(\alpha _{0}-\alpha _{1})(\alpha
_{0^{\prime }}-\alpha _{1^{\prime }})}  \label{scriptC} \\
& <(1-\alpha _{i}\alpha _{i^{\prime }}p^{2-2k}\phi
)^{-1}e_{g}e_{ord}j_{X^{\prime }(N)}^{\ast }\varrho (\frac{\theta ^{\prime
-1-t}f^{[\pi ^{\prime }]}}{2}-\frac{\theta ^{-1-t}f^{[\pi ]}}{2}),\eta
_{g}^{u-r}> \\
& =(-1)^{t}t!\frac{\mathcal{C}_{t}(f,g)}{(\alpha _{0}-\alpha _{1})(\alpha
_{0^{\prime }}-\alpha _{1^{\prime }})}\text{ }\dsum\limits_{i,i^{\prime
}}(-1)^{i+i^{\prime }}\frac{\alpha _{i}\alpha _{i^{\prime }}}{1-\alpha
_{i}\alpha _{i^{\prime }}\beta _{1}p^{-t}p^{2-2k}}
\end{align}%
where

\begin{equation}
\mathcal{C}_{t}(f,g):=<e_{g}e_{ord}\text{ }j_{X^{\prime}(N)}^{\ast}\varrho(%
\frac{\theta^{\prime-1-t}f^{[\pi^{\prime}]}}{2}-\frac{\theta ^{-1-t}f^{[\pi]}%
}{2}),\eta_{g}^{u-r}>  \label{true scriptC}
\end{equation}
and $\beta_{1}=\beta_{1}(g)$. Taking 
\begin{equation*}
\mathcal{E}(f,g):=\dprod
\limits_{i,i^{\prime}}(1-\alpha_{i}\alpha_{i}\beta_{1}p^{2-2k-t})
\end{equation*}
as common denominator in the summation, two numerators add up%
\begin{equation*}
\mathit{\alpha}_{0}\mathit{(\alpha}_{0^{\prime}}\mathit{-\alpha}_{1^{\prime}}%
\mathit{)(1-\alpha}_{1}\mathit{\alpha}_{0^{\prime}}\mathit{\beta p}^{2-2k+t}%
\mathit{)(1-\alpha}_{1}\mathit{\alpha}_{1^{\prime}}\mathit{\beta p}^{2-2k+t}%
\mathit{)}
\end{equation*}
and the other two add up%
\begin{equation*}
-\alpha_{1}(\alpha_{0^{\prime}}-\alpha_{1^{\prime}})(1-\alpha_{0}\alpha_{0^{%
\prime}}\beta p^{2-2k+t})(1-\alpha_{0}\alpha_{1^{\prime}}\beta p^{2-2k+t})%
\text{,}
\end{equation*}
so all four numerators add up

\begin{align*}
& \mathit{(\alpha }_{0}\mathit{-\alpha }_{1}\mathit{)(\alpha }_{0^{\prime }}%
\mathit{-\alpha }_{1^{\prime }}\mathit{)(1-\alpha }_{0}\mathit{\alpha }_{1}%
\mathit{\alpha }_{0^{\prime }}\mathit{\alpha }_{1^{\prime }}\mathit{\beta }%
^{2}\mathit{p}^{4-4k+2t}\mathit{)} \\
& =\mathit{(\alpha }_{0^{\prime }}\mathit{-\alpha }_{1^{\prime }}\mathit{%
)(\alpha }_{0}\mathit{-\alpha }_{1}\mathit{)(1-\beta }^{2}\mathit{p}%
^{2-2k+2t}\mathit{)}
\end{align*}%
since $\mathit{\alpha }_{0}\mathit{\alpha }_{1}\mathit{\alpha }_{0^{\prime }}%
\mathit{\alpha }_{1^{\prime }}=\mathit{p}^{2k-2}$ (and the cancellation of
the missing terms has also used $\alpha _{0}\alpha _{1}=\alpha _{0^{\prime
}}\alpha _{1^{\prime }}=p^{k-1}$). Therefore the factor multiplying to $%
\mathcal{C}_{p}(f,g)$ in (\ref{scriptC}) is $\mathcal{E}_{1}(g)/\mathcal{E}%
(f,g)$ for 
\begin{equation*}
\mathcal{E}_{1}(g):=\mathit{(1-\beta }^{2}\mathit{p}^{2-2k+2t}\mathit{%
)=(1-\beta }^{2}\mathit{p}^{-k_{0}}\mathit{)}
\end{equation*}

\bigskip We finally check that $\mathcal{E}_{1}(g)$ does not vanish. Observe
that it vanishes if and only if $\beta _{1}=\pm p^{k-1}$\ . Recall that $%
\beta _{0}$\ is a $p$- adic unit,\ so $\beta _{1}$\ has valuation $p^{k-1}$\
as $\beta _{0}\beta _{1}=p^{k-1}$ . In case the Euler factor vanished, it
would be $b_{p}=\beta _{0}+\beta _{1}=\pm (1+p^{k-1})$, contradicting the
fact that $\left\vert b_{p}\right\vert \leq 2p^{\frac{k-1}{2}}$according to
the Ramanujan-Petersson conjecture for Hecke newforms of arbitrary weight,
proved by Deligne in \cite{Deligne on Ramanujan}.\textit{\ } $\square $

\bigskip

As pointed out in the introduction, the right hand side of the equality of
theorem 3.12 of \cite{DR} has the opposite sign to ours. This is becasue
each of the four summands $g_{\alpha }\otimes h_{\beta }$ in which the
authors split the product $gh$, are transformed, after Frobenius and
restriction, to the sum of $-g_{\alpha }^{[p]}h_{\beta }$ plus the
derivative of a $p$-depleted power series, which lies in the kernel of $%
e_{ord}.$

\bigskip

\section{\textbf{Hida families and a }$p-$\textbf{adic Gross-Zagier formula}}

\bigskip

Decompose the $p-$profinite group 
\begin{equation*}
Z(N_{2}):=\lim Cl_{\mathfrak{K}}(N_{2}p^{\alpha })\approx W\times
Z(N_{2})_{tor}
\end{equation*}%
in torsion free part $W$ and torsion part $Z(N_{2})_{tor}$. The action of $%
Z(N_{2}):=\lim Cl_{\mathfrak{K}}(N_{2}p^{\alpha })$ on the $p^{\alpha }$
-roots of unity yields a homomorphism $Z(N_{2})\longrightarrow Gal(\mathfrak{%
K}(\mu _{p^{\alpha }}),\mathfrak{K})$ , and taking limit we obtain the
cyclotomic character 
\begin{equation*}
\chi :Z(N_{2})\longrightarrow Gal(\mathfrak{K}(\mu _{p^{\infty }}),\mathfrak{%
K})\longrightarrow Gal(\mathbb{Q}(\mu _{p^{\infty }}),\mathbb{Q})\approx
Aut(\mu _{p^{\infty }})=\mathbb{Z}_{p}^{\ast }
\end{equation*}%
This restricts to isomorphisms between their torsion free parts%
\begin{equation}
\chi \lfloor :W\overset{\approx }{\longrightarrow }Gal(\mathfrak{K}(\mu
_{p^{\infty }}),\mathfrak{K(}\mu _{p}))\overset{\approx }{\longrightarrow }%
\Gamma _{p}:=1+p\mathbb{Z}_{p}\text{ ,}  \label{isomorfismo dando alfa=1}
\end{equation}%
the one in the middle being justified because

\begin{equation}
1\longrightarrow Gal(\mathfrak{K}(\mu _{p^{\infty }}),\mathfrak{K(}\mu
_{p}))\longrightarrow Gal(\mathfrak{K}(\mu _{p^{\infty }}),\mathfrak{K}%
)\longrightarrow Gal(\mathfrak{K}(\mu _{p}),\mathfrak{K})\longrightarrow 1
\label{sucesion corta dando alfa=1}
\end{equation}%
is the projection to the torsion part. The fact that $\chi \lfloor $ is an
isomorphism is a particular case of the\ Leopoldt conjecture, which holds
for totally real fields, in particular for $\mathfrak{K}$. The classical
characters $Z(N_{2})\longrightarrow \mathbb{Z}_{p}^{\ast }\subseteq \mathbb{C%
}_{p}^{\ast }$ are the integer powers $\chi ^{l}$ of this cyclotomic
character.

The Iwasawa $\mathbb{Z}_{p}-$algebra $\Lambda \mathbb{\ }$is $\mathbb{Z}%
_{p}[[W]]$ . Inside the group of characters%
\begin{equation*}
\mathcal{\mathfrak{X}(\Lambda )=}Hom_{\mathbb{Z}_{p}}(\Lambda ,\overline{%
\mathbb{Q}}_{p})\approx Hom_{cts}(\Gamma _{p},\overline{\mathbb{Q}}%
_{p}^{\ast })\text{ ,}
\end{equation*}%
denote $\mathcal{\mathfrak{X}}_{\text{alg}}\mathcal{(\Lambda )\subseteq 
\mathfrak{X}(\Lambda )}$ the dense subset of characters of the form $\chi
^{l}\varepsilon $ for a torsion character $\varepsilon :\Gamma
_{p}\longrightarrow \overline{\mathbb{Q}}_{p}^{\ast }$ , so the classical
characters are those in 
\begin{equation*}
\mathcal{\mathfrak{X}}_{\text{clas}}\mathcal{(\Lambda )}=\mathcal{\mathfrak{X%
}}_{\text{alg}}\mathcal{(\Lambda )\cap \mathcal{\mathfrak{X}(\Lambda )}}_{tf}
\end{equation*}%
i.e. algebraic with torsion factor $\varepsilon =1$. This is dense in the
open torsion free part $\mathcal{\mathfrak{X}(\Lambda )}_{tf}\subseteq 
\mathcal{\mathfrak{X}(\Lambda )}$.

We recall next the following definitions from \cite{DR} and \cite{Bertoli y
Darmon}.

\textbf{Definition}.\label{HIda} A Hida family of tame level $N>1$ , i.e.
not divisible by $p$ , is $\mathbf{g}=(\Lambda _{\mathbf{g}},\Omega _{%
\mathbf{g}},\Omega _{\mathbf{g}}^{cl},\mathbf{g}(q))$ , where

$\bullet $ $\Lambda _{\mathbf{g}}$ is a finite flat extension of $\Lambda $,
thus inducing $\zeta _{\mathbf{g}}:Hom(\Lambda _{\mathbf{g}},\mathbb{C}%
_{p})\longrightarrow \mathcal{\mathfrak{X}(\Lambda )}$

$\bullet $ $\Omega _{\mathbf{g}}$ is a non-empty open subset of $Hom(\Lambda
_{\mathbf{g}},\mathbb{C}_{p})$.

$\bullet $ $\Omega _{\mathbf{g}}^{cl}$ is a $p-$adically dense subset of $%
\Omega _{\mathbf{g}}$ whose image under $\zeta _{\mathbf{g}}$ lies in $%
\mathcal{\mathfrak{X}}_{\text{clas}}\mathcal{(\Lambda )}$.

$\bullet $ $\mathbf{g}(q)=\dsum\limits_{n\geq 1}^{{}}\mathbf{a}_{n}q^{n}$ $%
\in \Lambda _{\mathbf{g}}[[q]]$ \ is such that, for all $r\in \Omega _{%
\mathbf{g}}^{cl}$ , the power series $g_{r}^{(p)}(q):=\dsum\limits_{n\geq
1}^{{}}\mathbf{a}_{n}(r)q^{n}$ is the $q-$expansion of the ordinary $p-$%
stabilization $g_{r}^{(p)}$ of a normalised newform $g_{r}$ of weight $\zeta
(r)$ for $\Gamma _{0}(N)$ .

\textbf{Definition}. \label{Landa-adic}A $\Lambda -$adic modular form of
tame level $N$ is a quadruple $(\Lambda _{\phi },\Omega _{\phi },\Omega
_{\phi }^{cl},\phi (q))$ where

$\bullet \Lambda _{\phi }$ is a complete finitely generated (but not
necessarily finite), flat extension of $\Lambda $, thus inducing $\zeta
_{\phi }:Hom(\Lambda _{\phi },\mathbb{C}_{p})\longrightarrow \mathcal{%
\mathfrak{X}(\Lambda )}$

$\bullet $ $\Omega _{\phi }$ is a non-empty open subset of $Hom(\Lambda
_{\phi },\mathbb{C}_{p})$

$\bullet $ $\Omega _{\phi }^{cl}$ is a $p-$adically dense subset of $\Omega
_{\phi }$ whose image under $\zeta _{\phi }$ lies in $\mathcal{\mathfrak{X}}%
_{\text{clas}}\mathcal{(\Lambda )}$.

$\bullet $ $\phi (q)=\dsum\limits_{n\geq 1}^{{}}\mathbf{a}_{n}q^{n}\in
\Lambda _{\phi }[[q]]$ is such that, for all $r\in \Omega _{\phi }^{cl}$ ,
the power series $\phi _{r}^{(p)}(q):=$ $\dsum\limits_{n\geq 1}^{{}}\mathbf{a%
}_{n}(r)q^{n}$ is the $q-$expansion of a classical ordinary cusp form in $%
S_{\zeta _{\phi }(r)}(\Gamma _{1}(N)\cap \Gamma _{0}(p),\overline{\mathbb{Q}}%
)\otimes \mathbb{C}_{p}$ .

\bigskip\ Next, we present the following material on $p-$adic families of
non-ordinary $\mathfrak{a-}$Hilbert cuspidal forms, as in \cite{Yamagami}
sections 3.2 and 4. First, consider the increasing function $\mathfrak{%
\vartheta }:\mathbb{Q}_{\geq 0}\longrightarrow \mathbb{Q}_{\geq 0}$ \ given
by

\begin{equation*}
\mathfrak{\vartheta }(x):=[2xh^{+}\left( \frac{18(p-1)}{2(p-2)}x+2\right) ]%
\text{\ }
\end{equation*}%
where \ $h^{+}$ is the narrow class number of $\mathfrak{K}$ and the bracket
is the integer part. For a positive rational $\sigma $ denote 
\begin{equation*}
\Lambda _{\sigma }=R_{\mathfrak{m}}<\frac{X}{p^{\mathfrak{\vartheta }(\sigma
)}}>:=\{\dsum\limits_{n\geq 0}^{{}}a_{n}\text{ \ \ }\ (\frac{X}{p^{\mathfrak{%
\vartheta }(\sigma )}})^{n}\text{\ \ }\mid \lim \left\vert a_{n}\right\vert
_{p}=0\text{\ \ }\}\approx \mathbb{Z}_{p}[[X]]_{p^{\mathfrak{\vartheta }%
(\sigma )}}
\end{equation*}%
and for $z\in \mathcal{O}_{\mathbb{C}_{p}}$, consider the $p-$adic ball%
\begin{equation*}
B_{z,\mathfrak{\vartheta }(\sigma )}=\{w\in \mathcal{O}_{\mathbb{C}_{p}}\mid
\left\vert w-z\right\vert _{p}\leq p^{-\mathfrak{\vartheta }(\sigma )}\}%
\text{.}
\end{equation*}%
Any element of $\Lambda _{\sigma }$ defines a rigid analytic function on $%
B_{z,\mathfrak{\vartheta }(\sigma )}$ defined over $K_{\mathfrak{m}}$ .

\textbf{Definition }(Yamagami) \label{Definicion familia yamagami} A $p-$%
adic family of $\mathfrak{a-}$Hilbert cuspidal forms of tame level $%
N_{2}\geq 4$ , i.e. not divisible by $p$ , and positive rationals $\sigma
\geq \sigma ^{\prime }$ as slope, is $\mathbf{f}$ $=(\Lambda _{\mathbf{f}},$%
\ $\Omega _{\mathbf{f}},\Omega _{\mathbf{f}}^{cl},\mathbf{f}(q))$ with

$\bullet $ $\Lambda _{\mathbf{f}}=\Lambda _{\sigma }$

$\bullet $ $\Omega _{\mathbf{f}}=B_{n_{0},\mathfrak{\vartheta }(\sigma )}$
for some $n_{0}\in \mathbb{N}$

$\bullet $ $\Omega _{\mathbf{f}}^{cl}=\{n\in B_{n_{0},\mathfrak{\vartheta }%
(\sigma )}\cap \mathbb{\mathbb{N}}\mid n\equiv n_{0}$ \ $(\func{mod}2)\}$

$\bullet $ $\mathbf{f}(q)=$\ $\dsum\limits_{\nu \in \mathfrak{a}%
_{+}^{-1}}^{{}}\mathbf{a}_{\nu }q^{\nu }\in \Lambda _{\mathbf{f}}[[q^{\nu }]]
$\ \ specializing to each $n\in \Omega _{\mathbf{f}}^{cl}$ with $n+1>\sigma $
to the $q-$expansion $f_{n}(q)=\dsum\limits_{\nu \in \mathfrak{a}%
_{+}^{-1}}^{{}}\mathbf{a}_{\nu }(n)q^{\nu }$\ of an eigenform $f_{n}$ $\in
S_{n+2,n+2}(\mathfrak{a,}N,\overline{\mathbb{Q}})^{\sigma ,\sigma ^{\prime }}
$of the full Hecke algebra.

The assumption in the following theorem is not essentially restrictive, as
it dissapears while working adelicaly and imposing on $f$ the condition of
being $\pi $ and $\pi ^{\prime }$-new, which is a mild condition playing in
the finite slope setting the role of the primitivity condition in the
ordinary setting (for details cf. Theorem 4.5 in \cite{Yamagami}).

\textbf{Theorem. }\ \ Suppose $\sigma $ and $\sigma ^{\prime }$ are chosen
so that 
\begin{equation}
\dim _{\mathbb{C}_{p}}S_{k,k}(\mathfrak{a,}N_{2},\mathbb{C}_{p})^{\sigma
,\ast }=\dim _{\mathbb{C}_{p}}S_{k,k}(\mathfrak{a,}N,\overline{\mathbb{Q}}%
)^{\ast ,\sigma ^{\prime }}=1  \label{monodim. condition}
\end{equation}%
and let $f\in S_{k,k}(\mathfrak{a,}N_{2},\overline{\mathbb{Q}})^{\sigma
,\sigma ^{\prime }}$ be a Hecke-eigenform (i.e. for the full Hecke algebra
of operators). There exists a $p-$adic family $\mathbf{f}$ $=(\Lambda _{%
\mathbf{f}},$\ $\Omega _{\mathbf{f}}$,$\Omega _{\mathbf{f}}^{cl},\mathbf{f}%
(q))$ of $\mathfrak{a-}$Hilbert cuspidal forms of level $N_{2}$ and slope $%
(\sigma ,\sigma ^{\prime })$ such that $f_{k-2}=f$.\footnote{%
In \cite{Yamagami} theor. 4.4, $p$-adic families of Hilbert modular forms
are presented as collections of rigid analytic functions $\{a_{T}\}$ on $p$%
-adic balls of the form $B_{n_{0},\vartheta (\sigma )}$, indexed by all the
Hecke operators, such that their evaluation at classical weights are Hecke
eigenvalues of (adelic) Hilbert modular forms.
\par
We restrict ourselves to $\mathfrak{a}$-Hilbert modular forms, embedded into
adelic ones as $f\mapsto (0,...,0,f,0,...,0)$ (in the $\mathfrak{a}$%
-argument). Yamagami%
\'{}%
s $p$-adic families applied to our particular case are equivalent to ours.
Indeed, it is enough to keep only the sub-family of rigid analytic functions 
$\{a_{\nu }\}$ with $\nu \in (\mathfrak{a}^{-1})^{+}$, (since the other
operators act trivially on our $\mathfrak{a}$-Hilbert modular forms embedded
into the direct sum, as we explain section 3) and to consider the $\mathfrak{%
a}$-component of Yamagami%
\'{}%
s specialisations. The functions $\{a_{\nu }\}$ belong to the Tate algebra
which we have taken as first input in the quadruple defining $\mathbf{f}$ \
(def. \ref{Definicion familia yamagami}).
\par
Finally, we pack together these functions into a formal power series to be
in accordance with the presentation in \cite{DR} of Hida families and $%
\Lambda $-adic modular forms (recalled in def. \ref{HIda} and \ref%
{Landa-adic}). With this we produce, at the end of our construction, a $%
\Lambda $-adic modular form which is taken as an input together with the
Hida family of modular cusp forms passing through $g$, to apply \cite{DR}
prop. 2.19 and obtain our $p$-adic $L$-function.}

Given a rational cuspidal Hecke eigenform $g$, normalized and ordinary at $p$%
, new of weight\textit{\ }$\mathit{2}$ and of $\Gamma _{0}-$ level $N=d_{%
\mathfrak{K}}N_{1}N_{2}$ the one-dimensional Hida theorem (quoted in \cite%
{DR} 2.15, cf. also \cite{Darmon y Bertolini}(12)) associates a Hida family $%
\mathbf{g=(}\Lambda _{\mathbf{g}},\Omega _{\mathbf{g}},\Omega _{\mathbf{g}%
}^{cl},\mathbf{g}(q))$ of $\Gamma _{0}-$ level $N$ such that, for some
classical weight $r$ with $\zeta _{\mathbf{g}}(r)=2$, it is $g_{r}=g$ .

Consider now the $\Lambda -$adic family $\mathbf{h=(}\Lambda _{\mathbf{h}%
},\Omega _{\mathbf{h}},\Omega _{\mathbf{h}}^{cl},\mathbf{h}(q))$ of cuspidal
forms of $\Gamma _{0}-$ level $N$ , given in the following way:

$\bullet $ $\Lambda _{\mathbf{h}}$ is the flat $\Lambda -$algebra 
\begin{equation*}
\Lambda \longrightarrow \Lambda \otimes _{\mathbb{Z}_{p}}\Lambda
\longrightarrow \Lambda \otimes \Lambda _{\sigma }=\Lambda _{\mathbf{h}}
\end{equation*}%
the left homomorphism being given by $[z]\mapsto \lbrack z^{2}]\otimes
\lbrack z^{2}]$

$\bullet $ $\Omega _{\mathbf{h}}\subseteq Hom_{\mathbb{Z}_{p}}(\Lambda _{%
\mathbf{h}},\mathbb{Q}_{p})$ is the $p-$adic open 
\begin{eqnarray*}
Hom_{cts}(\Gamma _{p},\overline{\mathbb{Q}}^{\ast })\times D_{0,\mu (\sigma
)} &\subseteq &Hom_{cts}(\Gamma _{p},\overline{\mathbb{Q}}^{\ast })\times
B_{0,\mathfrak{\vartheta }(\sigma )} \\
&\subseteq &Hom_{\mathbb{Z}_{p}}(\Lambda \otimes _{\mathbb{Z}_{p}}\Lambda
_{\sigma },\overline{\mathbb{Q}}_{p})\text{ ,}
\end{eqnarray*}%
where $D_{0,\mathfrak{\vartheta }(\sigma )}=\{z\in \mathbb{Q}_{p}\mid
\left\vert z\right\vert _{p}\leq p^{-\mathfrak{\vartheta }(\sigma )}\}$ ,
and the last inclusion is given by multiplication in $\overline{\mathbb{Q}}%
_{p}$ \ .

\bigskip $\bullet $ To define $\Omega _{\mathbf{h}}^{cl}$, recall that any $%
z\in \mathbb{Z}_{p}^{\ast }$ is decomposable as 
\begin{equation*}
z=\mu (z)<z>\text{,}
\end{equation*}%
with $<z>\in \Gamma _{p}$ and $\mu (z)$ the only (after Hensel lemma) $p-1$
root of $1$ with $\mu (z)$ congruent to $z$ \ mod. $p$ .\ Our dense subset $%
\Omega _{\mathbf{h}}^{cl}\subseteq \Omega _{\mathbf{h}}$ will be the subset
of\textbf{\ } $(\overline{-1})$ $\times \Omega _{\mathbf{f}}^{cl}$
consisting of pairs $(j,s)$ , with $\zeta _{\mathbf{f}}$ $(s)$ classical
lying over $k$, such that $k+j\geq 1$ (Here we identify an integer $k$\ with
the continuous homomorphism $\Gamma _{p}\longrightarrow \overline{\mathbb{Q}}%
_{p}^{\ast }$\ rising to the $k-$th power, and $\overline{-1}$ denotes the
class of $-1$ $\func{mod}.$ $p-1$ ). The image by the epimorphism $\zeta _{%
\mathbf{h}}$ lies indeed in the set $\mathbb{Z}^{\geq 2}$\ of classical
characters $\Gamma _{p}\longrightarrow \overline{\mathbb{Q}}_{p}^{\ast }$\ .
It is a consequence of the fact that, for $(j,s)\in \Omega _{\mathbf{h}%
}^{cl} $ and $\zeta _{\mathbf{f}}(s)=\chi ^{k}$\ , it is 
\begin{equation*}
\zeta _{\mathbf{h}}(j,s)(z)=j(z^{2})s(z^{2})=z^{2j+2k}\text{ }
\end{equation*}

$\bullet$ \ First, recall that 
\begin{equation*}
f_{s}^{(\pi,\pi^{\prime})[\pi]}+f_{s}^{(\pi,\pi^{\prime})[\pi^{\prime}]}%
\approx f_{s}^{[\pi]}\ +f_{s}^{[\pi^{\prime}]}
\end{equation*}

\bigskip by the lemma inside the proof of \ \ref{TEOREMA PRINCIPAL
INTRODUCCION}, so that 
\begin{align}
& e_{ord}(\frac{\theta ^{\prime -1-t}f_{s}^{[\pi ^{\prime }]}\ -\theta
^{-1-t}f_{s}^{[\pi ]}}{2})\lfloor (q)  \label{expansion de fpis+fpi's} \\
& =\frac{1}{2}e_{ord}(\dsum\limits_{\nu \in (\mathfrak{a}^{-1})^{+}\text{, }%
\pi ^{\prime }\nmid \nu }^{{}}\nu ^{\prime -1-t}\mathbf{a}_{\nu
}(s)q^{Tr(\nu )\text{ \ }}-\dsum\limits_{\nu \in (\mathfrak{a}^{-1})^{+}%
\text{, }\pi \nmid \nu }^{{}}\nu ^{-1-t}\mathbf{a}_{\nu }(s)q^{Tr(\nu )\text{
\ }})\text{ .}
\end{align}

With this in mind, we define $\mathbf{h}(q)$ as 
\begin{align}
\frac{1}{2}e_{ord}(\dsum \limits_{\nu\in(\mathfrak{a}^{-1})^{+}\text{, }%
\pi^{\prime }\nmid\nu\text{ i.e.}\pi\nmid\nu^{\prime}\text{ }%
}^{{}}\mu(\nu^{\prime})^{-1}[ & <\nu^{\prime}>]\otimes\mathbf{a}_{\nu}\text{ 
}q^{Tr(\nu)})  \notag \\
-\frac{1}{2}e_{ord}(\dsum \limits_{\nu\in(\mathfrak{a}^{-1})^{+}\text{, }%
\pi\nmid\nu}^{{}}\mu(\nu)^{-1}[ & <\nu>]\otimes\mathbf{a}_{\nu}\text{ }%
q^{Tr(\nu)})
\end{align}
It has, indeed, all coefficients in $\Lambda_{\mathbf{h}}$ , since $\pi
\nmid\nu,\nu^{\prime}$ implies $\pi\nmid\mu(\nu),\mu(\nu^{\prime}) $ so that 
$\mu(\nu)^{j},\mu(\nu^{\prime})^{j}\in R_{\mathfrak{m}}\approx\mathbb{Z}_{p}$
also for negative values of $j$ (as $\mathfrak{m}$ lies over $\pi $)\textbf{%
\ The key point is} that, for $(j,s)\in\Omega_{\mathbf{h}}^{cl}$ with $j=-1$ 
$+(p-1)m\geq-1$ (so that $j$ is $-1$ in $(\mathbb{Z}/p)^{\ast}$) \textbf{and 
}$\zeta_{\mathbf{f}}(s)=\chi^{k}$, it is%
\begin{align*}
(\mu(\nu)^{-1}[ & <\nu>]\otimes\mathbf{a}_{\nu})(j,s)=\mu(\nu)^{-1}<\nu >^{j}%
\mathbf{a}_{\nu}(s) \\
& =\mu(\nu)^{j}<\nu>^{j}\mathbf{a}_{\nu}(s)=\nu^{j}a_{\nu}(s)
\end{align*}
and analogously $(\mu(\nu^{\prime})^{-1}[<\nu^{\prime}>]\otimes\mathbf{a}%
_{\nu})=\nu^{\prime j}\mathbf{a}_{\nu}$, so that

\begin{equation*}
h_{(j,s)}^{(p)}(q)=-\frac{1}{2}e_{ord}(\dsum \limits_{\nu\in(\mathfrak{a}%
^{-1})^{+}\text{, }\pi\nmid\nu}^{{}}\nu^{j}\mathbf{a}_{\nu}(s)q^{Tr(\nu)})+%
\frac{1}{2}e_{ord}(\dsum \limits_{\nu\in(\mathfrak{a}^{-1})^{+}\text{, }%
\pi\nmid \nu^{\prime}}^{{}}\nu^{\prime j}\mathbf{a}_{\nu}(s)q^{Tr(\nu)}
\end{equation*}
The overconvergent cuspidal form $h_{(j,s)}^{(p)}$ has $\zeta_{\mathbf{h}%
}(j,s)=2n+2j\geq2$ and, being overconvergent of level $N$ and ordinary, it
is classical of level $Np$: 
\begin{equation*}
h_{(j,s)}^{(p)}\in S_{2j+2n}(\Gamma_{0}(Np),\overline{\mathbb{Q}})\subseteq
S_{2j+2n}(\Gamma_{1}(N)\cap\Gamma_{0}(p),\overline{\mathbb{Q}})\text{ .}
\end{equation*}

Denote $\Omega_{\mathbf{f,g}}$ and $\Omega_{\mathbf{f,g}}^{cl}$ the two
fibre products%
\begin{equation}
\begin{array}{ccc}
\Omega_{\mathbf{f,g}} & \longrightarrow & \Omega_{\mathbf{g}} \\ 
\downarrow &  & \downarrow\zeta_{\mathbf{g}} \\ 
\text{\ }\Omega_{\mathbf{h}} & \overset{\zeta_{\mathbf{h}}}{\longrightarrow}
& \mathcal{\mathfrak{X}(\Lambda)}%
\end{array}
\text{ and }%
\begin{array}{ccc}
\Omega_{\mathbf{f,g}}^{cl} & \longrightarrow & \Omega_{\mathbf{g}}^{cl} \\ 
\downarrow &  & \downarrow\zeta_{\mathbf{g}} \\ 
\text{\ }\Omega_{\mathbf{h}}^{cl} & \overset{\zeta_{\mathbf{h}}}{%
\longrightarrow} & \mathcal{\mathfrak{X}}_{\text{clas}}\mathcal{(\Lambda)}%
\end{array}
\text{ }  \label{fibre product}
\end{equation}

For any pair of such families $\mathbf{f}$ and $\mathbf{g}$, and classical
weights $r\in\Omega_{\mathbf{g}}^{cl}$\ lying over $k_{0}$ and $s\in \Omega_{%
\mathbf{f}}^{cl}$\ lying over $k$, with 
\begin{equation*}
k_{0}=\zeta_{\mathbf{g}}(r)=\zeta_{\mathbf{h}}(j,s)=2k+2j\geq2
\end{equation*}
i.e. with $(r,s)\in\Omega_{\mathbf{f,g}}^{cl}$ in \ref{fibre product}, the
ordinary cuspidal form $e_{g_{r}}h_{r,s}^{(p)}$ of level $Np$ is, by the
comment preceding lemma 2.19 in \cite{DR}, the ordinary stabilization of an
ordinary cuspidal form $h_{r,s}$\ \ of level $N$, i.e.

\begin{equation}
e_{g_{r}}h_{r,s}^{(p)}=h_{r,s}-\beta_{1}(g_{r})V_{p}h_{r,s}
\label{defino hrs}
\end{equation}

\bigskip Lemma 2.19 of \cite{DR} tells us that there is a partial \
"Garrett-Hida $p-$adic $L-$function", i.e.a function $L_{p}(\mathbf{f},%
\mathbf{g)}:\Omega_{\mathbf{f,g}}\dashrightarrow\mathbb{C}_{p}$ having poles
at finitely many nonclassical weights and satisfying at classical weights $%
(r,s)\in\Omega_{\mathbf{f,g}}^{cl}$,

\qquad%
\begin{align}
L_{p}(\mathbf{f},\mathbf{g)}(r,s) & \mathbf{=}\frac{%
[e_{g_{r}}h_{r,s}^{(p)},g_{r}^{(p)}]_{Np}}{[g_{r}^{(p)},g_{r}^{(p)}]_{Np}}=
\label{PW product} \\
\frac{\lbrack h_{r,s},g_{r}]_{N}}{[g_{r},g_{r}]_{N}} &
=<\omega_{h_{r,s}},\eta_{g_{r}}^{u-r}>\in\mathbb{C}_{p}\text{ .}
\end{align}
The brackets here indicate the Petersson inner product of cuspidal forms of
weight $k_{0}$ and indicated level, and $<\omega_{h_{r,s}},%
\eta_{g_{r}}^{u-r}>$ denotes, as always, the Poincar\'{e} pairing in $%
H_{dR}^{1}(X(N),\mathcal{L}_{0}^{k_{0}-2})$ of the classes of the
overconvergent $\mathcal{L}_{0}^{k_{0}-2}$-differentials $\omega_{h_{r,s}}$
and $\eta_{g_{r}}^{u-r}$ . this last defined in \ref{unit-root in higher
weights} by $<\omega_{g_{r}},\eta_{g_{r}}^{u-r}>=1$ (This class is presented
in \cite{DR} as lift of the class in $H_{dR}^{1}(X(N),L_{0}^{-(k_{0}-2)})$
represented by the $L_{0}^{-(k_{0}-2)}$-antiholomorphic form hermitian-dual
to the $L_{0}^{k_{0}-2}$- holomorphic form $\omega_{g_{r}}$ , and last
equality is just \cite{DR} \ (111), formula \ref{PW product}).

We are now in position to prove the second theorem stated at the
introduction.

\textbf{Proof of theorem \ref{Main theorem copy(1)}.} The proof that $%
\mathcal{E}_{0}(g_{r})$ is non-null is analogous to the proof of the
analogous statement for $\mathcal{E}_{1}(g_{r})$. The last assert belongs to
the theorem \ref{TEOREMA PRINCIPAL INTRODUCCION} already proved.

Denote $\alpha_{1}(g_{r})$ and $\beta_{1}(g_{r})$ as $\alpha_{1}$ and $%
\beta_{1}$. Applying the Hida projector $e_{ord}$ to \ref{defino hrs}, we
obtain%
\begin{equation*}
e_{g_{r}}h_{r,s}^{(p)}=e_{ord}h_{r,s}-e_{ord}\beta_{1}V_{p}h_{r,s}\text{ .}
\end{equation*}

\bigskip This represents the same cohomology class as 
\begin{align*}
e_{ord}h_{r,s}-\beta_{1}V_{p}e_{ord}h_{r,s} &
=(1-\beta_{1}\beta_{0}^{-1})e_{ord}h_{r,s} \\
& =(1-\beta_{1}^{2}p^{1-k_{0}})e_{ord}h_{r,s} \\
& =\mathcal{E}_{0}\mathcal{(}g_{r})e_{ord}h_{r,s}\text{ , }
\end{align*}
\textbf{\ }

since $V_{p}$ commutes with $U_{p}$ , thus with $e_{ord}$ , in cohomology.
Now, by \ref{PW product}, 
\begin{equation}
\text{ }L_{p}(\mathbf{f,g)(}r,s)=<\omega _{e_{ord}h_{r,s}},\eta
_{g_{r}}^{u-r}>=\mathcal{E}_{0}(g_{r})^{-1}<\omega
_{e_{g_{r}}h_{r,s}^{(p)}},\eta _{g_{r}}^{u-r}>\in \mathbb{C}_{p}\text{ .}
\label{calculo L familiar}
\end{equation}

By theorem \ref{TEOREMA PRINCIPAL INTRODUCCION},

\begin{align*}
& AJ_{p}(\Delta _{k-2,k_{0}-2})(\omega _{f}\otimes \eta _{g}^{u-r}) \\
& \mathcal{=}(-1)^{t}t!\frac{\mathcal{E}_{1}(g_{r})}{\mathcal{E}(f_{s},g_{r})%
}<e_{g}e_{ord}\text{ }j_{X^{\prime }(N)}\frac{1}{2}(\varrho (\theta ^{\prime
-1-t}f_{s}^{[\pi ^{\prime }]}-\theta ^{-1-t}f_{s}^{[\pi ]}),\eta
_{g_{r}}^{u-r}> \\
& =(-1)^{t}t!\frac{\mathcal{E}_{1}(g_{r})}{\mathcal{E}(f_{s},g_{r})}<\omega
_{e_{g_{r}}h_{r,s}^{(p)}},\eta _{g_{r}}^{u-r}>\text{,}
\end{align*}%
.

so that \ref{calculo L familiar} proves \ref{formula teorema cap 6} .$%
\square $

\bigskip

We end by justifying our denomination of $L_{p}(\mathbf{f,\mathbf{g)}}$ as a 
$p-$adic $L-$function. First, we recall that at classical weights $r$ and $s$
above $(k,k)$ and $k_{0}$ with $k-(k_{0}/2)\geq2,$ the value $L_{p}(\mathbf{f%
}$,$\mathbf{g})(r,s)$ computes the $p$-adic Abel-Jacobi map at the cycle $%
\Delta_{n,n_{0}}$ evaluated at suitable differential forms. The case $%
k-(k_{0}/2)=1$ is excluded unless $k=k_{0}=2$, in which case $p$-adic
Abel-Jacobi map at the cycle $\Delta_{0,0}$ is computed. This range of
classical weights plays here the role of the set of balanced weights in \cite%
{DR} .

Second, we observe that , in general, at any pair of classical weights $%
(r,s) $, the value $L_{p}(\mathbf{f},\mathbf{g})(r,s)=\frac{%
[g_{r},h_{r,s}]_{N}}{[g_{r},g_{r}]_{N}}$ belongs to $\mathbb{\overline{%
\mathbb{Q}}}$ , which we view in $\mathbb{C}_{p}$, via the $p-$adic
embedding fixed beforehand.

Third, for classical weights $r,s$ above $(k,k)$ and $k_{0}$ with $%
k-(k_{0}/2)\leq0$, it is 
\begin{equation*}
L_{p}(\mathbf{f},\mathbf{g})=\frac{[g_{r}^{(p)},h_{r,s}^{(p)}]_{N}}{%
[g_{r}^{(p)},g_{r}^{(s)}]_{N}}=\frac{1}{2}\frac{%
[g_{r}^{(p)},e_{g_{r}}e_{ord}j_{X^{\prime}(N)}^{\ast}(\theta^{%
\prime-1-t}f_{s}^{[\pi^{\prime}]}-\theta^{-1-t}f_{s}^{[\pi]})]_{Np}}{%
[g_{r}^{(p)},g_{r}^{(p)}]_{N}}
\end{equation*}

which is proved, as in the balanced case, to be equal to%
\begin{align*}
& \frac{1}{2}\frac{(-1)^{t}}{t!}\frac{\mathcal{E}_{0}(f_{s},g_{r})}{\mathcal{%
E}_{1}(g_{r})}\frac{[g_{r}^{(p)},e_{g_{r}}j_{X^{\prime}(N)}^{\ast
}(\theta^{\prime-1-t})f_{s}^{[\pi^{\prime}]}-\theta^{-1-t}f_{s}^{[%
\pi]})]_{Np}}{[g_{r}^{(p)},g_{r}^{(p)}]_{N}} \\
& =\frac{1}{2}\frac{(-1)^{t}}{t!}\frac{\mathcal{E}_{0}(f_{s},g_{r})}{%
\mathcal{E}_{1}(g_{r})\mathcal{E}_{1}(g_{r})}\frac{[g_{r}^{(p)},\text{ }%
j_{X^{\prime}(N)}^{\ast}(\theta^{\prime-1-t}f_{r}^{[\pi^{\prime}]}-%
\theta^{-1-t}f_{s}^{[\pi]})]_{N}}{[g_{r}^{(p)},g_{r}^{(p)}]_{N}}
\end{align*}

By theorem \ 1.1. of \cite{Michino}, this is a special value of the complex $%
L-$function of the representation associated to $g_{r}$ tensored by the
representation associated to $\theta^{\prime-1-t}f_{s}-\theta^{-1-t}f_{s}$
restricted to $X(N)$. This corresponds with the second of the\ three cases
dealt with in the quoted theorem. Our ongoing reserch tackles the
construction of an Euler system out of $L_{p}(\mathbf{f}$,$\mathbf{\mathbf{g}%
}),$ in the same spirit as in \cite{DR II}.

\end{document}